\documentclass[12pt, a4paper]{amsart}
\usepackage{amsmath,amsthm,amsfonts,amscd,amssymb,eucal,latexsym,mathrsfs}
\usepackage[all]{xy}
\usepackage{latexsym, amssymb, amscd, mathrsfs, graphics, fullpage}
\usepackage[backref=page,colorlinks,bookmarksopen=true]{hyperref}

\def\gG{\mathfrak{G}}
\def\gH{\mathfrak{H}}

\def\gM{\mathfrak{M}}
\def\gN{\mathfrak{N}}

\begin{document}

\title[Measured quantum groupoids with a central basis] {Measured quantum groupoids with a central basis}
\author{Michel Enock}
\address{Institut de Math\'ematiques de Jussieu, Unit\'{e} Mixte Paris 6 / Paris 7 /
CNRS de Recherche 7586 \\175, rue du Chevaleret, Plateau 7E, F-75013 Paris}
 \email{enock@math.jussieu.fr}

\begin{abstract}
Mimicking the von Neumann version of Kustermans and Vaes' locally compact quantum groups, Franck Lesieur had introduced a notion of measured quantum groupoid, in the setting of von Neumann algebras. In this article, we suppose that the basis of the measured quantum groupoid is central; in that case, we prove that a specific sub-${\bf C}^*$ algebra is invariant under all the data of the measured quantum groupoid; moreover, this sub-${\bf C}^*$-algebra is a continuous field of ${\bf C}^*$-algebras; when the basis is central in both the measured quantum groupoid and its dual, we get that the measured quantum groupoid is a continuous field of locally compact quantum groups. On the other hand, using 
this sub-${\bf C}^*$-algebra, we prove that any abelian measured quantum groupoid comes from a locally compact groupoid.  \end{abstract}

\keywords
{Measured quantum groupoids, continuous fields of ${\bf C}^*$-algebras, locally compact quantum groups}
\maketitle
\newpage
\section{Introduction}
\label{intro}
\subsection{}
 In two articles (\cite{Val2}, \cite{Val3}), J.-M. Vallin has introduced two notions (pseudo-multiplicative
unitary, Hopf-bimodule), in order to generalize, up to the groupoid
case, the classical notions of multiplicative unitary \cite{BS}, and of Hopf-von Neumann algebras \cite{ES}, which were introduced to describe and explain duality of groups, and leaded to appropriate notions
of quantum groups (\cite{ES}, \cite{W1}, \cite{W2}, \cite{BS}, \cite{MN}, \cite{W3}, \cite{KV1}, \cite{KV2}, \cite{MNW}). 
\\ In another article \cite{EV}, J.-M. Vallin and the author have constructed, from a depth 2 inclusion of
von Neumann algebras $M_0\subset M_1$, with an operator-valued weight $T_1$ verifying a regularity
condition, a pseudo-multiplicative unitary, which leaded to two structures of Hopf bimodules, dual
to each other. Moreover, we have then
construct an action of one of these structures on the algebra $M_1$ such that $M_0$
is the fixed point subalgebra, the algebra $M_2$ given by the basic construction being then
isomorphic to the crossed-product. We construct on $M_2$ an action of the other structure, which
can be considered as the dual action.
\\  If the inclusion
$M_0\subset M_1$ is irreducible, we recovered quantum groups, as proved and studied in former papers
(\cite{EN}, \cite{E1}).
\\ Therefore, this construction leads to a notion of "quantum groupoid", and a construction of a
duality within "quantum groupoids". 
\subsection{}
In a finite-dimensional setting, this construction can be
mostly simplified, and is studied in \cite{NV1}, \cite{BSz1},
\cite{BSz2}, \cite{Sz}, \cite{Val4}, \cite{Val5}, and examples are described. In \cite{NV2}, the link between these "finite quantum
groupoids" and depth 2 inclusions of
$II_1$ factors is given. 
\subsection{}
Franck Lesieur introduced (\cite{L1}) a notion of "measured quantum groupoids", in which a modular hypothesis on the basis is required. Mimicking in a wider setting the technics of Kustermans and Vaes \cite{KV1}, he obtained then a pseudo-multiplicative unitary, which, as in the quantum group case, "contains" all the information of the object (the von Neuman algebra, the coproduct, the antipode, the co-inverse). Unfortunately, the axioms chosen by Lesieur don't fit perfectely with the duality (namely, the dual object does not fit the modular condition on the basis chosen in \cite{L1}), and, in order to get a perfect symmetry, Lesieur gave the name of "measured quantum groupoid" to a wider class (\cite{L2}). In \cite{E3} had been shown that, with suitable conditions, the objects constructed in \cite{EV} from depth 2 inclusions, were "measured quantum groupoids" in this new sense. The axioms given in \cite{L2} were very complicated, and there was a serious need for simplification. This was made in \cite{E5}, and recalled in \cite{E6} in an appendix. 
\subsection{}
All these constructions had been made in a von Neumann setting, which was natural, once we are dealing with (or thinking of) depth 2 inclusions of von Neumann algebras. But, as for quantum groups, a 
 ${\bf C}^*$-version of this theory is to be done, at least to obtain quantum objects similar to locally compact groupoids. Many difficulties exist on that direction : how to define a relative ${\bf C}^*$-tensor product ? how to define the analog of operator-valued weights at the ${\bf C}^*$ level ?
 \newline
 A first attempt in that direction is due to T. Timmermann who defined a relative ${\bf C}^*$-tensor product and ${\bf C}^*$-pseudo-multiplicative unitaries (\cite{Ti}). 
 \newline
 This article is another step in that direction, and is devoted only to the special case when the basis of the measured quantum groupoid is central. A first version had been given in \cite{E4}. 
 \newline
 In that case, we get closed links with the theory of continuous fields of ${\bf C}^*$-algebras, as studied by Etienne Blanchard; using this theory and formalism, we then obtain some results at the ${\bf C}^*$ level for measured quantum groupoids having a central basis. This will allow to prove that measured quantum groupoid whose underlying von Neumann algebra is commutative comes from a locally compact groupoid. Applying this result to the measured quantum groupoid obtained from a measured groupoid, we obtain Ramsay's theorem on measured groupoids (which says, roughly speaking, that a measured quantum groupoid is equivalent, with respect to the product, the inverse and the measure) to a locally compact groupoid. A similar result is obtained for measured fields of locally compact quantum groups. 
 \subsection{}
The paper is organized as follows : in chapter \ref{pr}, we give all the preliminaries needed for that theory, mostly Connes-Sauvageot relative tensor product, weights on ${\bf C}^*$-algebras and continuous fields of ${\bf C}^*$-algebras; in chapter \ref{quantum} is recalled the notion of pseudo-multiplicative unitary, and the Hopf-bimodules associated, and the notion of measured quantum groupoid. In chapter \ref{anw}, we construct a sub-${\bf C}^*$-algebra of a measured quantum groupoid, which is a  invariant by all the data of the measured quantum groupoid. 
 \newline
 In chapters \ref{central} and \ref{mcfield}, we deal with the particular case of a measured quantum groupoid whose basis is central; in that case, we obtain, in chapter \ref{central}, properties of the restrictions of the coproduct and the weights to this sub-${\bf C}^*$-algebra; in \ref{mcfield}, we prove that this sub-${\bf C}^*$-algebra is, in two different ways,  a continuous field of ${\bf C}^*$-algebras, and that the restriction of the coproduct sends this sub-${\bf C}^*$-algebra into the multiplier
algebra of the min tensor product of these continuous fields, as introduced by Blanchard in \cite{Bl1}. 
\newline
In particular, in chapter \ref{abelian}, we look after a measured quantum groupoid, whose underlying von Neuman algebra itself is abelian; it is then proved that we obtain, in that case, a locally compact groupoid. Applying that result to the abelian measured quantum groupoid constructed from a measured groupoid, we recover Ramsay's theorem. 
\newline
In chapter \ref{fieldqg}, we define a notion of a measured field of locally compact quantum groups, and use this construction to get that it is equivalent to a continuous one, in a way which is similar to Ramsay's theorem; all that was underlying in Blanchard's work \cite{Bl2}; these are exactly the measured quantum groupoids with central basis, and with a dual which has also a central basis. Blanchard's examples are recalled.  
\newline
We finish this article (chapter \ref{Example}) by giving De Commer's example (\cite{DC}) of a measure quantum groupoid with a central basis ${\bf C}^2$, which is not central in the dual. 
 \subsection{}
The author is mostly indebted to E. Blanchard, F. Lesieur, S. Vaes, 
L. Va\u{\i}nerman and J.-M. Vallin for many fruitful conversations.

\section{Preliminaries}
\label{pr}
In this chapter are mainly recalled definitions and notations about Connes' spatial
theory and the fiber product construction (\ref{spatial})
which are the main technical tools of the theory of measured quantum groupoids. In \ref{C*} are recalled classical results about weights on ${\bf C}^*$-algebras, and a standard procedure for going from ${\bf C}^*$-algebra weight theory to von Neumann weight theory and vice versa. In \ref{defCXC*} is recalled the definition of a continuous field of ${\bf C}^*$-algebras, and E. Blanchard's results on the minimal tensor product of two continuous fields of ${\bf C}^*$-algebras.

 
\subsection{Spatial theory and relative tensor products of Hilbert spaces (\cite{C1},\cite{S},\cite{T})}
\label{spatial}
 Let $N$ a von neumann algebra, $\nu$ a normal semi-finite faithful weight on $N$; we shall denote by $H_\nu$, $\gN_\nu$, etc the canonical objects of the Tomita-Takesaki theory associated to the weight $\nu$; let $\alpha$ be a non degenerate faithful representation of $N$ on a Hilbert space $\mathcal H$; the set of $\nu$-bounded elements of the left-module $_\alpha\mathcal H$ is :
\[D(_\alpha\mathcal{H}, \nu)= \lbrace \xi \in \mathcal{H};\exists C < \infty ,\| \alpha (y) \xi\|
\leq C \| \Lambda_{\nu}(y)\|,\forall y\in \gN_{\nu}\rbrace\]
Then, for any $\xi$ in $D(_\alpha\mathcal{H}, \nu)$, there exists a bounded operator
$R^{\alpha,\nu}(\xi)$ from $H_\nu$ to $\mathcal{H}$,  defined, for all $y$ in $\gN_\nu$ by :
\[R^{\alpha,\nu}(\xi)\Lambda_\nu (y) = \alpha (y)\xi\]
which intertwines the representions of $N$. 
\newline
If $\xi$, $\eta$ are bounded vectors, we define the operator product :
\[\langle\xi,\eta\rangle_{\alpha,\nu} = R^{\alpha,\nu}(\eta)^* R^{\alpha,\nu}(\xi)\]
which belongs to $\pi_{\nu}(N)'$, which, thanks to Tomita-Takesaki theory, will be identified to the opposite von Neumann algebra $N^o$, on which is defined a canonical weight $\nu^o$. 
\newline
If now $\beta$ is a non degenerate faithful antirepresentation of $N$ on a Hilbert space $\mathcal K$, we define the relative tensor product $\mathcal K\underset{\nu}{_\beta\otimes_\alpha}\mathcal H$ as the completion of the algebraic tensor product $\mathcal K\odot D(_\alpha\mathcal{H}, \nu)$ by the scalar product defined,  if $\xi_1$, $\xi_2$ are in $\mathcal{K}$, $\eta_1$, $\eta_2$ are in $D(_\alpha\mathcal{H},\nu)$, by the following formula :
\[(\xi_1\odot\eta_1 |\xi_2\odot\eta_2 )= (\beta(\langle\eta_1, \eta_2\rangle_{\alpha,\nu})\xi_1 |\xi_2)\]
If $\xi\in \mathcal{K}$, $\eta\in D(_\alpha\mathcal{H},\nu)$, we shall denote $\xi\underset{\nu}{_\beta\otimes_\alpha}\eta$ the image of $\xi\odot\eta$ into $\mathcal K\underset{\nu}{_\beta\otimes_\alpha}\mathcal H$, and, writing $\rho^{\beta, \alpha}_\eta(\xi)=\xi\underset{\nu}{_\beta\otimes_\alpha}\eta$, we get a bounded linear operator from $\mathcal H$ into $\mathcal K\underset{\nu}{_\beta\otimes_\alpha}\mathcal H$, which is equal to $1_\mathcal K\otimes_\nu R^{\alpha, \nu}(\eta)$. 
\newline
One should bear in mind that, if we start from another faithful semi-finite normal weight $\nu '$, we
get another Hilbert space $\mathcal{H}\underset{\nu'}{_\beta\otimes_\alpha}\mathcal{K}$; there exists an isomorphism $U^{\nu, \nu'}_{\beta, \alpha}$ from $\mathcal{H}\underset{\nu}{_\beta\otimes_\alpha}\mathcal{K}$ to $\mathcal{H}\underset{\nu'}{_\beta\otimes_\alpha}\mathcal{K}$, which is unique up to some functorial property (\cite{S}, 2.6) (but this isomorphism
does not send  $\xi\underset{\nu}{_\beta\otimes_\alpha}\eta$ on
$\xi\underset{\nu'}{_\beta\otimes_\alpha}\eta$ !). 
\newline
The relative tensor product $\mathcal K\underset{\nu}{_\beta\otimes_{\pi_\nu}}H_\nu$ is canonically identified with $\mathcal K$ (\cite{S},2.4(a)). 
\newline
The linear set generated by operators $\theta^{\alpha, \nu}(\xi, \eta)=R^{\alpha, \nu}(\xi)R^{\alpha, \nu}(\eta)^*$, for all $\xi$, $\eta$ in $D(_\alpha\mathcal H, \nu)$, is a weakly dense ideal in $\alpha(N)'$. We shall denote by $\mathcal K_{\alpha, \nu}$ the norm closure of this set of operators, which is a ${\bf C}^*$-algebra, and also a a weakly dense ideal of $\alpha(N)'$. 
\newline
Moreover, there exists a family $(e_i)_{i\in I}$ of vectors in $D(_\alpha\mathcal H, \nu)$ such that the operators $\theta^{\alpha, \nu}(e_i, e_i)$ are 2 by 2 orthogonal projections ($\theta^{\alpha, \nu}(e_i, e_i)$ being then the projection on the closure of $\alpha(N)e_i$). Such a family is called an orthogonal $(\alpha, \nu)$-basis of $\mathcal H$. 
\newline
We shall denote $\sigma_\nu$ the relative flip, which is a unitary sending $\mathcal{K}\underset{\nu}{_\beta\otimes_\alpha}\mathcal{H}$ onto $\mathcal{H}\underset{\nu^o}{_\alpha\otimes _\beta}\mathcal{K}$, defined, for any $\xi$ in $D(\mathcal {K}_\beta ,\nu^o )$, $\eta$ in $D(_\alpha \mathcal {H},\nu)$, by :
\[\sigma_\nu (\xi\underset{\nu}{_\beta\otimes_\alpha}\eta)=\eta\underset{\nu^o}{_\alpha\otimes_\beta}\xi\]
If $\xi\in D(\mathcal H_\beta, \nu^o)$ and $\eta\in\mathcal K$, we can then define a bounded linear operator $\lambda_\xi^{\beta, \alpha}$ from $\mathcal K$ into $\mathcal K\underset{\nu}{_\beta\otimes_\alpha}\mathcal H$ such that $\lambda_\xi^{\beta, \alpha}=\xi\underset{\nu}{_\beta\otimes_\alpha}\eta$. 
\newline
If $x\in \beta(N)'$, $y\in \alpha(N)'$, it is possible to define an operator $x\underset{\nu}{_\beta\otimes_\alpha}y$ on $\mathcal K\underset{\nu}{_\beta\otimes_\alpha}\mathcal H$, with natural values on the elementary tensors. It is easy to get that this operator does not depend upon the weight $\nu$ and it will be denoted $x\underset{N}{_\beta\otimes_\alpha}y$. Let $A$ be a ${\bf C}^*$-algebra of operators acting on $\mathcal H$, such that $A\subset \alpha(N)'$, and $B$ a ${\bf C}^*$-algebra of operators acting on $\mathcal K$, such that $B\subset \beta(N)'$; the linear space generated by  the set of operators $x\underset{N}{_\beta\otimes_\alpha}y$, with $x\in B$ and $y\in A$, is clearly an involutive algebra, and its norm closure a ${\bf C}^*$-algebra, that we shall denote by $B\underset{N}{_\beta\otimes_\alpha}A$ . 
\newline
Let us suppose now that $\mathcal{H}$ is a $N-N_1$ bimodule; that means that there exists a von
Neumann algebra $N_1$, and a non-degenerate normal anti-representation $\epsilon$ of $N_1$ on
$\mathcal{H}$, such that
$\epsilon (N_1)\subset\alpha (N_1)'$. We shall write then $_\alpha\mathcal{H}_\epsilon$. If $y$ is in $N_1$, we
have seen that it is possible to define then the operator
$1_{\mathcal{K}}\underset{N}{_\beta\otimes_\alpha}\epsilon (y)$ on
$\mathcal{K}\underset{\nu}{_\beta\otimes_\alpha}\mathcal{H}$, and we define this way a
non-degenerate normal antirepresentation of $N_1$ on
$\mathcal{K}\underset{\nu}{_\beta\otimes_\alpha}\mathcal{H}$, we shall call again $\epsilon$ for
simplification. If $\mathcal K$ is a $N_2-N$ bimodule, then
$\mathcal{K}\underset{\nu}{_\beta\otimes_\alpha}\mathcal{H}$ becomes a $N_2-N_1$ bimodule (Connes'
fusion of bimodules).
\newline
Taking a faithful semi-finite normal weight
$\nu_1$  on $N_1$, and a left $N_1$-module $_{\gamma}\mathcal{L}$ (i.e. a Hilbert space $\mathcal{L}$ and a normal
non-degenerate representation $\gamma$ of $N_1$ on $\mathcal{L}$), it is possible then to define
$(\mathcal{K}\underset{\nu}{_\beta\otimes_\alpha}\mathcal{H})\underset{\nu_1}{_\epsilon\otimes_\gamma}\mathcal{L}$.
Of course, it is possible also to consider the Hilbert space
$\mathcal{K}\underset{\nu}{_\beta\otimes_\alpha}(\mathcal{H}\underset{\nu_1}{_\epsilon\otimes_\gamma}\mathcal{L})$.
It can be shown that these two Hilbert spaces are isomorphic as $\beta (N)'-\gamma
(N_1)^{'o}$-bimodules. (In (\cite{Val2} 2.1.3), the proof, given for $N=N_1$ abelian can be used, without
modification, in that wider hypothesis). We shall write then
$\mathcal{K}\underset{\nu}{_\beta\otimes_\alpha}\mathcal{H}\underset{\nu_1}{_\epsilon\otimes_\gamma}\mathcal{L}$
without parenthesis, to emphazise this coassociativity property of the relative tensor
product.
\newline
Dealing now with that Hilbert space $\mathcal{K}\underset{\nu}{_\beta\otimes_\alpha}\mathcal{H}\underset{\nu_1}{_\epsilon\otimes_\gamma}\mathcal{L}$, there exist different flips, and it is necessary to be careful with notations. For instance, $1\underset{\nu}{_\beta\otimes_\alpha}\sigma_{\nu_1}$ is the flip from this Hilbert space onto $\mathcal{K}\underset{\nu}{_\beta\otimes_\alpha}(\mathcal L\underset{\nu_1^o}{_\gamma\otimes_\epsilon}\mathcal H)$, where $\alpha$ is here acting on the second leg of $\mathcal L\underset{\nu^o}{_\gamma\otimes_\epsilon}\mathcal H$ (and should therefore be written $1\underset{\nu^o}{_\gamma\otimes_\epsilon}\alpha$, but this will not be done for obvious reasons). Here, the parenthesis remains, because there is no associativity rule, and to remind that $\alpha$ is not acting on $\mathcal L$. The adjoint of $1\underset{\nu}{_\beta\otimes_\alpha}\sigma_{\nu_1}$ is $1\underset{\nu}{_\beta\otimes_\alpha}\sigma_{\nu_1^o}$. 
\newline
The same way, we can consider $\sigma_\nu\underset{\nu_1}{_\epsilon\otimes_\gamma}1$ from $\mathcal{K}\underset{\nu}{_\beta\otimes_\alpha}\mathcal{H}\underset{\nu_1}{_\epsilon\otimes_\gamma}\mathcal{L}$ onto $(\mathcal H\underset{\nu^o}{_\alpha\otimes_\beta}\mathcal K)\underset{\nu_1}{_\epsilon\otimes_\gamma}\mathcal L$. 
 \newline
 Another kind of flip sends $\mathcal{K}\underset{\nu}{_\beta\otimes_\alpha}(\mathcal L\underset{\nu_1^o}{_\gamma\otimes_\epsilon}\mathcal H)$ onto $\mathcal L\underset{\nu_1^o}{_\gamma\otimes_\epsilon}(\mathcal{K}\underset{\nu}{_\beta\otimes_\alpha}\mathcal{H})$. We shall denote this application $\sigma^{1,2}_{\alpha, \epsilon}$ (and its adjoint $\sigma^{1,2}_{\epsilon, \alpha}$), in order to emphasize that we are exchanging the first and the second leg, and the representations $\alpha$ and $\epsilon$ on the third leg. 

\subsection{Operator-valued weights}
\label{ovw}
Let $M_0\subset M_1$ be an inclusion of von Neumann algebras  (for simplification, these algebras will be supposed to be $\sigma$-finite), equipped with a normal faithful semi-finite operator-valued weight $T_1$ from $M_1$ to $M_0$ (to be more precise, from $M_1^{+}$ to the extended positive elements of $M_0$ (cf. \cite{T} IX.4.12)). Let $\psi_0$ be a normal faithful semi-finite weight on $M_0$, and $\psi_1=\psi_0\circ T_1$; for $i=0,1$, let $H_i=H_{\psi_i}$, $J_i=J_{\psi_i}$, $\Delta_i=\Delta_{\psi_i}$ be the usual objects constructed by the Tomita-Takesaki theory associated to these weights. 
 \newline
Following (\cite{EN} 10.6), for $x$ in $\gN_{T_1}$, we shall define $\Lambda_{T_1}(x)$ by the following formula, for all $z$ in $\gN_{\psi_{0}}$ :
\[\Lambda_{T_1}(x)\Lambda_{\psi_{0}}(z)=\Lambda_{\psi_1}(xz)\]
This operator belongs to $Hom_{M_{0}^o}(H_{0}, H_1)$; if $x$, $y$ belong to $\gN_{T_1}$, then $\Lambda_{T_1}(x)\Lambda_{T_1}(y)^*$ belongs to the von Neumann algebra $M_{2}=J_1M'_0J_1$, which is called the basic construction made from the inclusion $M_0\subset M_1$, and
$\Lambda_{T_1}(x)^*\Lambda_{T_1}(y)=T_1(x^*y)\in M_0$. 
\newline
By Tomita-Takesaki theory, the Hilbert space $H_1$ bears a natural structure of $M_1-M_1^o$-bimodule, and, therefore, by restriction, of $M_0-M_0^o$-bimodule. Let us write $r$ for the canonical representation of $M_0$ on $H_1$, and $s$ for the canonical antirepresentation given, for all $x$ in $M_0$, by $s(x)=J_1r(x)^*J_1$. Let us have now a closer look to the subspaces $D(H_{1s}, \psi_0^o)$ and $D(_rH_1, \psi_0)$. If $x$ belongs to $\gN_{T_1}\cap\gN_{\psi_1}$, we easily get that $J_1\Lambda_{\psi_1}(x)$ belongs to $D(_rH_1, \psi_0)$, with :
\[R^{r, \psi_0}(J_1\Lambda_{\psi_1}(x))=J_1\Lambda_{T_1}(x)J_0\]
and $\Lambda_{\psi_1}(x)$ belongs to $D(H_{1s}, \psi_0)$, with :
\[R^{s, \psi_0^o}(\Lambda_{\psi_1}(x))=\Lambda_{T_1}(x)\]
The subspace $D(H_{1s}, \psi_0^o)\cap D(_rH_1, \psi_0)$ is dense in $H_1$; more precisely, let $\mathcal T_{\psi_1, T_1}$ be the algebra made of elements $x$ in $\gN_{\psi_1}\cap\gN_{T_1}\cap\gN_{\psi_1}^*\cap\gN_{T_1}^*$, analytical with respect to $\psi_1$, and such that, for all $z$ in ${\bf C}$, $\sigma^{\psi_1}_z(x_n)$ belongs to $\gN_{\psi_1}\cap\gN_{T_1}\cap\gN_{\psi_1}^*\cap\gN_{T_1}^*$. Then (\cite{E6}, 2.2.1):
\newline
(i) the algebra $\mathcal T_{\psi_1, T_1}$ is weakly dense in $M_1$; it will be called Tomita's algebra with respect to $\psi_1$ and $T_1$; 
\newline
(ii) for any $x$ in  $\mathcal T_{\psi_1, T_1}$, $\Lambda_{\psi_1}(x)$ belongs to $D(H_{1s}, \psi_0)\cap D(_rH_1, \psi_0)$;
\newline
(iii) for any $\xi$ in $D(H_{1s}, \psi_0^o))$, there exists a sequence $x_n$ in $\mathcal T_{\psi_1, T_1}$ such that $\Lambda_{T_1}(x_n)=R^{s, \psi_0^o}(\Lambda_{\psi_1}(x_n))$ is weakly converging to $R^{s, \psi_0^o}(\xi)$ and $\Lambda_{\psi_1}(x_n)$ is converging to $\xi$. 
\newline
More precisely, in (\cite{E3}, 2.3) was constructed an increasing sequence of projections $p_n$ in $M_1$, converging to $1$, and elements $x_n$ in $\mathcal T_{\psi_1, T_1}$ such that $\Lambda_{\psi_1}(x_n)=p_n\xi$. We then get that :
\begin{align*}
T_1(x_n^*x_n)
&=
\langle R^{s, \psi_0^o}(\Lambda_{\psi_1}(x_n), R^{s, \psi_0^o}(\Lambda_{\psi_1}(x_n)\rangle_{s, \psi_0^o}\\
&=
\langle p_n\xi, p_n\xi\rangle_{s, \psi_0^o}\\
&=
R^{s, \psi_0^o}(\xi)^*p_nR^{s, \psi_0^o}(\xi)
\end{align*}
which is increasing and weakly converging to $\langle\xi, \xi\rangle_{s, \psi_0^o}$. 
Moreover, if $M_0$ is abelian, and if we write $X$ for the spectrum of the ${\bf C}^*$-algebra generated by all elements of the form $<\eta_1, \eta_2>_{s, \psi_0^o}$, we can identify $\psi_0$ as a positive Radon measure on $X$, and $M_0$ with $L^\infty(X, \psi_0)$; using now Dini's theorem on $C_0(X)$, we get that $T_1(x_n^*x_n)$ is norm converging to $\langle\xi, \xi\rangle_{s, \psi_0^o}$, and that :
\[\|\Lambda_{T_1}(x_n)-R^{s, \psi_0^o}(\xi)\|^2=\|T_1(x_n^*x_n)-\langle\xi, \xi\rangle_{s, \psi_0^o}\|\]
is converging to $0$.

\subsection{Fiber products of von Neumann algebras and slice maps (\cite{EV},\cite{E2}).}
\label{fiber}
Let's go on with the notations of \ref{spatial}.  If $P$ is a von Neumann algebra on $\mathcal H$, with $\alpha(N)\subset P$, and $Q$ a von Neumann algebra on $\mathcal K$, with $\beta(N)\subset Q$, then we define the fiber product $Q\underset{N}{_\beta*_\alpha}P$ as $\{x\underset{N}{_\beta\otimes_\alpha}y, x\in Q', y\in P'\}'$. 
\newline
Moreover, this von Neumann algebra can be defined independently of the Hilbert spaces on which $P$ and $Q$ are represented; if $(i=1,2)$, $\alpha_i$ is a faithful non degenerate homomorphism from $N$ into $P_i$, $\beta_i$ is a faithful non degenerate antihomomorphism from $N$ into $Q_i$, and $\Phi$ (resp. $\Psi$) an homomorphism from $P_1$ to $P_2$ (resp. from $Q_1$ to $Q_2$) such that $\Phi\circ\alpha_1=\alpha_2$ (resp. $\Psi\circ\beta_1=\beta_2$), then, it is possible to define an homomorphism $\Psi\underset{N}{_{\beta_1}*_{\alpha_1}}\Phi$ from $Q_1\underset{N}{_{\beta_1}*_{\alpha_1}}P_1$ into $Q_2\underset{N}{_{\beta_2}*_{\alpha_2}}P_2$. 
\newline
Let $A$ be in $Q\underset{N}{_\beta *_\gamma}P$, and let $\xi_1$, $\xi_2$ be in
$D(\mathcal{H}_\beta,\nu^o)$; let us define $(\omega_{\xi_1, \xi_2}\underset{\nu}{_\beta*_\gamma}id)(A)$ as a bounded operator on $\mathcal{K}$,
which belongs to $P$, such that :
\[((\omega_{\xi_1, \xi_2}\underset{\nu}{_\beta*_\gamma}id)(A)\eta_1|\eta_2)=
(A(\xi_1\underset{\nu}{_\beta\otimes_\gamma}\eta_1)|
\xi_2\underset{\nu}{_\beta\otimes_\gamma}\eta_2)\]
One should note that $(\omega_{\xi_1, \xi_2}\underset{\nu}{_\beta*_\gamma}id)(1)=\gamma (\langle\xi_1, \xi_2 \rangle_{\beta, \nu^o})$. 
\newline
Let us define the same way, for any $\eta_1$, $\eta_2$ in
$D(_\gamma\mathcal{K}, \nu)$:
\[(id\underset{\nu}{_\beta*_\gamma}\omega_{\eta_1, \eta_2})(A)=(\rho^{\beta, \gamma}_{\eta_2})^*A\rho^{\beta, \gamma}_{\eta_1}\]
which belongs to $Q$. 
\newline
Let $\phi$ be a normal semi-finite weight on
$Q^+$; we may define an element of the extended positive part
of $P$, denoted
$(\phi\underset{\nu}{_\beta*_\gamma}id)(A)$, such that, for all $\eta$ in $D(_\gamma \mathcal K, \nu)$, we have :
\[\|(\phi\underset{\nu}{_\beta*_\gamma}id)(A)^{1/2}\eta\|^2=\phi(id\underset{\nu}{_\beta*_\gamma}\omega_\eta)(A)\]
Moreover, if $\psi$ is a normal semi-finite weight on $P^+$, we have then :
\[\psi(\phi\underset{\nu}{_\beta*_\gamma}id)(A)=\phi(id\underset{\nu}{_\beta*_\gamma}\psi)(A)\]
and if $\omega_i$ be in $Q_{*}$ such that $\phi_1=sup_i\omega_i$, we
have $(\phi_1\underset{\nu}{_\beta*_\gamma}id)(A)=sup_i(\omega_i\underset{\nu}{_\beta*_\gamma}id)(A)$.
\newline
Let now $Q_1$ be a von Neuman algebra such that $\beta(N)\subset Q_1\subset Q$, 
and $P_1$ be a von Neuman algebra such that $\gamma(N)\subset P_1\subset P$
and let $T$ (resp. $T'$) 
be a normal faithful semi-finite operator valued weight from $Q$ to $Q_1$ (resp. from $P$ to $P_1$); then, there exists an element $(T\underset{\nu}{_\beta*_\gamma}id)(A)$
of the extended positive part
of $Q_1\underset{N}{_\beta*_\gamma}P$, such that (\cite{E2}, 3.5), for all
$\eta$ in $D(_\gamma \mathcal K, \nu)$, and $\xi$ in $\mathcal H$, we have :
\[\|(T\underset{\nu}{_\beta*_\gamma}id)(A)^{1/2}(\xi\underset{\nu}{_\beta\otimes_\gamma}\eta)\|^2=
\|T[(id\underset{\nu}{_\beta*_\gamma}\omega_\eta)(A)]^{1/2}\xi\|^2\]
If $\phi_1$ is a normal semi-finite weight on $P_1$, we have :
\[(\phi_1\circ T\underset{\nu}{_\beta*_\gamma}id)(A)=(\phi_1\underset{\nu}{_\beta*_\gamma}id)
(T\underset{\nu}{_\beta*_\gamma}id)(A)\]
We define the same way an element $(id\underset{\nu}{_\beta*_\gamma}T')(A)$ of the extended positive part
of
$Q\underset{N}{_\gamma*_\beta}P_1$, and we have :
\[(id\underset{\nu}{_\beta*_\gamma}T')((T\underset{\nu}{_\beta*_\gamma}id)(A))=
(T\underset{\nu}{_\beta*_\gamma}id)((id\underset{\nu}{_\beta*_\gamma}T')(A))\]
 Considering now an element $x$ of $Q\underset{\nu}{_\beta*_{\pi_\nu}} \pi_\nu(N)$, which can be identified 
 to $Q\cap\beta(N)'$ (thanks to the identification of $\mathcal K\underset{\nu}{_\beta\otimes_{\pi_\nu}}H_\nu$ with $\mathcal K$), we get that, for $e$ in $\gN_\nu$, we have  \[(id_\beta\underset{\nu}{*}{}_{\pi_\nu}\omega_{J_\nu \Lambda_{\nu}(e)})(x)=\beta(ee^*)x\]
 Therefore, by increasing limits, we get that $(id_\beta\underset{\nu}{*}{}_{\pi_\nu}\nu)$ is the injection of $Q\cap\beta(N)'$ into $Q$.  More precisely, if $x$ belongs to $Q\cap\beta(N)'$, we have :
 \[(id\underset{\nu}{_\beta*_{\pi\nu}}\nu)(x{}_\beta\underset{\nu}{\otimes}{}_{\pi_\nu} 1)=x\]
 Therefore, if $T'$ is a normal faithful semi-finite operator-valued weight from $P$ onto $\gamma(N)$, we get that we have :
 \[(id_\beta\underset{\nu}{*}{}_\gamma\nu\circ T')(A){}_\beta\underset{\nu}{\otimes}{}_\gamma 1=
 (id_\beta\underset{\nu}{*}{}_\gamma T')(A)\]
If $\alpha(N)\subset Z(P)$, and $\beta(N)\subset Z(Q)$, the von Neumann algebra $Q\underset{N}{_\beta*_\alpha}P$ is clearly the weak closure of the ${\bf C}^*$-algebra $Q\underset{N}{_\beta\otimes_\alpha}P$ we have defined in \ref{spatial}. 

\subsection{Notations} 
\label{aut}
Let $M$ be a von Neuman algebra, and $\alpha$ an action from a locally compact group $G$ on $M$, i.e. a homomorphism from $G$ into $Aut M$, such that, for all $x\in M$, the function $g\mapsto \alpha_g(x)$ is $\sigma$-weakly continuous. Let us denote by ${\bf C}^*(\alpha)$ the set of elements $x$ of $M$, such that this function $t\mapsto \alpha_g(x)$ is norm continuous. It is (\cite{Pe}, 7.5.1) a sub-${\bf C}^*$-algebra of $M$, invariant under the $\alpha_g$, generated by the elements ($x\in N$,$f\in L^1(G)$):
\[\alpha_f(x)=\int_{\mathbb{R}}f(s)\alpha_s(x)ds\]
More precisely, we get that, for any $x$ in $M$, $\alpha_f(x)$ is $\sigma$-weakly converging to $x$ when $f$ goes in an approximate unit of $L^1(G)$, which proves that ${\bf C}^*(\alpha)$ is $\sigma$-weakly dense in $M$, and that $x\in M$ belongs to ${\bf C}^*(\alpha)$ if and only if this file is norm converging. 
\newline
If $\alpha_t$ and $\gamma_s$ are two one-parameter automorphism groups of $M$, such that, for all $s$, $t$ in $\mathbb{R}$, we have $\alpha_t\circ\gamma_s=\gamma_s\circ\alpha_t$, by considering the action of $\mathbb{R}^2$ given by $(s,t)\mapsto \gamma_s\circ\alpha_t$, we obtain a dense sub-${\bf C}^*$-algebra of $M$, on which both $\alpha$ and $\gamma$ are norm continuous, we shall denote ${\bf C}^*(\alpha, \gamma)$.

\subsection{Weights on ${\bf C}^*$-algebras}
\label{C*}
Let $A$ be a ${\bf C}^*$-algebra, and $\varphi$ a lower semi-continuous, densely defined non zero weight on $A$ (\cite{Co}). We shall use all classical notations, and, in particular, we shall denote $(H_\varphi, \Lambda_\varphi, \pi_\varphi)$ the GNS construction for $\varphi$; if $\varphi$ is faithful, so is $\pi_\varphi$; let us denote $M=\pi_\varphi(A)''$ and $\overline{\varphi}$ the semi-finite normal weight on $M^+$, constructed by  (\cite{B}, cor. 9), which verify $\overline{\varphi}\circ\pi_\varphi = \varphi$. 
\newline
Let us recall that if the ${\bf C}^*$-algebra $A$ is unital, any densely defined weight $\varphi$ is everywhere defined, and therefore finite. 
\newline
Following \cite{Co}, we shall say that $\varphi$ is KMS if there exists a norm-continuous one parameter group of automorphisms $\sigma_t$ of $A$ such that, for all $t\in \mathbb{R}$, $\varphi=\varphi\circ\sigma_t$, and such that $\varphi$ verifies the KMS conditions with respect to $\sigma$. (For an equivalent definition of these conditions, see \cite{KV1}, 1.3). One can find in (\cite{KV1}, 1.35) the proof that every KMS weight extends to a faithful extension $\overline{\varphi}$ on $M^+$, and that we have then $\pi_\varphi\circ\sigma_t=\sigma^{\overline{\varphi}}_t\circ\pi_\varphi$, where $\sigma^{\overline{\varphi}}_t$ is the modular automorphism group of $M$ given by the Tomita-Takesaki theory of the faithful semi-finite weight $\overline{\varphi}$ on $M$. This leads easily to the uniqueness of the one-parameter group $\sigma_t$, which we shall emphasize by writing it $\sigma_t^\varphi$.
\newline
Moreover, it is well known that the set of elements $x$ in $A$ such that the function $t\mapsto\sigma_t^\varphi(x)$ extends to an analytic function in $A$ is a dense involutive subalgebra of $A$ (see for instance \cite{Val1} 0.3.2 and 0.3.4). 
\newline

Let now $\psi$ be a normal semi-finite faithful weight on a von Neumann algebra $N$. We shall write ${\bf C}^*(\psi)$ the sub-${\bf C}^*$-algebra of ${\bf C}^*(\sigma^\psi)$ generated by elements $\sigma^\psi_f(x)$, with $f\in L^1(\mathbb{R})$ and $x\in\gM_\psi$. The weak closure of ${\bf C}^*(\psi)$ contains $\gM_\psi$, and, therefore, ${\bf C}^*(\psi)$ is weakly dense in $N$; moreover, it is straightforward to see that the restriction of $\psi$ to that ${\bf C}^*$-algebra is densely defined, lower semi-continuous and KMS. If $1\in{\bf C}^*(\psi)$, then the restriction of $\psi$ to  ${\bf C}^*(\psi)$ is finite, so $\psi(1)<\infty$ and ${\bf C}^*(\psi)={\bf C}^*(\sigma^\psi)$. If $\psi$ is a trace, then ${\bf C}^*(\psi)$ is the norm closure of $\gM_\psi$, and $M({\bf C}^*(\psi))=N$. 
\newline
If $\gamma_t$ is one-parameter group of $N$, such that $\psi\circ \gamma_t=\psi$, for all $t\in \mathbb{R}$, we may as well define the ${\bf C}^*$-algebra ${\bf C}^*(\psi, \gamma)$ generated by all elements :
\[\int_{\mathbb{R}^2}f(s)g(t)\sigma_s^\psi\circ\gamma_t(x)dsdt\]
where $f,g$ belong to $L^1(\mathbb{R})$, and $x$ belongs to $\gM_\psi$; this ${\bf C}^*$-algebra ${\bf C}^*(\psi, \gamma)$ is weakly dense in $N$, invariant under $\gamma$, the restriction of $\psi$ to this ${\bf C}^*$-algebra is densely defined, lower semi-continuous and KMS, and the restriction of $\gamma$ to this ${\bf C}^*$-algebra is norm-continuous. If $\psi$ is a trace, we have $M({\bf C}^*(\psi, \gamma))={\bf C}^*(\gamma)$.

\subsection{Continuous fields of ${\bf C}^*$-algebras}
\label{defCXC*}
Let $X$ be a locally compact space; following \cite{Ka}, we shall say that a ${\bf C}^*$-algebra $A$ is a $C_0(X)-{\bf C}^*$-algebra if there exists an injective non degenerate $*$-homomorphism $\alpha$ from $C_0(X)$ into $Z(M(A))$. If $x\in X$, let us write $C_x(X)$ for the ideal of $C_0(X)$ made of all functions in $C_0(X)$ with value $0$ at $x$, and let us consider $\alpha(C_x(X))A$, which is an ideal in $A$; let us write $A^x$ for the quotient ${\bf C}^*$-algebra $A/\alpha(C_x(X))A$. For any $a$ in $A$, let us write $a^x$ for its image in $A^x$. Then, we have (\cite{Bl2}, 2.8) :
\[\|a\|=sup_{x\in X}\|a^x\|\]
By definition (\cite{D}), we shall say that $A$ is a continuous field over $X$ if, for all $a$ in $A$, the function $x\mapsto \| a^x\|$ is continuous. 
\newline
Let $A$ be a $C_0(X)-{\bf C}^*$-algebra, and $\mathcal E$ be a $C_0(X)$-Hilbert module; let $\pi$ be a $C_0(X)$-linear morphism $\pi$ from $A$ to $\mathcal L(\mathcal E)$, which means that the specialization $\pi_x$ is a representation of $A$ on the Hilbert space $\mathcal E_x$ whose kernel contains $\alpha(C_x(X))A$. We say that $\pi$ is a continuous field of faithful representations if, for all $x\in X$, we have $Ker \pi_x=\alpha(C_x(X))A$. We may then, for all $x\in X$,  consider $\pi_x$ as a faithful representation of the ${\bf C}^*$-algebra $A^x$ on $\mathcal E_x$. It is proved in (\cite{Bl2}, 3.3) that, if $A$ is a separable $C_0(X)-{\bf C}^*$-algebra, the following are equivalent :
\newline
(i) $A$ is a continuous field over $X$ of ${\bf C}^*$-algebras;
\newline
(ii) there exists a continuous field of faithful representations of $A$. 
\newline
A continuous field of states on $A$ is a positive $C_0(X)$-linear application $\omega$ from $A$ into $C_0(X)$, such that, for all $x\in X$, the specialization $\omega_x$ is a state on $A^x$\vspace{5mm}. 
\newline
Given two $C_0(X)$-algebras $A_1$ and $A_2$, Blanchard (\cite{Bl1}, 2.9) had defined the minimal ${\bf C}^*$-norm on the involutive algebra $(A_1\odot A_2)/J(A_1, A_2)$, where $A_1\odot A_2$ is the algebraic tensor product of the algebras $A_1$ and $A_2$, and $J(A_1, A_2)$ the involutive ideal in $A_1\odot A_2$ made of finite sums $\sum_{i=1}^n a_i\otimes b_i$, with $a_i\in A_1$, $b_i\in A_2$, such that $\sum_{i=1}^n a_i^x\otimes b_i^x=0$, for all $x\in X$. 
This ${\bf C}^*$-norm is given by :
\[\|\sum_{i=1}^n a_i\otimes b_i\|_m=sup_{x\in X}\|\sum_{i=1}^n a_i^x\otimes b_i^x\|_m\]
where, on the right hand on the formula, is taken the minimal tensor product of the ${\bf C}^*$-algebras $A_1^x$ and $A_2^x$. The completion with respect to that norm will be called the minimal tensor product of the $C_0(X)$-algebras $A_1$ and $A_2$, and will be denoted $A_1\otimes_{C_0(X)}^m A_2$. 
\newline
In the case of continuous fields of ${\bf C}^*$-algebras, it is proved in (\cite{Bl2}, 3.21) that this ${\bf C}^*$-algebra, equipped with  the  morphism $f\mapsto f\otimes_{C_0(X)}^m 1=1\otimes_{C_0(X)}^m f$ is equal to the $C_0(X)-{\bf C}^*$-algebra $A_1\otimes_{C_0(X)}A_2$. (If $\pi_1$ (resp. $\pi_2$) is a faithful non degenerate $C_0(X)$-representation on a $C_0(X)$-Hilbert module $\mathcal E_1$ (resp. $\mathcal E_2$), $A_1\otimes_{C_0(X)}A_2$ is defined as operators on $\mathcal E_1\otimes_{C_0(X)}\mathcal E_2$ (\cite{Bl2}, 3.18) and does not depend on the choice of the $C_0(X)$-representations (\cite{Bl2}, 3.20).  Therefore, the tensor product $\otimes_{C_0(X)}^m$ is then associative (\cite{Bl1}, 4.1).

\section{Measured quantum groupoids}
\label{quantum}
In this section, we give a summary of the theory of "Hopf-bimodules" (\ref{Hbimod}), "pseudo-multiplicative unitaries" (\ref{defmult}), and  "measured quantum groupoids" (\cite{L1}, \cite{L2}, \cite{E5}, \cite{E6}) (\ref{LW}, \ref{thL1}, \ref{defMQG}, \ref{thL2}). We describe the canonical example of measured groupoids (\ref{gd}, \ref{gd2}). Proposition \ref{prop2Gamma} will be used in Theorem \ref{thcentral1}. 

\subsection{Definition}
\label{Hbimod}
A quintuplet $(N, M, \alpha, \beta, \Gamma)$ will be called a Hopf-bimodule, following (\cite{Val2}, \cite{EV} 6.5), if
$N$,
$M$ are von Neumann algebras, $\alpha$ a faithful non-degenerate representation of $N$ into $M$, $\beta$ a
faithful non-degenerate anti-representation of
$N$ into $M$, with commuting ranges, and $\Gamma$ an injective involutive homomorphism from $M$
into
$M\underset{N}{_\beta *_\alpha}M$ such that, for all $X$ in $N$ :
\newline
(i) $\Gamma (\beta(X))=1\underset{N}{_\beta\otimes_\alpha}\beta(X)$
\newline
(ii) $\Gamma (\alpha(X))=\alpha(X)\underset{N}{_\beta\otimes_\alpha}1$ 
\newline
(iii) $\Gamma$ satisfies the co-associativity relation :
\[(\Gamma \underset{N}{_\beta *_\alpha}id)\Gamma =(id \underset{N}{_\beta *_\alpha}\Gamma)\Gamma\]
This last formula makes sense, thanks to the two preceeding ones and
\ref{fiber}. The von Neumann algebra $N$ will be called the basis of $(N, M, \alpha, \beta, \Gamma)$\vspace{5mm}.\newline
If $(N, M, \alpha, \beta, \Gamma)$ is a Hopf-bimodule, it is clear that
$(N^o, M, \beta, \alpha,
\varsigma_N\circ\Gamma)$ is another Hopf-bimodule, we shall call the symmetrized of the first
one. (Recall that $\varsigma_N\circ\Gamma$ is a homomorphism from $M$ to
$M\underset{N^o}{_r*_s}M$).
\newline
If $N$ is abelian, $\alpha=\beta$, $\Gamma=\varsigma_N\circ\Gamma$, then the quadruplet $(N, M, \alpha, \alpha,
\Gamma)$ is equal to its symmetrized Hopf-bimodule, and we shall say that it is a symmetric
Hopf-bimodule\vspace{5mm}.\newline
Let $\mathcal G$ be a measured groupoid, with $\mathcal G^{(0)}$ as its set of units, and let us denote
by $r$ and $s$ the range and source applications from $\mathcal G$ to $\mathcal G^{(0)}$, given by
$xx^{-1}=r(x)$ and $x^{-1}x=s(x)$. As usual, we shall denote by $\mathcal G^{(2)}$ (or $\mathcal
G^{(2)}_{s,r}$) the set of composable elements, i.e. 
\[\mathcal G^{(2)}=\{(x,y)\in \mathcal G^2; s(x)=r(y)\}\]
Let $(\lambda^u)_{u\in \mathcal G^{(0)}}$ be a Haar system on $\mathcal G$ and $\nu$ a measure $\nu$ on $\mathcal G^{(0)}$. Let us write $\mu$ the measure on $\mathcal G$ given by integrating $\lambda^u$ by $\nu$ :
\[\mu=\int_{{\mathcal G}^{(0)}}\lambda^ud\nu\]
By definition,  $\nu$ is said quasi-invariant if $\mu$ is equivalent to its image under the inverse $x\mapsto x^{-1}$ of $\mathcal G$ (see \cite{Ra}, \cite{R},
 \cite{C2} II.5, \cite{P} for more details and examples of groupoids). 
\newline
In \cite{Y1} and \cite{Val2} were associated to a measured groupoid $\mathcal G$, equipped with a Haar system $(\lambda^u)_{u\in \mathcal G ^{(0)}}$ and a quasi-invariant measure $\nu$ on $\mathcal G ^{(0)}$ two
Hopf-bimodules : 
\newline
The first one is $(L^\infty (\mathcal G^{(0)}, \nu), L^\infty (\mathcal G, \mu), r_{\mathcal G}, s_{\mathcal G}, \Gamma_{\mathcal
G})$, where we define $r_{\mathcal G}$ and $s_{\mathcal G}$ by writing , for $g$ in $L^\infty (\mathcal G^{(0)})$ :
\[r_{\mathcal G}(g)=g\circ r\]
\[s_{\mathcal G}(g)=g\circ s\]
 and where
$\Gamma_{\mathcal G}(f)$, for $f$ in $L^\infty (\mathcal G)$, is the function defined on $\mathcal G^{(2)}$ by $(s,t)\mapsto f(st)$;
$\Gamma_{\mathcal G}$ is then an involutive homomorphism from $L^\infty (\mathcal G)$ into $L^\infty
(\mathcal G^2_{s,r})$ (which can be identified to
$L^\infty (\mathcal G){_s*_r}L^\infty (\mathcal G)$).
\newline
The second one is symmetric; it is $(L^\infty (\mathcal G^{(0)}, \nu), \mathcal L(\mathcal G), r_{\mathcal G}, r_{\mathcal G},
\widehat{\Gamma_{\mathcal G}})$, where
$\mathcal L(\mathcal G)$ is the von Neumann algebra generated by the convolution algebra associated to the
groupoid
$\mathcal G$, and $\widehat{\Gamma_{\mathcal G}}$ has been defined in \cite{Y1} and
\cite{Val2}. 

\subsection{Definition}
\label{defmult}
Let $N$ be a von Neumann algebra; let
$\gH$ be a Hilbert space on which $N$ has a non-degenerate normal representation $\alpha$ and two
non-degenerate normal anti-representations $\hat{\beta}$ and $\beta$. These 3 applications
are supposed to be injective, and to commute two by two.  Let $\nu$ be a normal semi-finite faithful weight on
$N$; we can therefore construct the Hilbert spaces
$\gH\underset{\nu}{_\beta\otimes_\alpha}\gH$ and
$\gH\underset{\nu^o}{_\alpha\otimes_{\hat{\beta}}}\gH$. A unitary $W$ from
$\gH\underset{\nu}{_\beta\otimes_\alpha}\gH$ onto
$\gH\underset{\nu^o}{_\alpha\otimes_{\hat{\beta}}}\gH$
will be called a pseudo-multiplicative unitary over the basis $N$, with respect to the
representation $\alpha$, and the anti-representations $\hat{\beta}$ and $\beta$ (we shall write it is an $(\alpha, \hat{\beta}, \beta)$-pseudo-multiplicative unitary), if :
\newline
(i) $W$ intertwines $\alpha$, $\hat{\beta}$, $\beta$  in the following way :
\[W(\alpha
(X)\underset{N}{_\beta\otimes_\alpha}1)=
(1\underset{N^o}{_\alpha\otimes_{\hat{\beta}}}\alpha(X))W\]
\[W(1\underset{N}{_\beta\otimes_\alpha}\beta
(X))=(1\underset{N^o}{_\alpha\otimes_{\hat{\beta}}}\beta (X))W\]
\[W(\hat{\beta}(X) \underset{N}{_\beta\otimes_\alpha}1)=
(\hat{\beta}(X)\underset{N^o}{_\alpha\otimes_{\hat{\beta}}}1)W\]
\[W(1\underset{N}{_\beta\otimes_\alpha}\hat{\beta}(X))=
(\beta(X)\underset{N^o}{_\alpha\otimes_{\hat{\beta}}}1)W\]
(ii) The operator $W$ satisfies :
\[(1_\gH\underset{N^o}{_\alpha\otimes_{\hat{\beta}}}W)
(W\underset{N}{_\beta\otimes_\alpha}1_{\gH})
=(W\underset{N^o}{_\alpha\otimes_{\hat{\beta}}}1_{\gH})
\sigma^{2,3}_{\alpha, \beta}(W\underset{N}{_{\hat{\beta}}\otimes_\alpha}1)
(1_{\gH}\underset{N}{_\beta\otimes_\alpha}\sigma_{\nu^o})
(1_{\gH}\underset{N}{_\beta\otimes_\alpha}W)\]
Here, $\sigma^{2,3}_{\alpha, \beta}$
goes from $(H\underset{\nu^o}{_\alpha\otimes_{\hat{\beta}}}H)\underset{\nu}{_\beta\otimes_\alpha}H$ to $(H\underset{\nu}{_\beta\otimes_\alpha}H)\underset{\nu^o}{_\alpha\otimes_{\hat{\beta}}}H$, 
and $1_{\gH}\underset{N}{_\beta\otimes_\alpha}\sigma_{\nu^o}$ goes from $H\underset{\nu}{_\beta\otimes_\alpha}(H\underset{\nu^o}{_\alpha\otimes_{\hat{\beta}}}H)$ to $H\underset{\nu}{_\beta\otimes_\alpha}H\underset{\nu}{_{\hat{\beta}}\otimes_\alpha}H$. 
\newline
All the properties supposed in (i) allow us to write such a formula, which will be called the
"pentagonal relation". 
\newline
If $W$ is an $(\alpha, \hat{\beta}, \beta)$-pseudo-multiplicative unitary, then the unitary $\sigma_\nu W^*\sigma_\nu$ from $\gH\underset{\nu}{_{\hat{\beta}}\otimes_\alpha}\gH$ to $\gH\underset{\nu^o}{_\alpha\otimes_\beta}\gH$ is an $(\alpha, \beta, \hat{\beta})$-pseudo-multiplicative unitary, called the dual of $W$ and denoted $\widehat{W}$.

\subsection{Algebras and Hopf-bimodules associated to a pseudo-multiplicative unitary}
\label{AW}
For $\xi_2$ in $D(_\alpha\gH, \nu)$, $\eta_2$ in $D(\gH_{\hat{\beta}}, \nu^o)$, the operator $(\rho_{\eta_2}^{\alpha,
\hat{\beta}})^*W\rho_{\xi_2}^{\beta, \alpha}$ will be written $(id*\omega_{\xi_2, \eta_2})(W)$; we have, therefore, for all
$\xi_1$, $\eta_1$ in $\gH$ :
\[((id*\omega_{\xi_2, \eta_2})(W)\xi_1|\eta_1)=(W(\xi_1\underset{\nu}{_\beta\otimes_\alpha}\xi_2)|
\eta_1\underset{\nu^o}{_\alpha\otimes_{\hat{\beta}}}\eta_2)\]
and, using the intertwining property of $W$ with $\hat{\beta}$, we easily get that $(id*\omega_{\xi_2, \eta_2})(W)$ belongs
to $\hat{\beta} (N)'$. 
\newline
If $x$ belongs to $N$, we have :
\[(id*\omega_{\xi_2, \eta_2})(W)\alpha (x)=(id*\omega_{\xi_2, \alpha(x^*)\eta_2})(W)\]
\[\beta(x)(id*\omega_{\xi_2, \eta_2})(W)=(id*\omega_{\hat{\beta}(x)\xi_2, \eta_2})(W)\]
We shall write $A_n(W)$ (resp. $A_w(W)$) the norm (resp. weak) closure of the linear span of these operators; $A_n(W)$ and $A_w(W)$ are right $\alpha(N)$-modules and left $\beta(N)$-modules. Applying (\cite{E2} 3.6), we get that $A_n(W)$, $A_n(W)^*$, $A_w(W)$ and $A_w(W)^*$ are non-degenerate algebras (one should note that the notations of (\cite{E2}) had been changed in order to fit with Lesieur's notations). 
We shall write $\mathcal A(W)$ the von Neumann algebra generated by $A_w(W)$ .
We then have $\mathcal A(W)\subset\hat{\beta}(N)'$.
\newline
For $\xi_1$ in $D(\gH_\beta,\nu^o)$, $\eta_1$ in $D(_\alpha\gH, \nu)$, we shall write $(\omega_{\xi_1,
\eta_1}*id)(W)$ for the operator $(\lambda_{\eta_1}^{\alpha,
\hat{\beta}})^*W\lambda_{\xi_1}^{\beta, \alpha}$; we have,
therefore, for all
$\xi_2$,
$\eta_2$ in
$\gH$ :
\[((\omega_{\xi_1,\eta_1}*id)(W)\xi_2|\eta_2)=(W(\xi_1\underset{\nu}{_\beta\otimes_\alpha}\xi_2)|
\eta_1\underset{\nu^o}{_\alpha\otimes_{\hat{\beta}}}\eta_2)\]
and, using the intertwining property of $W$ with $\beta$, we easily get that $(\omega_{\xi_1,
\eta_1}*id)(W)$ belongs to $\beta(N)'$. 
\newline
We shall write $\widehat{A_n(W)}$ (resp. $\widehat{A_w(W)}$) the norm (resp. weak) closure of the linear span of these operators. As $\widehat{A_n(W)}=A_n(\widehat{W})^*$, it is clear that these subspaces are non degenerate algebras; following (\cite{EV} 6.1 and 6.5), we shall write $\widehat{\mathcal A(W})$ the von Neumann algebra generated by  $\widehat{A_w(W)}$. We then have $\widehat{\mathcal A(W)}\subset\beta(N)'$. 
\newline
In (\cite{EV} 6.3 and 6.5), using the pentagonal equation, we got
that
$(N,\mathcal A(W),\alpha,\beta,\Gamma)$, and
$(N,\widehat{\mathcal A(W)}, \alpha, \hat{\beta}, \widehat{\Gamma})$ are Hopf-bimodules, where $\Gamma$ and
$\widehat{\Gamma}$ are defined, for any $x$ in $\mathcal A(W)$ and $y$ in $\widehat{\mathcal
A(W)}$, by :
\[\Gamma(x)=W^*(1\underset{N^o}{_\alpha\otimes_{\hat{\beta}}}x)W\]
\[\widehat{\Gamma}(y)=\sigma_{\nu^o}W(y\underset{N}{_\beta\otimes_\alpha}1)W^*\sigma_\nu\]
(Here also, we have changed the notations of \cite{EV}, in order to fit with Lesieur's notations). 
In (\cite{EV} 6.1(iv)), we had obtained that $x$ in $\mathcal L(\gH)$ belongs to $\mathcal A(W)'$ if and only if $x$ belongs to $\alpha(N)'\cap
\beta(N)'$ and verifies $(x\underset{N^o}{_\alpha\otimes_{\hat{\beta}}}1)W=W(x\underset{N}{_\beta\otimes_\alpha}1)$. 
We obtain the same way
that $y$ in $\mathcal L(\gH)$ belongs to $\widehat{\mathcal A(W)}'$ if and only if $y$ belongs to $\alpha(N)'\cap
\hat{\beta}(N)'$ and verify $(1\underset{N^o}{_\alpha\otimes_{\hat{\beta}}}y)W=W(1\underset{N}{_\beta\otimes_\alpha}y)$.  \newline
Moreover, we get that $\alpha(N)\subset\mathcal A\cap\widehat{\mathcal A}$, $\beta(N)\subset\mathcal A$,
$\hat{\beta}(N)\subset\widehat{\mathcal A}$, and, for all $x$ in $N$ :
\[\Gamma (\alpha (x))=\alpha (x)\underset{N}{_\beta\otimes_\alpha}1\]
\[\Gamma (\beta (x))=1\underset{N}{_\beta\otimes_\alpha}\beta (x)\]
\[\widehat{\Gamma}(\alpha(x))=\alpha (x)\underset{N}{_{\hat{\beta}}\otimes_\alpha}1\]
\[\widehat{\Gamma}(\hat{\beta}(x))=1\underset{N}{_{\hat{\beta}}\otimes_\alpha}\hat{\beta}(x)\]

\subsection{Fundamental example}
\label{gd}
Let $\mathcal G$ be a measured groupoid; let's use all notations introduced in \ref{Hbimod}. Let us
note :
\[\mathcal G^2_{r,r}=\{(x,y)\in \mathcal G^2, r(x)=r(y)\}\]
Then, it has been shown \cite{Val2} that the formula $W_{\mathcal G}f(x,y)=f(x,x^{-1}y)$, where $x$, $y$ are
in
$\mathcal G$, such that $r(y)=r(x)$, and $f$ belongs to $L^2(\mathcal G^{(2)})$ (with respect to an
appropriate measure, constructed from $\lambda^u$ and $\nu$), is a unitary from $L^2(\mathcal G^{(2)})$ to $L^2(\mathcal G^2_{r,r})$ (with respect also to another
appropriate measure, constructed from $\lambda^u$ and $\nu$). 
\newline
Let us define $r_{\mathcal G}$ and $s_{\mathcal G}$ from
$L^\infty (\mathcal G^{(0)}, \nu)$ to $L^\infty (\mathcal G, \mu)$ (and then considered as representations on $\mathcal L(L^2(\mathcal
G, \mu))$, for any
$f$ in
$L^\infty (\mathcal G^{(0)}, \nu)$, by
$r_{\mathcal G}(f)=f\circ r$ and $s_{\mathcal G}(f)=f\circ s$.
\newline
We shall identify (\cite{Y1}, 3.2.2) the Hilbert space $L^2(\mathcal G^{(2)})$ with the relative Hilbert tensor product $L^2(\mathcal G, \mu)\underset{L^{\infty}(\mathcal G^{(0)}, \nu)}{_{s_{\mathcal G}}\otimes_{r_{\mathcal G}}}L^2(\mathcal G, \mu)$, and the Hilbert space $L^2(\mathcal G^2_{r,r})$ with the relative Hilbert tensor product $L^2(\mathcal G, \mu)\underset{L^{\infty}(\mathcal G^{(0)}, \nu)}{_{r_{\mathcal G}}\otimes_{r_{\mathcal G}}}L^2(\mathcal G, \mu)$. Moreover, the unitary $W_{\mathcal G}$ can be then interpreted \cite{Val3} as a pseudo-multiplicative unitary over the basis
$L^\infty (\mathcal G^{(0)}, \nu)$, with respect to the representation $r_{\mathcal G}$, and anti-representations
$s_{\mathcal G}$ and
$r_{\mathcal G}$ (as here the basis is abelian, the notions of representation and anti-representations are the same, and the commutation property is fulfilled). So, we get that $W_{\mathcal G}$ is a $(r_\mathcal G, s_\mathcal G, r_\mathcal G)$ pseudo-multiplicative unitary. 
\newline
Let us take the notations of \ref{AW}; the von Neumann algebra $\mathcal A(W_{\mathcal G})$ is equal to the von Neumann algebra $L^{\infty}(\mathcal G, \nu)$ (\cite{Val3}, 3.2.6 and 3.2.7); using (\cite{Val3}
3.1.1), we get that the Hopf-bimodule homomorphism
$\Gamma$ defined on
$L^{\infty}(\mathcal G, \mu)$ by $W_{\mathcal G}$ is equal to the usual Hopf-bimodule homomorphism $\Gamma_{\mathcal G}$ studied in \cite{Val2}, and
recalled in
\ref{Hbimod}.
Moreover, the von Neumann algebra $\widehat{\mathcal A(W_{\mathcal G})}$ is equal to the von Neumann algebra $\mathcal
L(\mathcal G)$ (\cite{Val3}, 3.2.6 and 3.2.7); using (\cite{Val3} 3.1.1), we get that the Hopf-bimodule homomorphism $\widehat{\Gamma}$ defined on $\mathcal L(\mathcal
G)$ by
$W_{\mathcal G}$ is the usual Hopf-bimodule homomorphism $\widehat{\Gamma_{\mathcal G}}$ studied in \cite{Y1} and \cite{Val2}. 

\subsection{Lemma}
\label{lem1}
{\it Let $W$ be an $(\alpha, \hat{\beta}, \beta)$-pseudo-multiplicative unitary, $\xi_1$ in $D(\gH_\beta, \nu^o)$, $\xi_2$ in $D(_\alpha \gH, \nu)$, $\eta$ in $\gH$; let $\zeta_i$ in $D(\gH_\beta, \nu^o)$ and $\zeta'_i$ in $\gH$ such that $W^*(\xi_2\underset{\nu^o}{_\alpha\otimes_{\hat{\beta}}}\eta)=\sum_i \zeta_i\underset{\nu}{_\beta\otimes_\alpha}\zeta'_i$; then we have :}
\[\sum_i\alpha(\langle\zeta_i, \xi_1\rangle_{\beta, \nu^o})\zeta'_i=(\omega_{\xi_1, \xi_2}*id)(W)^*\eta\]

\begin{proof}
Let $\theta$ in $\gH$; we have :
\begin{eqnarray*}
((\omega_{\xi_1, \xi_2}*id)(W)^*\eta|\theta)
&=&(W^*(\xi_2\underset{\nu^o}{_\alpha\otimes_{\hat{\beta}}}\eta)|\xi_1\underset{\nu}{_\beta\otimes_\alpha}\theta)\\
&=&(\sum_i \zeta_i\underset{\nu}{_\beta\otimes_\alpha}\zeta'_i|\xi_1\underset{\nu}{_\beta\otimes_\alpha}\theta)\\
&=&(\sum_i\alpha(\langle\zeta_i, \xi_1\rangle_{\beta, \nu^o})\zeta'_i|\theta)
\end{eqnarray*}
from which we get the result.  \end{proof}

\subsection{Lemma}
\label{lem2}
{\it Let $W$ be an $(\alpha, \hat{\beta}, \beta)$-pseudo-multiplicative unitary, $\xi_1$, $\zeta_1$ in $D(\gH_\beta, \nu^o)$, $\xi$ in $D(_\alpha\gH, \nu)$ and $\eta_1$, $\eta_2$ in $\gH$. Let us consider the flip $\sigma^{1,2}_{\hat{\beta}, \alpha}$ from $H\underset{\nu}{_\beta\otimes_\alpha}(H\underset{\nu^o}{_\alpha\otimes_{\hat{\beta}}}H)$ onto $H\underset{\nu^o}{_\alpha\otimes_{\hat{\beta}}}(H\underset{\nu}{_\beta\otimes_\alpha}H)$. Then, we have :}
\begin{multline*}
(\sigma^{1,2}_{\hat{\beta}, \alpha}
(1_{\gH}\underset{N}{_\beta\otimes_\alpha}W)(\xi_1\underset{\nu}{_\beta\otimes_\alpha}\eta_1\underset{\nu}{_\beta\otimes_\alpha}\xi)|\eta_2\underset{\nu^o}{_\alpha\otimes_{\hat{\beta}}}(\zeta_1\underset{\nu}{_\beta\otimes_\alpha}\zeta_2))=\\
(W(\eta_1\underset{\nu}{_\beta\otimes_\alpha}\xi)|\eta_2\underset{\nu^o}{_\alpha\otimes_{\hat{\beta}}}\alpha(\langle\zeta_1, \xi_1\rangle_{\beta, \nu^o})\zeta_2)
\end{multline*}
\begin{proof}
The scalar product 
\[(\sigma^{1,2}_{\hat{\beta}, \alpha}
(1_{\gH}\underset{N}{_\beta\otimes_\alpha}W)(\xi_1\underset{\nu}{_\beta\otimes_\alpha}\eta_1\underset{\nu}{_\beta\otimes_\alpha}\xi)|\eta_2\underset{\nu^o}{_\alpha\otimes_{\hat{\beta}}}(\zeta_1\underset{\nu}{_\beta\otimes_\alpha}\zeta_2))\]
is equal to :
\[(\xi_1\underset{\nu}{_\beta\otimes_\alpha}W(\eta_1\underset{\nu}{_\beta\otimes_\alpha}\xi)|
\zeta_1\underset{\nu}{_\beta\otimes_\alpha}(\eta_2\underset{\nu^o}{_\alpha\otimes_{\hat{\beta}}}\zeta_2))\]
from which we get the result. \end{proof}

\subsection{Proposition}
\label{prop2Gamma}
{\it Let $W$ be an $(\alpha, \hat{\beta}, \beta)$-pseudo-multiplicative unitary, $\Gamma$ the coproduct constructed in \ref{AW}, $\xi$ in $D(_\alpha\gH, \nu)$, $\eta$ in $D(\gH_{\hat{\beta}}, \nu^o)$. Let $\xi_1$, $\eta_1$ in $D(\gH_\beta, \nu^o)$, $\xi_2$, $\eta_2$ in $D(_\alpha\gH, \nu)$; then, we have :}
\[(\Gamma((id*\omega_{\xi, \eta})(W))(\xi_1\underset{\nu}{_\beta\otimes_\alpha}\eta_1)|
\xi_2\underset{\nu}{_\beta\otimes_\alpha}\eta_2)=
((\omega_{\xi_1, \xi_2}*id)(W)(\omega_{\eta_1, \eta_2}*id)(W)\xi|\eta)\]
\begin{proof}
Using the definition of $\Gamma$ (\ref{AW}), we get that :
\[(\Gamma((id*\omega_{\xi, \eta})(W))(\xi_1\underset{\nu}{_\beta\otimes_\alpha}\eta_1)|
\xi_1\underset{\nu}{_\beta\otimes_\alpha}\eta_2)=
((1\underset{\nu^o}{_\alpha\otimes_{\hat{\beta}}}(id*\omega_{\xi, \eta})(W))W(\xi_1\underset{\nu}{_\beta\otimes_\alpha}\eta_1)|W(\xi_2\underset{\nu}{_\beta\otimes_\alpha}\eta_2))\]
which is equal to :
\[((1\underset{N^o}{_\alpha\otimes_{\hat{\beta}}}W)(W\underset{N}{_\beta\otimes_\alpha}1)(\xi_1\underset{\nu}{_\beta\otimes_\alpha}\eta_1\underset{\nu}{_\beta\otimes_\alpha}\xi)|(W\underset{N^o}{_\alpha\otimes_{\hat{\beta}}}1)((\xi_2\underset{\nu}{_\beta\otimes_\alpha}\eta_2)\underset{\nu^o}{_\alpha\otimes_{\hat{\beta}}}\eta)\]
which, using the pentagonal equation (\ref{defmult}), is equal to :
\[(\sigma^{2,3}_{\alpha, \beta}(W\underset{N}{_{\hat{\beta}}\otimes_\alpha}1)
(1_{\gH}\underset{N}{_\beta\otimes_\alpha}\sigma_{\nu^o})
(1_{\gH}\underset{N}{_\beta\otimes_\alpha}W)(\xi_1\underset{\nu}{_\beta\otimes_\alpha}\eta_1\underset{\nu}{_\beta\otimes_\alpha}\xi)|
(\xi_2\underset{\nu}{_\beta\otimes_\alpha}\eta_2)\underset{\nu^o}{_\alpha\otimes_{\hat{\beta}}}\eta)\] 
or, to :
\[((W\underset{N}{_{\hat{\beta}}\otimes_\alpha}1)
(1_{\gH}\underset{N}{_\beta\otimes_\alpha}\sigma_{\nu^o})
(1_{\gH}\underset{N}{_\beta\otimes_\alpha}W)(\xi_1\underset{\nu}{_\beta\otimes_\alpha}\eta_1\underset{\nu}{_\beta\otimes_\alpha}\xi)|(\xi_2\underset{\nu^o}{_\alpha\otimes_{\hat{\beta}}}\eta)\underset{\nu}{_\beta\otimes_\alpha}\eta_2)\]
which is equal to :
\[(\sigma^{1,2}_{\hat{\beta}, \alpha}
(1_{\gH}\underset{N}{_\beta\otimes_\alpha}W)(\xi_1\underset{\nu}{_\beta\otimes_\alpha}\eta_1\underset{\nu}{_\beta\otimes_\alpha}\xi)|\eta_2\underset{\nu^o}{_\alpha\otimes_{\hat{\beta}}}(W^*(\xi_2\underset{\nu}{_\alpha\otimes_{\hat{\beta}}}\eta)))\]
Defining now $\zeta_i$, $\zeta'_i$ as in \ref{lem1}, we get, using \ref{lem2}, that it is equal to :
\[(W(\eta_1\underset{\nu}{_\beta\otimes_\alpha}\xi)|\eta_2\underset{\nu^o}{_\alpha\otimes_{\hat{\beta}}}\sum_i\alpha(\langle\zeta_i, \xi_1\langle_{\beta, \nu^o})\zeta'_i)\]
which, thanks to \ref{lem1}, is equal to :
\[(W(\eta_1\underset{\nu}{_\beta\otimes_\alpha}\xi)|\eta_2\underset{\nu^o}{_\alpha\otimes_{\hat{\beta}}}(\omega_{\xi_1, \xi_2}*id)(W)^*\eta)\]
and, therefore, to 
\[((\omega_{\eta_1, \eta_2}*id)(W)\xi|(\omega_{\xi_1, \xi_2}*id)(W)^*\eta)\]
which finishes the proof.  \end{proof}

\subsection{Definitions (\cite{L1}, \cite{L2})}
\label{LW}
Let $(N, M, \alpha, \beta, \Gamma)$ be a Hopf-bimodule, as defined in \ref{Hbimod}; a normal, semi-finite, faithful operator valued weight $T$ from $M$ to $\alpha (N)$ is said to be left-invariant if, for all $x\in \gM_T^+$, we have :
\[(id\underset{N}{_\beta*_\alpha}T)\Gamma (x)=T(x)\underset{N}{_\beta\otimes_\alpha}1\]
or, equivalently (\ref{fiber}), if we write $\Phi=\nu\circ\alpha^{-1}\circ T$ :
\[(id\underset{N}{_\beta*_\alpha}\Phi)\Gamma (x)=T(x)\]
A normal, semi-finite, faithful operator-valued weight $T'$ from $M$ to $\beta (N)$ will be said to be right-invariant if it is left-invariant with respect to the symmetrized Hopf-bimodule, i.e., if, for all $x\in\gM_{T'}^+$, we have :
\[(T'\underset{N}{_\beta*_\alpha}id)\Gamma (x)=1\underset{N}{_\beta\otimes_\alpha}T'(x)\]
or, equivalently, if we write $\Psi=\nu\circ\beta^{-1}\circ T'$ : 
\[(\Psi\underset{N}{_\beta*_\alpha}id)\Gamma (x)=T'(x)\]

\subsection{Theorem(\cite{L1}, \cite{L2})}
\label{thL1}
{\it Let $(N, M, \alpha, \beta, \Gamma)$ be a Hopf-bimodule, as defined in \ref{Hbimod}, and let $T$ be a left-invariant normal, semi-finite, faithful operator valued weight from $M$ to $\alpha (N)$; let us choose a normal, semi-finite, faithful weight $\nu$ on $N$, and let us write $\Phi=\nu\circ\alpha^{-1}\circ T$, which is a normal, semi-finite, faithful weight on $M$; let us write $H_\Phi$, $J_\Phi$, $\Delta_\Phi$ for the canonical objects of the Tomita-Takesaki theory associated to the weight $\Phi$, and let us define, for $x$ in $N$, $\hat{\beta}(x)=J_\Phi\alpha(x^*)J_\Phi$.
\newline
(i) There exists an unique isometry $U$ from $H_\Phi\underset{\nu^o}{_\alpha\otimes_{\hat{\beta}}}H_\Phi$ to $H_\Phi\underset{\nu}{_\beta\otimes_\alpha}H_\Phi$, such that, for any $(\beta, \nu^o)$-orthogonal 
basis $(\xi_i)_{i\in I}$ of  $(H_\Phi)_\beta$, for any $a$ in $\gN_T\cap\gN_\Phi$ and for any $v$ in $D((H_\Phi)_\beta, \nu^o)$, we have 
\[U(v\underset{\nu^o}{_\alpha\otimes_{\hat{\beta}}}\Lambda_\Phi (a))=\sum_{i\in I} \xi_i\underset{\nu}{_\beta\otimes_\alpha}\Lambda_{\Phi}((\omega_{v, \xi_i}\underset{\nu}{_\beta*_\alpha}id)(\Gamma(a)))\]
(ii) Let us suppose there exists a right-invariant normal, semi-finite, faithful operator valued weight $T'$ from $M$ to $\beta (N)$; then this isometry is a unitary, and $W=U^*$ is an $(\alpha, \hat{\beta}, \beta)$-pseudo-multiplicative unitary from $H_\Phi\underset{\nu}{_\beta\otimes_\alpha}H_\Phi$ to $H_\Phi\underset{\nu^o}{_\alpha\otimes_{\hat{\beta}}}H_\Phi$ which verifies, for any $x$, $y_1$, $y_2$ in $\gN_T\cap\gN_\Phi$ :
\[(i*\omega_{J_\Phi\Lambda_\Phi (y_1^*y_2), \Lambda_\Phi (x)})(W)=
(id\underset{N}{_\beta*_\alpha}\omega_{J_\Phi\Lambda_\Phi(y_2), J_\Phi\Lambda_\Phi(y_1)})\Gamma (x^*)\]
Clearly, the pseudo-multplicative unitary $W$ does not depend upon the choice of the right-invariant operator-valued weight $T'$, and, for any $y$ in $M$, we have : }
\[\Gamma(y)=W^*(1\underset{N^o}{_\alpha\otimes_{\hat{\beta}}}y)W\]

\begin{proof} This is \cite{L2} 3.51 and 3.52. \end{proof}
\subsection{Definitions}
\label{defMQG}

Let us take the notations of \ref{thL1}; let us write $\Psi=\nu\circ\beta^{-1}\circ T'$. We shall say that $\nu$ is relatively invariant with respect to $T$ and $T'$ if the two modular automorphism groups associated to the two weights $\Phi$ and $\Psi$ commute; we then write down\vspace{5mm}: \newline 
{\bf Definition}
\newline
A measured quantum groupoid is an octuplet $(N, M, \alpha, \beta, \Gamma, T, T', \nu)$ such that :
\newline
(i) $(N, M, \alpha, \beta, \Gamma)$ is a Hopf-bimodule, as defined in \ref{Hbimod}, 
\newline
(ii) $T$ is a left-invariant normal, semi-finite, faithful operator valued weight $T$ from $M$ to $\alpha (N)$, as defined in \ref{LW}, 
\newline
(iii) $T'$ is a right-invariant normal, semi-finite, faithful operator-valued weight $T'$ from $M$ to $\beta (N)$, as defined in \ref{LW}, 
\newline
(iv) $\nu$ is normal semi-finite faitfull weight on $N$, which is relatively invariant with respect to $T$ and $T'$\vspace{5mm}. \newline 
{\bf Remark}
These axioms are not Lesieur's axioms, given in (\cite{L2}, 4.1). The equivalence of these axioms with Lesieur's axioms had been written down in \cite{E5}, and is recalled in the appendix of \cite{E6}. 

\subsection{Theorem (\cite{L2}, \cite{E6})}
\label{thL2}
{\it Let $\gG=(N, M, \alpha, \beta, \Gamma, T, T', \nu)$ be a measured quantum groupoid in the sense of \ref{defMQG}. Let us write $\Phi=\nu\circ\alpha^{-1}\circ T$, which is a normal, semi-finite faithful weight on $M$. Then 
\newline
(i) there exists a $*$-antiautomorphism $R$ on $M$, such that $R^2=id$, $R(\alpha(n))=\beta(n)$ for all $n\in N$, and :
\[\Gamma\circ R=\varsigma_{N^o}(R\underset{N}{_\beta*_\alpha}R)\Gamma\]
$R$ will be called the coinverse; 
\newline
(ii) there exists a one-parameter group $\tau_t$ of automorphisms of $M$, such that $R\circ\tau_t=\tau_t\circ R$ for all $t\in\mathbb{R}$, and, for all $t\in\mathbb{R}$ and $n\in N$, $\tau_t(\alpha(n))=\alpha(\sigma^\nu_t(n))$, $\tau_t(\beta(n))=\beta(\sigma^\nu_t(n))$ and :
\[\Gamma\circ\tau_t=(\tau_t\underset{N}{_\beta*_\alpha}\tau_t)\Gamma=(\sigma_t^\Phi\underset{N}{_\beta*_\alpha}\sigma_{-t}^{\Phi\circ R})\Gamma\]
\[\Gamma\circ\sigma_t^\Phi=(\tau_t\underset{N}{_\beta*_\alpha}\sigma_t^\Phi)\Gamma\]
$\tau_t$ will be called the scaling group;
\newline
(iii) the weight $\nu$ is relatively invariant with respect to $T$ and $RTR$; moreover, $R$ and $\tau_t$ are still the co-inverse and the scaling group of this new measured quantum groupoid, we shall denote :
\[\underline{\gG}=(N, M, \alpha, \beta, \Gamma, T, RTR, \nu) \]
\newline
(iv) for any $\xi$, $\eta$ in $D(_\alpha H_\Phi, \nu)\cap D((H_\Phi)_{\hat{\beta}}, \nu^o)$, $(id*\omega_{\xi, \eta})(W)$ belongs to $D(\tau_{i/2})$, and, if we define $S=R\tau_{-i/2}$, we have :
\[S((id*\omega_{\xi, \eta})(W))=(id*\omega_{\eta, \xi})(W)^*\]
More generally, for any $x$ in $D(S)=D(\tau_{-i/2})$, we get that $S(x)^*$ belongs to $D(S)$ and $S(S(x)^*)^*=x$; 
$S$ will be called the antipode of $\gG$ (or $\underline{\gG}$), and, therefore, the co-inverse and the scaling group, given by polar decomposition of the antipode, rely only on the pseudo-multiplicative $W$. 
\newline
(v) there exists a one-parameter group $\gamma_t$ of automorphisms of $N$ such that, for all $t\in\mathbb{R}$ and $n\in N$, we have :
\[\sigma_t^{T}(\beta(n))=\beta(\gamma_t(n))\]
Moreover, for all $t\in\mathbb{R}$, we have $\nu\circ\gamma_t=\nu$. 
\newline
(vi) there exists a positive non-singular operator $\lambda$ affiliated to $Z(M)$, and a positive non singular operator $\delta$ affiliated to $M$, such that :
\[(D\Phi\circ R: D\Phi)_t=\lambda^{it^2/2}\delta^{it}\]
and, therefore, we have :
\[(D\Phi\circ\sigma_s^{\Phi\circ R}:D\Phi)_t=\lambda^{ist}\]
The operator $\lambda$ will be called the scaling operator, and there exists a positive non-singular operator $q$ affiliated to $N$ such that $\lambda=\alpha(q)=\beta(q)$. We have $R(\lambda)=\lambda$. 
\newline
The operator $\delta$ will be called the modulus; we have $R(\delta)=\delta^{-1}$, and $\tau_t(\delta)=\delta$, for all $t\in\mathbb{R}$, and we can define a one-parameter group of unitaries $\delta^{it}\underset{N}{_\beta\otimes_\alpha}\delta^{it}$ which acts naturally on elementary tensor products, and verifies, for all $t\in\mathbb{R}$ :
\[\Gamma(\delta^{it})=\delta^{it}\underset{N}{_\beta\otimes_\alpha}\delta^{it}\]
(vii) we have $(D\Phi\circ\tau_t : D\Phi)_s=\lambda^{-ist}$, which leads to define a one-parameter group of unitaries $P^{it}$ by, for any $x\in\gN_\Phi$ :
\[P^{it}\Lambda_\Phi(x)=\lambda^{t/2}\Lambda_\Phi(\tau_t(x))\]
Moreover, for any $y$ in $M$, we get :
\[\tau_t(y)=P^{it}yP^{-it}\]
and it is possible to define one parameter groups of unitaries $P^{it}\underset{N}{_\beta\otimes_\alpha}P^{it}$ and $P^{it}\underset{N^o}{_\alpha\otimes_{\hat{\beta}}}P^{it}$ such that :
\[W(P^{it}\underset{N}{_\beta\otimes_\alpha}P^{it})=(P^{it}\underset{N^o}{_\alpha\otimes_{\hat{\beta}}}P^{it})W\]
Moreover, for all $v\in D(P^{-1/2})$, $w\in D(P^{1/2})$, $p$, $q$ in $D(_\alpha H_\Phi, \nu)\cap D((H_\Phi)_{\hat{\beta}}, \nu^o)$, we have:
\[(W^*(v\underset{\nu^o}{_\alpha\otimes_{\hat{\beta}}}q)|w\underset{\nu}{_\beta\otimes_\alpha}p)=
(W(P^{-1/2}v\underset{\nu}{_\beta\otimes_\alpha}J_\Phi p)|P^{1/2}w\underset{\nu^o}{_\alpha\otimes_{\hat{\beta}}}J_\Phi q)\]
We shall say that the pseudo-multiplicative unitary $W$ is "manageable", with "managing operator" $P$, which implies (with the notations of \ref{AW}) that $A_w(W)=\mathcal A(W)=M$ and $\widehat{A_w(W)}=\widehat{{\mathcal A}(W)}$
\newline
As, for all $s$, $t$ in $\mathbb{R}$, we have $\tau_s\circ\sigma_t^\Phi=\sigma_t^\Phi\circ\tau_s$, we get that $J_\Phi PJ_\Phi=P$. 
\newline
(viii) Let us write $\widehat{M}=\widehat{A_w(W)}=\widehat{{\mathcal A}(W)}$ and let us consider the coproduct $\widehat{\Gamma}$ on $\widehat{M}$, then, by \ref{AW}, $(N, \widehat{M} \alpha, \hat{\beta}, \widehat{\Gamma})$ is a Hopf-bimodule; moreover, there exists a $*$-antiautomorphism $\widehat{R}$ on $\widehat{M}$, such that $\widehat{R}^2=id$, $\widehat{R}\circ\alpha=\hat{\beta}$, and $\widehat{\Gamma}\circ\widehat{R}=\varsigma_{N^o}(\widehat{R}\underset{N}{_{\hat{\beta}}*_\alpha}\widehat{R})\Gamma$ and a left-invariant normal, semi-finite, faithful operator-valued weight $\hat{T}$ from $\widehat{M}$ to $\alpha (N)$.
\newline
(ix) $(N, \widehat{M}, \alpha, \hat{\beta}, \widehat{\Gamma}, \hat{T}, \widehat{R}\hat{T}\widehat{R}, \nu)$ is a measured quantum groupoid, called the dual measured quantum groupoid of $\gG$, we shall denote $\widehat{\gG}$. Moreover, we have $\widehat{\widehat{\gG}}=\underline{\gG}$. }

\begin{proof} This is (\cite{E6}, 3.8, 3.10 and 3.11). \end{proof}

\subsection{Example}
\label{gd2}
Let $\mathcal G$ be a measured groupoid; using the notations introduced in \ref{Hbimod} and \ref{gd}, we have seen that $(L^\infty(\mathcal G^{(0)}, \nu), L^\infty (\mathcal G, \mu), r_{\mathcal G}, s_{\mathcal G}, \Gamma_{\mathcal G})$ is a Hopf-bimodule; moreover, it is possible to prove that the formula which gives, for all positive $F$ in $L^\infty (\mathcal G, \mu)$ the image by $r_\mathcal G$ of the function $u\mapsto \int_\mathcal G Fd\lambda^u$ (resp. the image by $s_\mathcal G$ of the function $u\mapsto \int_\mathcal G Fd\lambda_u$) defines a normal semi-finite faithful operator-valued weight from $L^\infty (\mathcal G, \mu)$ onto $r_\mathcal G(L^\infty(\mathcal G^{(0)}, \nu))$ (resp. $s_\mathcal G(L^\infty(\mathcal G^{(0)}, \nu))$, which is left-invariant (resp. right-invariant) with respect to $\Gamma_\mathcal G$ as defined in \ref{LW}; moreover, as $L^\infty (\mathcal G, \mu)$ is abelian, the measure $\nu$ defines a relatively invariant weight as defined in \ref{defMQG}. Therefore, we obtain a measured quantum groupoid $(L^\infty(\mathcal G^{(0)}, \nu), L^\infty (\mathcal G, \mu), r_{\mathcal G}, s_{\mathcal G}, \Gamma_{\mathcal G}, T_\mathcal G, T'_\mathcal G, \nu)$, we shall denote $\gG(\mathcal G)$. The dual $\widehat{\gG(\mathcal G)}$ is symmetric : it is $(L^\infty(\mathcal G^{(0)}, \nu), \mathcal L(\mathcal G), r_\mathcal G, r_\mathcal G, \hat{T}_\mathcal G, \hat{T}_\mathcal G, \nu)$, where $\hat{T}_\mathcal G$ is a normal semi-finite faithful operator-valued weight from $\mathcal L(\mathcal G)$ to $r_\mathcal G(L^\infty(\mathcal G^{(0)}, \nu))$ which is both left and right-invariant. (\cite{L2}, 10). 

\section{A canonical sub-${\bf C}^*$-algebra of a measured quantum groupoid}
\label{anw}
Be given a measured quantum groupoid $\gG=(N, M, \alpha, \beta, \Gamma, T, T', \nu)$, we consider in this chapter the sub-${\bf C}^*$-algebra $A_n(W)\cap A_n(W)^*$ of $M$; we obtain that it is dense in $M$ (\ref{propmanag}(i)), invariant by $R$ (\ref{thL}(i)), by $\sigma_t^\Phi$ and $\tau_t$ (\ref{corPhi}(i)). The results are more precise when $N$ is abelian : in that case, this ${\bf C}^*$-algebra is equal to $A_n(W)$ (\ref{propmanag}(ii)) and the one-parameter group of automorphisms $\tau_t$ is norm continuous (\ref{corPhi}(iv)). Moreover, we obtain more complete results if the one-parameter group $\gamma_t$ is trivial (which is the case when $\alpha(N)\subset Z(M)$); then, the modular groups $\sigma_t^\Phi$ and $\sigma_t^{\Phi\circ R}$ are norm continuous on $A_n(W)$ (\ref{corPhi}(ii) and (iii)).

\subsection{Lemma}
\label{lemmanageable}
{\it Let $\gG=(N, M, \alpha, \beta, \Gamma, T, T', \nu)$ be a measured quantum groupoid, and let's use the notations of \ref{thL2}. Then, if $p$ belongs to  $D(_\alpha H_\Phi, \nu)\cap D((H_\Phi)_{\hat{\beta}}, \nu^o)\cap\mathcal D(P^{1/2})$ such that $P^{1/2}p$ belongs to $D(_\alpha H_\Phi, \nu)$, and $q$ belongs to $D(_\alpha H_\Phi, \nu)\cap D((H_\Phi)_{\hat{\beta}}, \nu^o)\cap\mathcal D(P^{-1/2})$ such that $P^{-1/2}q$ belongs to $D((H_\Phi)_{\hat{\beta}}, \nu^o)$, then we have :
\[(id*\omega_{J_\Phi p, J_\Phi q})(W)^*=(id*\omega_{P^{1/2}p, P^{-1/2}q})(W)\]
and, therefore, $(id*\omega_{J_\Phi p, J_\Phi q})(W)$ belongs to $A_n(W)\cap A_n(W)^*$.}
\begin{proof}
Let us take $v$ in $\mathcal D(P^{-1/2})$, $w$ in $\mathcal D(P^{1/2})$; then, we have, using \ref{thL2}(vii) :
\begin{eqnarray*}
((id*\omega_{J_\Phi p, J_\Phi q})(W)^*v|w)&=&(v|(id*\omega_{J_\Phi p, J_\Phi q})(W)w)\\
&=&(v\underset{\nu^o}{_\alpha\otimes_{\hat{\beta}}}Jq|W(w\underset{\nu}{_\beta\otimes_\alpha}Jp))\\
&=&
(W(P^{-1/2}v\underset{\nu}{_\beta\otimes_\alpha}p)|P^{1/2}w\underset{\nu^o}{_\alpha\otimes_{\hat{\beta}}}q)\\
&=&(W(v\underset{\nu}{_\beta\otimes_\alpha}P^{1/2}p)|w\underset{\nu^o}{_\alpha\otimes_{\hat{\beta}}}P^{-1/2}q)\\
&=&((id*\omega_{P^{1/2}p, P^{-1/2}q})(W)v|w)
\end{eqnarray*}
which, by density, gives the result. \end{proof}

\subsection{Lemma}
\label{lemP}
{\it Let $\gG=(N, M, \alpha, \beta, \Gamma, T, T', \nu)$ be a measured quantum groupoid, and let's use the notations of \ref{thL2}. Then :
\newline
(i) for any $p$ in $D(_\alpha H_\Phi, \nu)$, there exists a sequence $p_n$ in $D(_\alpha H_\Phi, \nu)\cap\mathcal D(P^{1/2})\cap\mathcal D(P^{-1/2})$, such that $P^{1/2}p_n$ belongs to $D(_\alpha H_\Phi, \nu)$, and such that $R^{\alpha, \nu}(p_n)$ is weakly converging to $R^{\alpha, \nu}(p)$.
\newline(ii) for any $q$ in $D((H_\Phi)_{\hat{\beta}}, \nu^o)$, there exists $q_n$ in $D((H_\Phi)_{\hat{\beta}}, \nu^o)\cap\mathcal D(P^{1/2})\cap\mathcal D(P^{-1/2})$, such that $P^{-1/2}p_n$ belongs to $D((H_\Phi)_{\hat{\beta}}, \nu^o)$, and such that $R^{\hat{\beta}, \nu^o}(p_n)$ is weakly converging to $R^{\hat{\beta}, \nu^o}(p)$. }
 \begin{proof}
Let us write :
\[p_n= \frac{\sqrt{n}}{\pi}\int_{-\infty}^{\infty}e^{-nt^2}P^{it}pdt\]
It is a usual calculation to prove that $p_n$ belongs to $\mathcal D(P^{1/2})\cap \mathcal D(P^{-1/2})$; moreover, we get, for any $a$ in $\gN_\nu$ :
\begin{eqnarray*}
\alpha(a)p_n&=&\frac{\sqrt{n}}{\pi}\int_{-\infty}^{\infty}e^{-nt^2}\alpha(a)P^{it}pdt\\
&=&\frac{\sqrt{n}}{\pi}\int_{-\infty}^{\infty}e^{-nt^2}P^{it}\alpha(\sigma^{\nu}_{-t}(a))pdt\\
&=&\frac{\sqrt{n}}{\pi}\int_{-\infty}^{\infty}e^{-nt^2}P^{it}R^{\alpha, \nu}(p)\Delta_{\nu}^{-it}\Lambda_{\nu}(a)dt
\end{eqnarray*}
from which we get that :
\[\|\alpha(a)p_n\|\leq\frac{\sqrt{n}}{\pi}\int_{-\infty}^{\infty}e^{-nt^2}\|R^{\alpha, \nu}(p)\|\|\Lambda_{\nu}(a)\|dt\]
which proves that $p_n$ belongs to $D(_\alpha H_\Phi, \nu)$ and that :
\[\|R^{\alpha, \nu}(p_n)\|\leq\|R^{\alpha, \nu}(p)\|\]
Moreover, we have, going on the same calculation :
\[R^{\alpha, \nu}(p_n)\Lambda_\nu (a)=\frac{1}{\pi}\int_{-\infty}^{\infty}e^{-t^2}P^{\frac{it}{\sqrt{n}}}R^{\alpha, \nu}(p)\Delta_{\nu}^{\frac{-it}{\sqrt{n}}}\Lambda_{\nu}(a)dt\]
which, using Lebesgue's theorem, is converging to $R^{\alpha, \nu}(p)\Lambda_{\nu}(a)$. With the norm majoration, we get this way the weak convergence of $R^{\alpha, \nu}(p_n)$ to $R^{\alpha, \nu}(p)$, which is (i). Part (ii) is obtained the same way.  \end{proof}

\subsection{Proposition}
\label{propmanag}
{\it Let $\gG=(N, M, \alpha, \beta, \Gamma, T, T', \nu)$ be a measured quantum groupoid, and let's use the notations of \ref{thL2}. Then :
\newline
(i) $A_n(W)\cap A_n(W)^*$ is a non degenerate ${\bf C}^*$-algebra, which is weakly dense in $M$. Moreover, if $y\in N$ is analytical with respect to $\nu$, then $\alpha(y)$ and $\beta(y)$ belong to the multipliers of this ${\bf C}^*$-algebra.
\newline
(ii) if $N$ is abelian, then $A_n(W)$ is a non degenerate ${\bf C}^*$-algebra, which is weakly dense in $M$; moreover, we have $\alpha(N)\subset M(A_n(W))$ and $\beta(N)\subset M(A_n(W))$. }

\begin{proof}
Let $\xi$ and $\eta$ in $D(_\alpha H_\Phi, \nu)\cap D((H_\Phi)_{\hat{\beta}}, \nu^o)$; then, using \ref{lemP}, it is possible to construct sequences $p_n$ and $q_n$ such that $R^{\hat{\beta}, \nu^o}(p_n)$ is weakly converging to $R^{\hat{\beta}, \nu^o}(\xi)$ (or, equivalently, $R^{\alpha, \nu}(J_\Phi p_n)$ is weakly converging to $R^{\alpha, \nu}(J_\Phi \xi)$) and $R^{\alpha, \nu}(q_n)$ is weakly converging to $R^{\alpha, \nu}(\eta)$ (or, equivalently, $R^{\hat{\beta}, \nu^o}(J_\Phi q_n)$ is weakly converging to $R^{\hat{\beta}, \nu^o}(J_\Phi \eta)$), and such that, using \ref{lemmanageable}, the operators $(id*\omega_{J_\Phi p_n, J_\Phi q_n})(W)$ belong to $A_n(W)\cap A_n(W)^*$. 
\newline
So, the element $(id*\omega_{J_\Phi\xi, J_\Phi\eta})(W)$ belongs to the weak closure of $A_n(W)\cap A_n(W)^*$. 
\newline
If $\xi$ belongs to $D(_\alpha H_\Phi, \nu)$, and $\eta$ to $D((H_\Phi)_{\hat{\beta}}, \nu^o)$, then, using \ref{ovw}, we obtain that the element $(id*\omega_{\xi, \eta})(W)$ belongs also to the weak closure of $A_n(W)\cap A_n(W)^*$. So, with the notations of \ref{AW}, we get that $A_w(W)$ is included in the weak closure of $A_n(W)\cap A_n(W)^*$, and, using now \ref{thL2}(vii), we get that $M$ is equal to the weak closure of $A_n(W)\cap A_n(W)^*$, which is (i). 
\newline
Let us now suppose that $N$ is abelian; then the weight $\nu$ is a trace, and the managing operator $P$ defined in \ref{thL2}(vii) is affiliated to $\alpha (N)'\cap \hat{\beta}(N)'$. Let us write :
\[P=\int_0^{\infty}e_\lambda de_\lambda\]
and let us define $p_n=\int_{1/n}^n de_\lambda$.
Then $p_n$ is an increasing sequence of projections, weakly converging to $1$, in $\alpha(N)'\cap \hat{\beta}(N)'$. Let us take $x$ in $\mathcal T_{\Phi, T}$ (with the notations of \ref{ovw}); then the vectors $p_n\Lambda_{\Phi}(x)$ belong to $D(_\alpha H_\Phi, \nu)\cap D((H_\Phi)_{\hat{\beta}}, \nu^o)\cap \mathcal D(P^{1/2})\cap D(P^{-1/2})$ and both $P^{1/2}p_n\Lambda_\Phi(x)$  and $P^{-1/2}p_n\Lambda_{\Phi}(x)$ satisfy the hypothesis of \ref{lemmanageable}. So, using \ref{lemmanageable}, we get that, for $x$, $y$ in $\mathcal T_{\Phi, T}$, the operator $(id*\omega_{J_\Phi p_n\Lambda_\Phi(x), J_\Phi p_n\Lambda_{\Phi}(y)})(W)$ belongs to $A_n(W)\cap A_n(W)^*$. 
\newline
Using now (\cite{E2}, 10.5), we get, taking the norm limit, that $(id*\omega_{J_\Phi\Lambda_{\Phi}(x), J_\Phi\Lambda_{\Phi}(y)})(W)$ belongs to $A_n(W)\cap A_n(W)^*$ for any $x$, $y$ in $\mathcal T_{\Phi, T}$, or that $(id*\omega_{\Lambda_{\Phi}(x), \Lambda_{\Phi}(y)})(W)$ belongs to $A_n(W)\cap A_n(W)^*$. Using now \ref{ovw}, we get that $(id*\omega_{\xi, \eta})(W)$ belongs to $A_n(W)\cap A_n(W)^*$, for any $\xi$ in $D(_\alpha H, \nu)$ and $\eta$ in $D(H_{\hat{\beta}}, \nu^o)$, and, therefore, we get that its norm closure $A_n(W)$ is also included in $A_n(W)\cap A_n(W)^*$, which finishes the proof. \end{proof}

\subsection{Theorem}
\label{thL}
{\it Let $\gG=(N, M, \alpha, \beta, \Gamma, T, T', \nu)$ be a measured quantum groupoid; then, for any $\xi$, $\eta$ in $D(_\alpha H_\Phi, \nu)$, for all $t$ in $\mathbb{R}$, we have :
\newline
(i) $R((i*\omega_{\xi, J_\Phi\eta})(W))=(i*\omega_{\eta, J_\Phi\xi})(W)$
\newline
(ii) $\tau_t((i*\omega_{\xi, J_\Phi\eta})(W))=(i*\omega_{\Delta_\Phi^{-it}\xi, \Delta_\Phi^{-it}J_\Phi\eta})(W)$
\newline
(iii)   $\sigma_t^\Phi((i*\omega_{\xi, J_\Phi\eta})(W))=
(i*\omega_{\delta^{it}J_\Phi\delta^{-it}J_\Phi\Delta_\Phi^{-it}\xi, P^{it}J_\Phi\eta})(W)$
\newline
$\sigma_t^{\Phi\circ R}((i*\omega_{\xi, J_\Phi\eta})(W))=
(i*\omega_{P^{it}\xi, \delta^{it}J_\Phi\delta^{-it}J_\Phi\Delta_\Phi^{-it}J_\Phi\eta})(W)$}
\begin{proof}
Results (i) and (ii) are (\cite{L2} 4.6). 
\newline
Let us take $\xi=J_\Phi\Lambda_\Phi (y_1^*y_2)$, and $\eta=J_\Phi\Lambda_\Phi(x)$, with $x$, $y_1$, $y_2$ in  $\gN_T\cap\gN_\Phi$; then, using \ref{thL1}(ii) and \ref{thL2}(ii), we get :
\[\sigma_t^\Phi((i*\omega_{J_\Phi\Lambda_\Phi (y_1^*y_2), \Lambda_\Phi(x)})(W))
=(id\underset{N}{_\beta*_\alpha}\omega_{J_\Phi\Lambda_\Phi(y_2), J_\Phi\Lambda_\Phi(y_1)}\circ\sigma_t^{\Phi\circ R})\Gamma (\tau_t(x^*))\]
which is equal to :
\begin{multline*}
(id\underset{N}{_\beta*_\alpha}\omega_{J_\Phi\Lambda_\Phi(\lambda^{t/2}\sigma_{-t}^{\Phi\circ R}(y_2)), J_\Phi\Lambda_\Phi(\lambda^{t/2}\sigma_{-t}^{\Phi\circ R}(y_1))})\Gamma (\tau_t(x^*))\\
=(id\underset{N}{_\beta*_\alpha}\omega_{J_\Phi\Lambda_\Phi(\sigma_{-t}^{\Phi\circ R}(y_2)), J_\Phi\Lambda_\Phi(\sigma_{-t}^{\Phi\circ R}(y_1))})\Gamma (\lambda^t\tau_t(x^*))
\end{multline*}
which, using again \ref{thL1}(ii) and \ref{thL2}(ii), is equal to :
\[(i*\omega_{J_\Phi\Lambda_\Phi (\sigma_t^{\Phi\circ R}(y_1^*y_2)), \Lambda_\Phi(\lambda^t\tau_t(x))})(W)
=(i*\omega_{J_\Phi\delta^{-it}J_\Phi\delta^{it}J_\Phi\Delta_\Phi^{it}\Lambda_\Phi (y_1^*y_2), P^{it}\Lambda_\Phi(x)})(W)\]
which gives the first result of (iii), using \ref{ovw}. 
\newline
By similar calculations, we obtain :
\[\sigma_t^{\Phi\circ R}((i*\omega_{J_\Phi\Lambda_\Phi (y_1^*y_2), \Lambda_\Phi(x)})(W))
=(id\underset{N}{_\beta*_\alpha}\omega_{J_\Phi\Lambda_\Phi(y_2), J_\Phi\Lambda_\Phi(y_1)}\circ\tau_t)\Gamma (\sigma_t^{\Phi\circ R}(x^*))\]
which is equal to :
\[(id\underset{N}{_\beta*_\alpha}\omega_{J_\Phi\Lambda_\Phi(\lambda^{t/2}\tau_t(y_2), J_\Phi\Lambda_\Phi(\lambda^{t/2}\tau_t(y_1)})\Gamma (\sigma_t^{\Phi\circ R}(x^*))
=(i*\omega_{J_\Phi\Lambda_\Phi (\lambda^t\tau_t(y_1^*y_2)), \Lambda_\Phi(\sigma^{\Phi\circ R}(x))})(W)\]
from which we obtain the second result of (iii). \end{proof}

\subsection{Corollary}
\label{corPhi}
{\it Let $\gG=(N, M, \alpha, \beta, \Gamma, T, T', \nu)$ be a measured quantum groupoid, and $\gamma_t$ the one-parameter group of automorphisms of $N$ defined in \ref{thL2}(v);  let $\xi$, $\eta$ in $D(_\alpha H, \nu)$, then :
\newline 
(i) We have $R(A_n(W))=A_n(W)$, and, for any $t$ in $\mathbb{R}$, we have $\sigma_t^\Phi(A_n(W))=A_n(W)$ and $\tau_t(A_n(W))=A_n(W)$. 
\newline
(ii) if $\langle\xi, \xi\rangle_{\alpha, \nu}^o$ belongs to ${\bf C}^*(\gamma)$ and $\langle\eta, \eta\rangle_{\alpha, \nu}^o$ to ${\bf C}^*(\sigma^\nu)$, then $(i*\omega_{\xi, J_\Phi\eta})(W)$ belongs to ${\bf C}^*(\sigma^\Phi)$; so, if $N$ is abelian, and if $\gamma=id$, we have $A_n(W)\subset {\bf C}^*(\sigma^\Phi)$. 
\newline
(iii) if $\langle\xi, \xi\rangle_{\alpha, \nu}^o$ belongs to ${\bf C}^*(\sigma^\nu)$ and $\langle\eta, \eta\rangle_{\alpha, \nu}^o$ to ${\bf C}^*(\gamma)$, then $(i*\omega_{\xi, J_\Phi\eta})(W)$ belongs to ${\bf C}^*(\sigma^{\Phi\circ R})$; so, if $N$ is abelian, and if $\gamma=id$, we have $A_n(W)\subset {\bf C}^*(\sigma^{\Phi\circ R})$.
\newline
(iv) if $\langle\xi, \xi\rangle_{\alpha, \nu}^o$ and $\langle\eta, \eta\rangle_{\alpha, \nu}^o$ belong to ${\bf C}^*(\sigma^\nu)$, then $(i*\omega_{\xi, J_\Phi\eta})(W)$ belongs to ${\bf C}^*(\tau)$; so, if $N$ is abelian, we have $A_n(W)\subset {\bf C}^*(\tau)$.}
\begin{proof}
By the weak continuity of $x\mapsto \sigma_t^\Phi(x)$ and $x\mapsto \tau_t(x)$, the first results are simple corollaries of \ref{thL}(ii) and (iii). 
Let now be $\nabla$ the self-adjoint positive operator defined on $L^2(N)$ defined, for all $n$ in $\gN_\nu$ and $t$ in $\mathbb{R}$ by :
\[\nabla^{it}\Lambda_\nu(n)=\Lambda_\nu(\gamma_t(n))\]
We have then :
\begin{eqnarray*}
R^{\alpha, \nu}(\delta^{it}J_\Phi\delta^{-it}J_\Phi\Delta_\Phi^{-it}\xi)\Lambda_\nu(n)
&=&\alpha(n)\delta^{it}J_\Phi\delta^{-it}J_\Phi\Delta_\Phi^{-it}\xi\\
&=&\sigma_{-t}^{\Phi\circ R}(\alpha(n))\xi\\
&=&\alpha(\gamma_t(n))\xi\\
&=&R^{\alpha, \nu}(\xi)\nabla^{it}\Lambda_\nu(n)
\end{eqnarray*}
from which we get that $R^{\alpha, \nu}(\delta^{it}J_\Phi\delta^{-it}J_\Phi\Delta_\Phi^{-it}\xi)=R^{\alpha, \nu}(\xi)\nabla^{it}$ and that 
\[\langle\delta^{it}J_\Phi\delta^{-it}J_\Phi\Delta_\Phi^{-it}\xi, \delta^{it}J_\Phi\delta^{-it}J_\Phi\Delta_\Phi^{-it}\xi\rangle_{\alpha, \nu}^o=\gamma_{-t}(\langle\xi, \xi\rangle_{\alpha, \nu}^o)\]
Therefore, if the function $t\mapsto\gamma_t(\langle\xi, \xi\rangle^o)$ is norm continuous, so is the function $t\mapsto\|R^{\alpha}(\delta^{it}J_\Phi\delta^{-it}J_\Phi\Delta_\Phi^{-it}\xi)\|$. 
\newline
On the other hand, we have :
\[R^{\alpha, \nu}(\Delta_\Phi^{-it}\xi)\Lambda_\nu(n)=\alpha(n)\Delta_\Phi^{-it}\xi\\=\Delta_\Phi^{-it}\alpha(\sigma_t^\nu(n))\xi
=\Delta_\Phi^{-it}R^{\alpha, \nu}(\xi)\Delta_\nu^{it}\Lambda_\nu(n)\]
from which we get that $\langle\Delta_\Phi^{-it}\xi, \Delta_\Phi^{-it}\xi\rangle_{\alpha, \nu}^o=\sigma_{t}^\nu(\langle\xi, \xi\rangle_{\alpha, \nu}^o)$; from these results, using \ref{thL}(iii) and \ref{thL}(ii), we get easily (ii), (iii) and (iv).  \end{proof}

\subsection{Proposition}
\label{Phi}
{\it Let $\gG=(N, M, \alpha, \beta, \Gamma, T, T', \nu)$ be a measured quantum groupoid; then :
\newline
(i) if $x$ is in $\gN_T\cap\gN_\Phi$, and $y$ is in $\gN_T\cap\gN_{R\circ T\circ R}\cap\gN_\Phi$, then $(i*\omega_{J_\Phi\Lambda_\Phi (y^*y), \Lambda_\Phi (x^*x)})(W)$ belongs to $\gM_T^+\cap\gM_\Phi^+$ and we have :
}
\[\Phi( (i*\omega_{J_\Phi\Lambda_\Phi (y^*y), \Lambda_\Phi (x^*x)})(W))=(R\circ T\circ R(y^*y)J_\Phi\Lambda_\Phi(x)|J_\Phi\Lambda_\Phi(x))\]
\begin{align*}
T( (i*\omega_{J_\Phi\Lambda_\Phi (y^*y), \Lambda_\Phi (x^*x)})(W))
&=
\alpha (\langle T\circ R(y^*y)J_\Phi\Lambda_\Phi(x),J_\Phi\Lambda_\Phi(x)\rangle_{\alpha, \nu})\\
&\leq\|T\circ R(y^*y)\|T(x^*x)
\end{align*}
{\it (ii) if $x$, $y$ are in $\gN_T\cap\gN_{R\circ T\circ R}\cap\gN_\Phi$, then the operator $(i*\omega_{J_\Phi\Lambda_\Phi (y^*y), \Lambda_\Phi (x^*x)})(W)$ belongs to $\gM_T^+\cap\gM_\Phi^+\cap\gM_{R\circ T\circ T}^+\cap\gM_{\Phi\circ R}^+$. 
\newline
(iii) if $x$, $y$ are in $\gM_T\cap\gM_{R\circ T\circ R}\cap\gM_\Phi\cap\gM_{\Phi\circ R}$, then $(i*\omega_{J_\Phi\Lambda_\Phi (y), \Lambda_\Phi (x)})(W)$ belongs to $\gM_T\cap\gM_\Phi\cap\gM_{R\circ T\circ R}\cap\gM_{\Phi\circ R}$.}
\begin{proof}
Using \ref{LW}, we obtain that the operator $(i*\omega_{J_\Phi\Lambda_\Phi (y^*y), \Lambda_\Phi (x^*x)})(W)$ is positive, and that :
\[R\circ T( (i*\omega_{J_\Phi\Lambda_\Phi (y^*y), \Lambda_\Phi (x^*x)})(W))
= R\circ T\circ R( (i*\omega_{J_\Phi\Lambda_\Phi (x^*x), \Lambda_\Phi (y^*y)})(W))\]
which, using \ref{thL}(iii) and the right-invariance of $R\circ T\circ R$, is equal to :
\begin{multline*}
R\circ T\circ R((id\underset{N}{_\beta *_\alpha}\omega_{J_\Phi\Lambda_\Phi (x)})\Gamma(y^*y))\\
=(id\underset{N}{_\beta *_\alpha}\omega_{J_\Phi\Lambda_\Phi (x)})(1\underset{N}{_\beta\otimes_\alpha}R\circ T\circ R (y^*y))\\
=\beta(\langle R\circ T\circ R(y^*y)J_\Phi\Lambda_\Phi(x),J_\Phi\Lambda_\Phi(x)\rangle_{\alpha, \nu})
\end{multline*}
from which we get (i); we then get (ii) by using \ref{thL}(i), and (iii) is just given by linearity. \end{proof}

\subsection{Lemma}
\label{lemomega}
{\it Let $\gG=(N, M, \alpha, \beta, \Gamma, T, T', \nu)$ be a measured quantum groupoid;  let us define $\Phi=\nu\circ\alpha^{-1}\circ T$; then, we have, for all $x$ in $\gN_\Phi$ :}
\[\omega_{J_\Phi\Lambda_{\Phi}(x)}\circ R=\omega_{J_{\Phi\circ R}\Lambda_{\Phi\circ R}(R(x^*))}\]
\begin{proof}
Let $y$ be analytical with respect to both $\Phi$ and $\Phi\circ R$; we then get that :
\[\langle\omega_{J_\Phi\Lambda_{\Phi}(x)}, y\rangle=\Phi(\sigma_{i/2}^{\Phi}(y)x^*x)\]
and, therefore :
\begin{eqnarray*}
\langle\omega_{J_\Phi\Lambda_{\Phi}(x)}\circ R, y\rangle
&=&\Phi(\sigma_{i/2}^{\Phi}(R(y))x^*x)\\
&=&\Phi(R(\sigma_{-i/2}^{\Phi\circ R}(y)x^*x)\\
&=&\Phi\circ R(R(x^*x)\sigma_{-i/2}^{\Phi\circ R}(y))\\
&=&\overline{\Phi\circ R(\sigma_{i/2}^{\Phi\circ R}(y^*)R(x)R(x^*))}
\end{eqnarray*}
which, by a similar calculation, is equal to $\overline{<\omega_{J_{\Phi\circ R}\Lambda_{\Phi\circ R}(R(x^*))}, y^*>}$; which gives the result. \end{proof}

\subsection{Lemma}
\label{lemA}
{\it Let $\gG=(N, M, \alpha, \beta, \Gamma, T, T', \nu)$ be a measured quantum groupoid, and let us suppose that the von Neuman algebra $N$ is abelian;  let us use the notations of \ref{propmanag}(iv) and consider the ${\bf C}^*$- subalgebra $A_n(W)$ of $M$ (\ref{AW}, \ref{propmanag}(iv)); for any $x$ in $\gN_T\cap\gN_\Phi$, there exists $x_n$ in $A_n(W)\cap\gM_T\cap\gM_\Phi$ such that $\Lambda_T(x_n)$ is norm converging to $\Lambda_T(x)$. }

\begin{proof}
As $\alpha(N)\subset M(A_n(W))$ (\ref{propmanag}(iv)), we get that $T(A_n(W)\cap\gM_T)$ is an ideal of $\alpha (N)$, which, by normality of $T$, is weakly dense in $N$; let $e_n$ be a countable approximate unit of $T(A_n(W)\cap\gM_T)$; we have $e_nT(x^*x)=T(x^*e_nx)$, which is increasing to $T(x^*x)$, and, using Dini's theorem, is therefore norm converging to $T(x^*x)$. Let $f_n$ positive in $A_n(W)$ such that $e_n=T(f_n)$; we have $e_nT(x^*x)=T(T(x^*x)^{1/2}f_nT(x^*x)^{1/2})=T(x_n^*x_n)$, where $x_n=f_n^{1/2}T(x^*x)^{1/2}$ belongs to $A_n(W)\cap\gM_T$. We then get that $\Lambda_T(x_n)$ is norm converging to $\Lambda_T(x)$. \end{proof}

\subsection{Theorem}
\label{thgamma}
{\it Let $\gG=(N, M, \alpha, \beta, \Gamma, T, T', \nu)$ be a measured quantum groupoid, and let us suppose that the von Neuman algebra $N$ is abelian;  let us use the notations of \ref{propmanag}(iv) and consider the ${\bf C}^*$-subalgebra $A_n(W)$ of $M$ (\ref{AW}, \ref{propmanag}(iv)); then, for all $x_1$, $x_2$ in $A_n(W)\cap\gN_T\cap\gN_\Phi$, $y_1$, $y_2$ in $A_n(W)\cap\gN_{R\circ T\circ R}\cap\gN_{\Phi\circ R}$ : 
\newline
(i) we have :
\[(id\underset{N}{_\beta *_\alpha}\omega_{J_\Phi\Lambda_\Phi(x_1), J_\Phi\Lambda_\Phi(x_2)})\Gamma(A_n(W))\subset A_n(W)\]
and the closed linear set generated by all elements of the form :
\[(id\underset{N}{_\beta *_\alpha}\omega_{J_\Phi\Lambda_\Phi(x_1), J_\Phi\Lambda_\Phi(x_2)})\Gamma(x)\]
where $x$ is in $A_n(W)$, $x_1$, $x_2$ in $A_n(W)\cap\gN_T\cap\gN_\Phi$,  is equal to $A_n(W)$. 
\newline
(ii) we have :
\[(\omega_{J_{\Phi\circ R}\Lambda_{\Phi\circ R}(y_1), \Lambda_{\Phi\circ R}(y_2)}\underset{N}{_\beta *_\alpha}id)\Gamma(A_n(W))\subset A_n(W)\]
and the closed linear set generated by all elements of the form :
\[(\omega_{J_{\Phi\circ R}\Lambda_{\Phi\circ R}(y_1), \Lambda_{\Phi\circ R}(y_2)}\underset{N}{_\beta *_\alpha}id)\Gamma(y)\]
 where $y$ is in $A_n(W)$, $y_1$, $y_2$ in $A_n(W)\cap\gN_{R\circ T\circ R}\cap\gN_{\Phi\circ R}$, is equal to $A_n(W)$. }

\begin{proof}
Let us take $x$, $x_1$, $x_2$ in $\gN_T\cap\gN_\Phi$; we have, by  \ref{thL1}(ii) :
\[(id\underset{N}{_\beta *_\alpha}\omega_{J_\Phi\Lambda_\Phi(x_1), J_\Phi\Lambda_\Phi(x_2)})\Gamma(x^*)=(id*\omega_{J_\Phi\Lambda_\Phi(x_1^*x_2), \Lambda_\Phi(x)})(W)\]
If $x$ is in $A_n(W)\cap\gN_T\cap\gN_\Phi$; by the norm density of  $A_n(W)\cap\gN_T\cap\gN_\Phi$ into $A_n(W)$, we get that, for any $y$ in $A_n(W)$, $(id\underset{N}{_\beta *_\alpha}\omega_{J_\Phi\Lambda_\Phi(x_1), J_\Phi\Lambda_\Phi(x_2)})\Gamma(y)$ belongs to $A_n(W)$, from which we get the first result of (i). 
\newline
 Using \ref{thL1}(ii), we get that the first closed linear set contains all elements of the form 
$(id*\omega_{J_\Phi\Lambda_\Phi(x_1^*x_2), \Lambda_\Phi(x)})(W)$, where $x$, $x_1$, $x_2$ are in $A_n(W)\cap\gN_T\cap\gN_\Phi$, and, by linearity, all elements of the form $(id*\omega_{J_\Phi\Lambda_\Phi(y), \Lambda_\Phi(x)})(W)$, where $x$ is in $A_n(W)\cap\gN_T\cap\gN_\Phi$ and $y$ is in $A_n(W)\cap\gM_T\cap\gM_\Phi$; using then \ref{lemA}, we get it contains all elements of the form $(id*\omega_{J_\Phi\Lambda_\Phi(y), \Lambda_\Phi(x)})(W)$, where $x$, $y$ belong to $\gN_T\cap\gN_\Phi$; so, by prop \ref{ovw}, it contains all elements of the form $(i*\omega_{\xi, \eta})(W)$, where $\xi$ is in $D(_\alpha H, \nu)$ and $\eta$ is in $D(H_{\hat{\beta}}, \nu^o)$. Therefore, it contains $A_n(W)$, and, by the first result of (i), it is equal to $A_n(W)$, which finishes the proof of (i). 
\newline
We have now, using \ref{lemomega} :
\begin{multline*}
(\omega_{J_{\Phi\circ R}\Lambda_{\Phi\circ R}(x_1), \Lambda_{\Phi\circ R}(x_2)}\underset{N}{_\beta *_\alpha}id)\Gamma(x)\\
=(\omega_{J_\Phi\Lambda_\Phi(R(x_2^*), J_\Phi\Lambda_\Phi(R(x_1^*))}\circ R\underset{N}{_\beta *_\alpha}id)\Gamma(x)\\
=R((id\underset{N}{_\beta *_\alpha}\omega_{J_\Phi\Lambda_\Phi(R(x_2^*)), J_\Phi\Lambda_\Phi(R(x_1^*))})\Gamma(R(x))
\end{multline*}
which gives (ii). \end{proof}

\section{Measured quantum groupoids with a central basis}
\label{central}
We deal now with a measured quantum groupoid $\gG=(N, M, \alpha, \beta, \Gamma, T, T', \nu)$, such that the von Neuman algebra $\alpha(N)$ is included into the center $Z(M)$. Then, we obtain first some results about the restrictions of $\Phi$, $\Phi\circ R$, $T$ and $RTR$ to the ${\bf C}^*$-algebra $A_n(W)$ (\ref{Phi2}), and a Plancherel-like formula for the coproduct $\Gamma$ (\ref{thcentral1}), which gives that the coproduct sends the ${\bf C}^*$-algebra $A_n(W)$ in the multiplier algebra of the ${\bf C}^*$-algebra $A_n(W)\underset{N}{_\beta\otimes_\alpha}A_n(W)$. A summary of all these properties of $A_n(W)$ is given in \ref{thresum}.

\subsection{Lemma}
\label{propcentral}

{\it Let $\gG=(N, M, \alpha, \beta, \Gamma, T, T', \nu)$ be a measured quantum groupoid; then, are equivalent :
\newline
(i) the von Neuman algebra $\alpha(N)$ is included into the center $Z(M)$; 
\newline
(ii) the von Neuman algebra $\beta(N)$ is included into the center $Z(M)$;
\newline
(iii) the representation $\hat{\beta}$ is equal to $\alpha$.}

\begin{proof}
As $\beta=R\circ\alpha$, with $R$ an anti-$*$-isomomorphism of $M$, we see trivially that (i) and (ii) are equivalent. Moreover, as, by definition, we have $\hat{\beta}(n)=J_\Phi\alpha(n^*)J_\Phi$, where $J_\phi$ is the canonical antilinear bijective and involutive isometry on $H_\Phi$ constructed by the Tomita-Takesaki theory asociated to the weight $\Phi=\nu\circ\alpha^{-1}\circ T$ on $M$, we get that (i) and (iii) are equivalent. 
\end{proof}

\subsection{Lemma}
\label{lemLT}
{\it Let $\gG=(N, M, \alpha, \beta, \Gamma, T, T', \nu)$ be a measured quantum groupoid; then :
(i) $\Lambda_\Phi(\gN_\Phi\cap\gN_T\cap\gN_{\Phi\circ R}\cap\gN_{RTR})$ is dense in $H_\Phi$ (\cite{L2}, 6.5), and the subset of elements $x$ in $\gM_\Phi\cap\gM_T\cap\gM_{\Phi\circ R}\cap\gM_{RTR}$ which are analytic with respect both $\Phi$ and $\Phi\circ R$, and such that $\sigma_z^\Phi\circ\sigma_{z'}^{\Phi\circ R}(x)$ belongs to $\gM_\Phi\cap\gM_T\cap\gM_{\Phi\circ R}\cap\gM_{RTR}$, for all $z,z'\in {\bf C}$, is a $*$-algebra dense in $M$ (\cite{L2}, 6.6).
\newline
(ii) let us suppose that $\alpha(N)$ is central in $M$; then, for any $x\in\gN_\Phi\cap\gN_T$, there exists $x_n$ in $\gN_\Phi\cap\gN_T\cap\gN_{\Phi\circ R}\cap\gN_{RTR}$ such that $\Lambda_T(x_n)$ is norm converging to $\Lambda_T(x)$. }

\begin{proof} Thanks to (\cite{L2}, 6.6), let's take $h_n$ increasing to $1$, with $h_n$ in $\gM_\Phi\cap\gM_T\cap\gM_{\Phi\circ R}\cap\gM_{RTR}$, $h_n$ analytic with respect both $\Phi$ and $\Phi\circ R$, and such that $\sigma_z^\Phi\circ\sigma_{z'}^{\Phi\circ R}(h_n)$ belongs to $\gM_\Phi\cap\gM_T\cap\gM_{\Phi\circ R}\cap\gM_{RTR}$, for all $z,z'\in {\bf C}$. For any $x$ in $\gN_\Phi$, we get that $x_n=x\sigma_{-i/2}^\Phi(h_n)$ belongs to $\gN_\Phi\cap\gN_T\cap\gN_{\Phi\circ R}\cap\gN_{RTR}$, and we get that $\Lambda_\Phi(x\sigma_{-i/2}^\Phi(h_n))=J_\Phi h_nJ_\Phi\Lambda_\Phi(x)$ is converging to $\Lambda_\Phi(x)$ (which gives an alternate proof of \cite{L2}, 6.5, mostly inspired from the initial one, that we shall use in the sequel of the proof of this lemma). Let us suppose now that $\alpha(N)$ is central in $M$, and let's take $x\in\gN_\Phi\cap\gN_T$, and $p\in\gN_\nu$; we have :
\[\Lambda_T(x_n)\Lambda_\nu(p)=\Lambda_\Phi(x\sigma_{-i/2}^\Phi(h_n)\alpha(p))=\Lambda_\Phi(x\alpha(p)\sigma_{-i/2}^\Phi(h_n))=J_\Phi h_nJ_\Phi\Lambda_T(x)\Lambda_\nu(p)\]
from which we get, by continuity, that $\Lambda_T(x_n)=J_\Phi h_nJ_\Phi\Lambda_T(x)$, and that $\Lambda_T(x_n)$ is weakly converging to $\Lambda_T(x)$; more precisely, we get that $T(x_n^*x_n)=\Lambda_T(x)J_\Phi h_n^2J_\Phi\Lambda_T(x)$ is increasing to $T(x^*x)$; if we write $X$ for the spectrum of the ${\bf C}^*$-algebra generated by $T(\gM_T)$, using Dini's theorem in $C_0(X)$, we get that $T(x_n^*x_n)$ is norm converging to $T(x^*x)$; more precisely, as :
\[\|\Lambda_T(x_n)-\Lambda_T(x)\|^2=\|T(x_n^*x_n)-T(x^*x_n)+T(x_n^*x)-T(x^*x)\|\]
and $T(x^*x_n)=\Lambda_T(x)^*J_\Phi h_nJ_\Phi\Lambda_T(x)=T(x_n^*x)$, we get the result. 
\end{proof}

\subsection{Theorem}
\label{Phi2}
{\it Let $\gG=(N, M, \alpha, \beta, \Gamma, T, T', \nu)$ be a measured quantum groupoid and let us suppose that $\alpha(N)$ is central in $M$; then the restrictions of $\Phi$ and $\Phi\circ R$ to the ${\bf C}^*$-algebra $A_n(W)$ are faithful lower semi-continuous densely defined KMS weights. Moreover, the restrictions of $T$ and $RTR$ to $A_n(W)$ are densely defined. }

\begin{proof}
We have got in \ref{Phi}(iii) that, for any $x,y$ in $\gM_\Phi\cap\gM_T\cap\gM_{\Phi\circ R}\cap\gM_{RTR}$, the operator $(i*\omega_{J_\Phi\Lambda_\Phi(y), \Lambda_\Phi(x)})(W)$ belongs to $\gM_T\cap\gM_\Phi\cap\gM_{\Phi\circ R}\cap\gM_{RTR}$; using now \ref{lemLT} and \ref{ovw}, we see that this set of operators are norm dense in the set of operators of the form $(i*\omega_{J_\Phi\Lambda_\Phi(y'), \Lambda_\Phi(x')})(W)$, for all $x'$, $y'$ in $\gN_\Phi\cap\gN_T$; using again \ref{ovw}, we see that it is norm dense again in the set of operators of the form $(i*\omega_{\xi, \eta})(W)$, for $\xi, \eta\in D(_\alpha H_\Phi, \nu)$, and, therefore, by definition, in $A_n(W)$. Moreover, by \ref{corPhi}(ii) and (iii), we get that the modular groups $\sigma_t^\Phi$ and $\sigma_t^{\Phi\circ R}$ are norm continuous on $A_n(W)$. \end{proof}

\subsection{Lemma}
\label{lemcentral}
{\it In the situation of \ref{propcentral}, 
let $(e_i)_{i\in I}$ be an $(\alpha, \nu)$-orthogonal basis of $H$; then, we have : 
\newline
(i) for all $\xi$, $\eta_2$ in $D(_\alpha H, \nu)$, and $\eta_1$ in $D(_\alpha H, \nu)\cap D(H_\beta, \nu)$ :
\[(\omega_{\eta_1, \eta_2}*id)(W)\xi=\sum_i\alpha(\langle (id*\omega_{\xi, e_i})(W)\eta_1, \eta_2\rangle_{\alpha, \nu})e_i\]
(ii) for all $\xi_1$ in $D(H_\beta, \nu)$, $\xi_2$ in $D(_\alpha H, \nu)\cap D(H_\beta, \nu)$ and $\eta$ in $D(_\alpha H, \nu)$ 
\[(\omega_{\xi_1, \xi_2}*id)(W)^*\eta=\sum_i\alpha(\langle\xi_2, (id*\omega_{e_i, \eta})(W)\xi_1\rangle_{\beta, \nu})e_i\]}
\begin{proof}
We have :
\[W(\eta_1\underset{\nu}{_\beta\otimes_\alpha}\xi)=\sum_i(id*\omega_{\xi, e_i})(W)\eta_1\underset{\nu}{_\alpha\otimes_\alpha}e_i\]
Let now $\zeta$ be in $H$; we have then :
\begin{eqnarray*}
((\omega_{\eta_1, \eta_2}*id)(W)\xi |\zeta)&=&
\sum_i((id*\omega_{\xi, e_i})(W)\eta_1\underset{\nu}{_\alpha\otimes_\alpha}e_i|\eta_2\underset{\nu}{_\alpha\otimes_\alpha}\zeta)\\
&=&(\sum_i\alpha(\langle (id*\omega_{\xi, e_i})(W)\eta_1, \eta_2\rangle_{\alpha, \nu})e_i|\zeta)
\end{eqnarray*}
from which we get (i).
\newline
We have :
\[W^*(\xi_2\underset{\nu}{_\alpha\otimes_\alpha}\eta)=\sum_i(id*\omega_{e_i, \eta})(W)^*\xi_2\underset{\nu}{_\beta\otimes_\alpha}e_i\]
Let now $\zeta$ be in $H$; we have then :
\begin{eqnarray*}
((\omega_{\xi_1, \xi_2}*id)(W)^*\eta|\zeta)
&=&(\sum_i(id*\omega_{e_i, \eta})(W)^*\xi_2\underset{\nu}{_\beta\otimes_\alpha}e_i|\xi_1\underset{\nu}{_\beta\otimes_\alpha}\zeta)\\
&=&(\sum_i\alpha(\langle\xi_2, (id*\omega_{e_i, \eta})(W)\xi_1\rangle_{\beta, \nu})e_i|\zeta)
\end{eqnarray*}
which finishes the proof. \end{proof}

\subsection{Lemma}
\label{lemp}
{\it In the situation of \ref{propcentral}, let $(e_i)_{i\in I}$ be an $(\alpha, \nu)$-orthogonal basis of $H$ and $J$ a finite subset of $I$; let us write $p_J=\Sigma_{i\in J}\theta^{\alpha, \nu}(e_i, e_i)$; then, for all $\Xi_1$, $\Xi_2$ in $H\underset{\nu}{_\beta\otimes_\alpha}H$, the finite sum :
\[\sum_{i\in J} ((id*\omega_{e_i, \eta})(W)\underset{N}{_\beta\otimes_\alpha}(id*\omega_{\xi, e_i})(W))\Xi_1|
\Xi_2)\]
is equal to :
\[((\sigma_{\nu}\underset{N}{_\alpha\otimes_\alpha}1_{\gH})
(1_{\gH}\underset{N}{_\alpha\otimes_\alpha}W)
\sigma^{1,2}_{\alpha, \alpha}
(1_{\gH}\underset{N}{_\beta\otimes_\alpha}(1\underset{N}{_\alpha\otimes_\alpha}p_J)W)(\Xi_1
\underset{\nu}{_\beta\otimes_\alpha}\xi)|
\Xi_2\underset{\nu}{_\alpha\otimes_\alpha}\eta)\]}

\begin{proof}
Let $\xi_1$ be in $D(H_\beta, \nu)$, $\xi_2$ in $D(_\alpha H, \nu)\cap D(H_\beta, \nu)$, $\eta_1$ in $D(_\alpha H, \nu)\cap D(H_\beta, \nu)$ and $\eta_2$ in $D(_\alpha H, \nu)$. If we take $\Xi_1=\xi_1\underset{\nu}{_\beta\otimes_\alpha}\eta_1$ and $\Xi_2=\xi_2\underset{\nu}{_\beta\otimes_\alpha}\eta_2$, the scalar product we are dealing with is equal to :
\[\sum_{i\in J}(\alpha(\langle(id*\omega_{e_i, \eta})(W)\xi_1, \xi_2\rangle_{\beta, \nu})(id*\omega_{\xi, e_i})(W)\eta_1|\eta_2)\]
Using the commutativity of $N$, we consider $\alpha(\langle e_i, e_i\rangle_{\alpha, \nu})(id*\omega_{\xi, e_i})(W)$, which is equal to $(id*\omega_{\xi, \alpha(\langle e_i, e_i\rangle_{\alpha, \nu})e_i})(W)$, thanks to the commutation relations of $W$, and the fact that $\alpha=\hat{\beta}$. But by \ref{spatial}, we know that $\alpha(\langle e_i, e_i\rangle_{\alpha, \nu})e_i=e_i$, and therefore : 
\[\alpha(\langle e_i, e_i\rangle_{\alpha, \nu})(id*\omega_{\xi, e_i})(W)=(id*\omega_{\xi, e_i})(W)\]
So, this scalar product is equal to :
\[\sum_{i\in J}(\alpha(\langle (id*\omega_{e_i, \eta})(W)\xi_1, \xi_2\rangle_{\beta, \nu})(id*\omega_{\xi, e_i})(W)\eta_1\underset{\nu}{_\alpha\otimes_\alpha}e_i|\eta_2\underset{\nu}{_\alpha\otimes_\alpha}e_i)\]
which is :
\[\sum_{i\in J}(\alpha(\langle(id*\omega_{e_i, \eta})(W)\xi_1, \xi_2\rangle_{\beta, \nu}\langle(id*\omega_{\xi, e_i})(W)\eta_1, \eta_2\rangle_{\alpha, \nu})e_i|e_i)\]
This is equal to :
\[\sum_{i\in J}(\alpha(\langle (id*\omega_{\xi, e_i})(W)\eta_1, \eta_2\rangle_{\alpha, \nu})e_i
|\alpha(\langle\xi_2, (id*\omega_{e_i, \eta})(W)\xi_1\rangle_{\beta, \nu})e_i)\]
which, thanks to \ref{lemcentral}(i) and (ii), is equal to :
\[((\omega_{\xi_1, \xi_2}*id)(W)p_J(\omega_{\eta_1, \eta_2}*id)(W)\xi|\eta)\]
Coming back to the calculations made in \ref{prop2Gamma}, we get it is equal to :
\[(W(\eta_1\underset{\nu}{_\beta\otimes_\alpha}\xi)|
\eta_2\underset{\nu}{_\alpha\otimes_\alpha}p_J(\omega_{\xi_1, \xi_2}*id)(W)^*\eta)\]
Defining now $\zeta_i$, $\zeta'_i$ as in \ref{lem1}, we get that it is equal to :
\[(W(\eta_1\underset{\nu}{_\beta\otimes_\alpha}\xi)|\eta_2\underset{\nu}{_\alpha\otimes_\alpha}\sum_i\alpha(\langle\zeta_i, \xi_1\rangle_{\beta, \nu})p_J\zeta'_i)\]
Using \ref{lem2}, it is equal to :
\[(\sigma^{1,2}_{\alpha, \alpha}
(1_{\gH}\underset{N}{_\beta\otimes_\alpha}W)(\xi_1\underset{\nu}{_\beta\otimes_\alpha}\eta_1\underset{\nu}{_\beta\otimes_\alpha}\xi)|\eta_2\underset{\nu}{_\alpha\otimes_\alpha}(1\underset{N}{_\alpha\otimes_\alpha}p_J)W^*(\xi_2\underset{\nu}{_\alpha\otimes_\alpha}\eta))\]
which is equal to :
\[(\sigma^{1,2}_{\alpha, \alpha}
(1_{\gH}\underset{N}{_\beta\otimes_\alpha}(1\underset{N}{_\alpha\otimes_\alpha}p_J)W)(\xi_1\underset{\nu}{_\beta\otimes_\alpha}\eta_1\underset{\nu}{_\beta\otimes_\alpha}\xi)|\eta_2\underset{\nu}{_\alpha\otimes_\alpha}W^*(\xi_2\underset{\nu}{_\alpha\otimes_\alpha}\eta))\]
from which we get the result, by linearity, continuity and density. 
\end{proof}

\subsection{Proposition}
\label{propc}
{\it Let $(N, M, \alpha, \beta, \Gamma, T, T', \nu)$ be a measured quantum groupoid, and let us suppose that the von Neuman algebra $\alpha(N)$ is included into the center $Z(M)$;
let $(e_i)_{i\in I}$ be an $(\alpha, \nu)$-orthogonal basis of $H$ and $J$ a finite subset of $I$; let us write $p_J=\Sigma_{i\in J}\theta^{\alpha, \nu}(e_i, e_i)$; let $k_1$, $k_2$ in $\mathcal{K}_{\alpha, \nu}$, $\xi$ in $D(_\alpha H, \nu)$, $\eta$ in $H$; then, we have :
\[lim_J\|(k_2\underset{N}{_\alpha\otimes_\alpha}(1-p_J))W(k_1\eta\underset{\nu}{_\beta\otimes_\alpha}\xi)\|=0\]}

\begin{proof}
Let $\eta_1$ in $D(_\alpha H, \nu)\cap D(H_\beta, \nu)$ and $\eta_2$ in $D(_\alpha H, \nu)$; we have :
\[R^{\alpha, \nu}(p_J(\omega_{\eta_1, \eta_2}*id)(W)\xi)=p_J(\omega_{\eta_1, \eta_2}*id)(W)R^{\alpha, \nu}(\xi)\]
and, therefore :
\begin{multline*}
\langle p_J(\omega_{\eta_1, \eta_2}*id)(W)\xi, p_J(\omega_{\eta_1, \eta_2}*id)(W)\xi\rangle_{\alpha, \nu}=\\
R^{\alpha, \nu}(\xi)^*(\omega_{\eta_1, \eta_2}*id)(W)^*p_J(\omega_{\eta_1, \eta_2}*id)(W)R^{\alpha, \nu}(\xi)
\end{multline*}
which is increasing with $J$ towards 
\[\langle (\omega_{\eta_1, \eta_2}*id)(W)\xi, (\omega_{\eta_1, \eta_2}*id)(W)\xi\rangle_{\alpha, \nu}\]
Let $X$ be the spectrum of ${\bf C}^*(\nu)$, and let us identify ${\bf C}^*(\nu)$ to $C_0(X)$;  using then Dini's theorem, we get it is norm converging, from which we infer that :
\[lim_J\|R^{\alpha, \nu}((1-p_J)(\omega_{\eta_1, \eta_2}*id)(W)\xi)\|=0\]
But, by \ref{lemcentral}(i), we have :
\[(1-p_J)(\omega_{\eta_1, \eta_2}*id)(W)\xi=\Sigma_{i\notin J}\alpha(\langle (id*\omega_{\xi, e_i})(W)\eta_1, \eta_2\rangle_{\alpha, \nu})e_i\]
and, therefore :
\[R^{\alpha, \nu}((1-p_J)(\omega_{\eta_1, \eta_2}*id)(W)\xi)=\Sigma_{i\notin J}R^{\alpha, \nu}(e_i)\langle (id*\omega_{\xi, e_i})(W)\eta_1, \eta_2\rangle_{\alpha, \nu}\]
and :
\begin{multline*}
\|R^{\alpha, \nu}((1-p_J)(\omega_{\eta_1, \eta_2}*id)(W)\xi)\|^2=\\
\|\Sigma_{i\notin J}\langle (id*\omega_{\xi, e_i})(W)\eta_1, \eta_2\rangle_{\alpha, \nu}^*\langle (id*\omega_{\xi, e_i})(W)\eta_1, \eta_2\rangle_{\alpha, \nu}\|
\end{multline*}
We have :
\begin{align*}
\langle (id*\omega_{\xi, e_i})(W)\eta_1, \eta_2\rangle_{\alpha, \nu}
&=
R^{\alpha, \nu}(\eta_2)^*(\rho^{\alpha, \alpha}_{e_i})^*W\rho^{\beta, \alpha}_\xi R^{\alpha, \nu}(\eta_1)\\
&=\
(\rho^{\alpha, \alpha}_{e_i})^*(R^{\alpha, \nu}(\eta_2)^*\underset{N}{_\alpha\otimes_\alpha}1)W\rho_\xi^{\beta, \alpha} R^{\alpha, \nu}(\eta_1)
\end{align*}
and, therefore :
\begin{multline*}
\Sigma_{i\notin J}\langle (id*\omega_{\xi, e_i})(W)\eta_1, \eta_2\rangle_{\alpha, \nu}^*\langle (id*\omega_{\xi, e_i})(W)\eta_1, \eta_2\rangle_{\alpha, \nu}=\\
R^{\alpha, \nu}(\eta_1)^*(\rho_\xi^{\beta, \alpha})^*W^*(\theta^{\alpha, \nu}(\eta_2, \eta_2)\underset{N}{_\alpha\otimes_\alpha}(1-p_J))W\rho_\xi^{\beta, \alpha} R^{\alpha, \nu}(\eta_1)
\end{multline*}
and its norm is equal to :
\[\|(\theta^{\alpha, \nu}(\eta_2, \eta_2)\underset{N}{_\alpha\otimes_\alpha}(1-p_J))W\rho_\xi^{\beta, \alpha} R^{\alpha, \nu}(\eta_1)\|^2\]
So, we have that :
\[lim_J\|(\theta^{\alpha, \nu}(\eta_2, \eta_2)\underset{N}{_\alpha\otimes_\alpha}(1-p_J))W\rho_\xi^{\beta, \alpha} R^{\alpha, \nu}(\eta_1)\|=0\]
and, therefore :
\[lim_J\|(\theta^{\alpha, \nu}(\eta_2, \eta_2)\underset{N}{_\alpha\otimes_\alpha}(1-p_J))W\rho_\xi^{\beta, \alpha} \theta^{\alpha, \nu}(\eta_1, \eta_1)\|=0\]
and we get :
\[lim_J\|(\theta^{\alpha, \nu}(\eta_2, \eta_2)\underset{N}{_\alpha\otimes_\alpha}(1-p_J))W( \theta^{\alpha, \nu}(\eta_1, \eta_1)\eta\underset{\nu}{_\beta\otimes_\alpha}\xi)\|=0\]
from which we get the result. \end{proof}

\subsection{Theorem}
\label{thcentral1}
{\it Let $(N, M, \alpha, \beta, \Gamma, T, T', \nu)$ be a measured quantum groupoid, and let us suppose that the von Neuman algebra $\alpha(N)$ is included into the center $Z(M)$;
let $(e_i)_{i\in I}$ be an $(\alpha, \nu)$-orthogonal basis of $H$; then, we have, for all $\xi$, $\eta$ in $D(_\alpha H, \nu)$ 
\[\Gamma ((id*\omega_{\xi, \eta})(W))=\sum_i (id*\omega_{e_i, \eta})(W)\underset{N}{_\beta\otimes_\alpha}(id*\omega_{\xi, e_i})(W)\]
where the sum is weakly and strictly convergent. 
\newline
Therefore, using the strict convergence, we get that 
$\Gamma(A_n(W))$ is included into the multiplier algebra of the ${\bf C}^*$-algebra $A_n(W)\underset{N}{_\beta\otimes_\alpha}A_n(W)$.}
\begin{proof}
Let $\xi$, $\eta$ in $D(_\alpha H, \nu)$. Using \ref{lemp}, for all finite $J\subset I$, we have :
\[\|\sum_{i\in J} ((id*\omega_{e_i, \eta})(W)\underset{N}{_\beta\otimes_\alpha}(id*\omega_{\xi, e_i})(W))\|
\leq \|R^{\alpha, \nu}(\xi)\|\|R^{\alpha, \nu}(\eta)\|\]
Let $\xi_1$ be in $D(H_\beta, \nu)$, $\xi_2$ in $D(_\alpha H, \nu)\cap D(H_\beta, \nu)$, $\eta_1$ in $D(_\alpha H, \nu)\cap D(H_\beta, \nu)$ and $\eta_2$ in $D(_\alpha H, \nu)$; using \ref{prop2Gamma} and \ref{lemp}, we have :
\[(\Gamma((id*\omega_{\xi, \eta})(W))(\xi_1\underset{\nu}{_\beta\otimes_\alpha}\eta_1)|
\xi_2\underset{\nu}{_\beta\otimes_\alpha}\eta_2)
=
((\omega_{\xi_1, \xi_2}*id)(W)(\omega_{\eta_1, \eta_2}*id)(W)\xi|\eta)\]
is equal to :
\[\sum_i ((id*\omega_{e_i, \eta})(W)\xi_1\underset{\nu}{_\beta\otimes_\alpha}(id*\omega_{\xi, e_i})(W)\eta_1|
\xi_2\underset{\nu}{_\beta\otimes_\alpha}\eta_2)\]
If we apply \ref{ovw} to the inclusion $\alpha(N)\subset \widehat{M}$ and the operator-valued weight $\hat{T}$, we get that $D(_\alpha H_\Phi, \nu)\cap D((H_\Phi)_\beta, \nu)$ is dense in $H_\Phi$, and we obtain the weak convergence of the sum $\sum_i((id*\omega_{e_i, \eta})(W)\underset{\nu}{_\beta\otimes_\alpha}(id*\omega_{\xi, e_i})(W)$ to $\Gamma((id*\omega_{\xi, \eta})(W))$.
\newline
Moreover, we get, using \ref{lemp} that, for any $k_1$, $k_2$ in $\mathcal K_{\alpha, \nu}$ :
\begin{multline*}
|(\Sigma_{i\notin J}(id*\omega_{e_i, \eta})(W)\xi_1\underset{\nu}{_\beta\otimes_\alpha}k_1^*(id*\omega_{\xi, e_i})(W)k_2\eta_1|\xi_2\underset{\nu}{_\beta\otimes_\alpha}\eta_2)|\leq\\
\|(k_2\underset{N}{_\alpha\otimes_\alpha}(1-p_J))W(k_1\eta\underset{\nu}{_\beta\otimes_\alpha}\xi)\|
\|\xi_1\underset{\nu}{_\beta\otimes_\alpha}\eta_1\|\|\xi_2\underset{\nu}{_\beta\otimes_\alpha}\eta_2\|
\end{multline*}
from which, thanks to \ref{propc}, we get that :
\[lim_J\|\sum_{i\notin J}(id*\omega_{e_i, \eta})(W)\underset{N}{_\beta\otimes_\alpha}k_1^*(id*\omega_{\xi, e_i})(W)k_2\|=0\]
Let now $y_1$, $y_2$, $y_3$, $y_4$ in $A_n(W)$, and $\epsilon$ positive; as $A_n(W)\subset \alpha(N)'=M(\mathcal K_{\alpha, \nu})$, there exists $k\in\mathcal K_{\alpha, \nu}$ such that $\|y_1k-y_1\|\leq \epsilon$, and $\|y_2k-y_2\|\leq\epsilon$. Moreover, there exists a finite subset $J\subset I$ such that :
\[\|\sum_{i\notin J}(id*\omega_{e_i, \eta})(W)\underset{N}{_\beta\otimes_\alpha}k^*y_1^*(id*\omega_{\xi, e_i})(W)y_2k\|\leq\epsilon\]
As, for any finite $J'$ such that $J\cap J'=\emptyset$, we have proved that :
\[\|\sum_{i\in J'} (id*\omega_{e_i, \eta})(W)\underset{N}{_\beta\otimes_\alpha}(id*\omega_{\xi, e_i})(W)\|
\leq \|R^{\alpha, \nu}(\xi)\|\|R^{\alpha, \nu}(\eta)\|\]
and, as $\sum_{i\in J'} y_3^*(id*\omega_{e_i, \eta})(W)y_4\underset{N}{_\beta\otimes_\alpha}y_1^*(id*\omega_{\xi, e_i})(W)y_2$ is equal to :
\begin{multline*}
[y_3^*\underset{N}{_\beta\otimes_\alpha}(y_1^*-k^*y_1^*)][\sum_{i\in J'}(id*\omega_{e_i, \eta})(W)\underset{N}{_\beta\otimes_\alpha}(id*\omega_{\xi, e_i})(W)][y_4\underset{N}{_\beta\otimes_\alpha}y_2]+\\
[y_3^*\underset{N}{_\beta\otimes_\alpha}k^*y_1^*][\sum_{i\in J'}(id*\omega_{e_i, \eta})(W)\underset{N}{_\beta\otimes_\alpha}(id*\omega_{\xi, e_i})(W)][y_4\underset{N}{_\beta\otimes_\alpha}(y_2-ky_2)]+\\
[y_3^*\underset{N}{_\beta\otimes_\alpha}k^*y_1^*][\sum_{i\in J'}(id*\omega_{e_i, \eta})(W)\underset{N}{_\beta\otimes_\alpha}(id*\omega_{\xi, e_i})(W)][y_4\underset{N}{_\beta\otimes_\alpha}ky_2]
\end{multline*}
we get that :
\begin{multline*}
\|\sum_{i\in J'} y_3^*(id*\omega_{e_i, \eta})(W)y_4\underset{N}{_\beta\otimes_\alpha}y_1^*(id*\omega_{\xi, e_i})(W)y_2\|\leq\\
\|y_2\|\|y_3\|\|y_4\|\|R^{\alpha, \nu}(\xi)\|\|R^{\alpha, \nu}(\eta)\|\epsilon+\|y_1\|\|y_3\|\|y_4\|\|R^{\alpha, \nu}(\xi)\|\|R^{\alpha, \nu}(\eta)\|\epsilon+\|y_3\|\|y_4\|\epsilon
\end{multline*}
which gives the result.  \end{proof}

\subsection{Theorem}
\label{thresum}
{\it Let $\gG=(N, M, \alpha, \beta, \Gamma, T, T', \nu)$ be a measured quantum groupoid, such that the von Neuman algebra $\alpha(N)$ is included into the center $Z(M)$; then, the ${\bf C}^*$-algebra $A_n(W)$ bear the following properties :
\newline
(i) we have : $\alpha(N)\subset Z(M(A_n(W)))$, and $\beta(N)\subset Z(M(A_n(W)))$;
\newline
(ii) we have : $\Gamma(A_n(W))\subset M(A_n(W)\underset{N}{_\beta\otimes_\alpha}A_n(W))$;
\newline
(iii) $A_n(W)$ is globally invariant under the co-inverse $R$ and the scaling group $\tau_t$; moreover, the restriction of $\tau_t$ to $A_n(W)$ is a one-parameter norm continuous group of $*$-automorphisms of $A_n(W)$; 
\newline
(iv) the restrictions of $\Phi$ and $\Phi\circ R$ to $A_n(W)$ are faithful lower semi-continuous densely defined KMS weights on $A_n(W)$; the restrictions of $T$ and $RTR$ to $A_n(W)$ are densely defined. }

\begin{proof}
Result (i) had been obtained in \ref{propmanag}(ii), result (ii) in \ref{thcentral1}, result (iii) in \ref{corPhi}(i) and (iv), and result (iv) in \ref{Phi2}. \end{proof}

\section{Measured Quantum Groupoids with a central basis and 
continuous fields of ${\bf C}^*$-algebras}
\label{mcfield}

In this chapter, we go on with a measured quantum groupoid $(N, M, \alpha, \beta, \Gamma, T, T', \nu)$ such that $\alpha (N)$ is included in the center $Z(M)$; writing $X$ for the spectrum of the ${\bf C}^*$-algebra ${\bf C}^*(\nu)$ (which is the norm closure of $\gM_\nu$, and whose multiplier algebra is the von Neumann algebra $N$), we show that the restrictions of $\alpha$ and $\beta$ to $A_n(W)$ make, in two different ways, $A_n(W)$ be a $C_0(X)-{\bf C}^*$-algebra (\ref{continuousfield}(i)), and, more precisely, a continuous field of ${\bf C}^*$-algebras (\ref{continuousfield}(v)), because the restriction of $T$ to $A_n(W)$ gives a field of lower semi-continuous faithful weights $\varphi^x$ (\ref{continuousfield}(ii)), whose representations $\pi_{\varphi^x}$ form a continuous field of faithful representations of $A_n(W)$(\ref{continuousfield}(iv)). Moreover, the ${\bf C}^*$-algebra $A_n(W)\underset{N}{_\beta\otimes_\alpha}A_n(W)$ can be interpreted as the Blanchard's min tensor product $A_n(W)\underset{C_0(X)}{_\beta\otimes^m_\alpha}A_n(W)$ of these two fields of ${\bf C}^*$-algebras (\ref{cortensor}), and the restriction of $\Gamma$ to $A_n(W)$ sends therefore $A_n(W)$ into the multiplier algebra of this min tensor product (\ref{thcentral2}), which is here associative. All these results are summarized in \ref{thresume}.

\subsection{Theorem}
\label{continuousfield}
{\it Let $\gG=(N, M, \alpha, \beta, \Gamma, T, T', \nu)$ be a measured quantum groupoid, such that the von Neuman algebra $\alpha(N)$ is included into the center $Z(M)$, and let $X$ be the spectrum of ${\bf C}^*(\nu)$; we shall identify $\nu$ with a positive Radon measure on $X$, and $N$ with $L^\infty(X,\nu)=C_b(X)$ (by \ref{C*}, we have $N=M({\bf C}^*(\nu))$), and the positive extension of $N$ can be identified with lower semi-continuous functions on $X$, with values in $[0, +\infty]$. Then:
\newline
(i) thanks to the $*$-homomorphism $\alpha_{|C_0(X)}$ (resp. $\beta_{|C_0(X)}$), the ${\bf C}^*$-algebra $A_n(W)$ is a $C_0(X)-{\bf C}^*$-algebra, in the sense of Kasparov-Blanchard (\cite{Ka}, \cite{Bl1}).
\newline
(ii) the restriction of the weight $\Phi$ to the ${\bf C}^*$-algebra $A_n(W)$ can be disintegrated into a measurable field of lower semi-continuous faithful weights $\varphi^x$, invariant under $\sigma_t^\Phi$, satisfying the KMS conditions for $\sigma_t^\Phi$, and such that, for any $a\in A_n(W)^+$ :
\[\Phi(a)=\int_X\varphi^x(a) d\nu(x)\]
Moreover, we can identify $T(a)$ with the (image by $\alpha$ of the) function $x\mapsto \varphi^x(a)$, which is therefore lower semi-continuous (and bounded continuous is $a\in\gM_T^+$).
\newline
(iii) for any $f\in C_b(X)^+$ and $a\in A_n(W)^+$, we have $\varphi^x(\alpha (f)a)=f(x)\varphi^x(a)$, and $\varphi^x(a)=0$ if and only if $a\in\alpha(C_x(X))A_n(W)$.
\newline
(iv) the representations $\pi_{\varphi^x}$ form a continuous field of faithful representations of $A_n(W)$.
\newline
(v)  thanks to the $*$-homomorphism $\alpha_{|C_0(X)}$ (resp. $\beta_{|C_0(X)}$), the ${\bf C}^*$-algebra $A_n(W)$ is a continuous field over $X$ of ${\bf C}^*$-algebras.
\newline
(vi) we have :
\[H_\Phi=\int_X^{\oplus}H_{\varphi^x}d\nu(x)\]
\[M=\int_X^{\oplus}\pi_{\varphi^x}(A_n(W)/\alpha(C_x(X))A_n(W))"d\nu(x)\]
\[\Phi=\int_X^{\oplus}\overline{\varphi^x}d\nu(x)\]
where $\overline{\varphi^x}$ is the faithful semi-finite normal extension to $\pi_{\varphi^x}(A_n(W)/\alpha(C_x(X))A_n(W))"$ recalled in \ref{C*}, and :
\[T(\int_X^{\oplus}a^xd\nu(x))=(x\mapsto\overline{\varphi^x}(a^x))\]
where $a=\int_X^{\oplus}a^xd\nu(x)\in M_T^+$.
\newline
(vii) let $R$ be the coinverse of $\gG$; then $R_{|A_n(W)}$ can be disintegrated into a continuous field $R^x:A_n(W)/\alpha(C_x(X)A_n(W)\rightarrow A_n(W)/\beta(C_x(X)A_n(W))$, and we have :
\[H_{\Phi\circ R}=\int_X^{\oplus}H_{\varphi^x\circ R^x}d\nu(x)\]
\[M=\int_X^{\oplus}\pi_{\varphi^x\circ R^x}(A_n(W)/\beta(C_x(X)A_n(W))"d\nu(x)\]
\[\Phi\circ R=\int_X^{\oplus}\overline{\varphi^x\circ R^x}d\nu(x)\]
where $\overline{\varphi^x\circ R^x}$ is the faithful semi-finite normal extension to $\pi_{\varphi^x\circ R^x}(A_n(W)/\beta(C_x(X)A_n(W))"$, and :
\[RTR(\int_X^{\oplus}b^xd\nu(x))=(x\mapsto \overline{\varphi^x\circ R^x}(b^x))\]
where $b=\int_X^{\oplus}b^xd\nu(x)\in\gM_{RTR}^+$. 
 }
\begin{proof}
By \ref{propmanag}(iv), we get that $\alpha({\bf C}^*(\nu))$$\subset M(A_n(W))$, and with the hypothesis, we get that $\alpha ({\bf C}^*(\nu))$$\subset Z(M(A_n(W))$, which gives the first part of (i); the same holds if we take $\beta$ instead of $\alpha$, which finishes the proof of (i). The first part of (ii) is given by (\cite{T}, 4.11); moreover, the application which sends $Y\in M^+$ on the image under $\alpha$ of the function $x\mapsto \overline{\varphi^x}(Y)$ (with the notations of \ref{C*}) is a normal semi-finite operator-valued weight $T'$ from $M$ onto $\alpha(N)$, such that $\nu\circ\alpha^{-1}\circ T'=\Phi=\nu\circ\alpha^{-1}\circ T$, from which we infer that $T=T'$; taking now the restrictions to $A_n(W)^+$, we finish the proof of (ii). 
\newline
The first result of (iii) is just the operator-valued weight property of $T'$ discussed in the proof of (ii); let now $a\in A_n(W)^+$ such that $\varphi^x(a)=0$; then $T(a)$ is a the image under $\alpha$ of the lower semi-continuous function $x\mapsto\varphi^x(a)=f(x)$; and let us write $f_p=[inf (1, f)]^{1/p}$; then $f_p\in N=C_b(X)$, $\alpha(f_p)\in M(A_n(W))$, and $\alpha(f_p)A$ is included into $ \alpha(C_x(X))A_n(W)$; but $\alpha(f_p)$ is increasing  to $supp T(a)$ (in $\alpha(N)$); therefore, $\alpha(f_p)a$ is increasing to $a \times supp T(a)$, which is less than $a$; but, as :
\[T(a \times supp T(a))=T(a)\times supp T(a)=T(a)\]
using the faithfullness of $T$, we get that $a \times supp T(a)=a$, and, therefore, that $\alpha(f_p)a$ is increasing to $a$; let $B$ be the abelian ${\bf C}^*$-algebra generated by $\alpha(C_b(X))$ and $a$, and let $Y$ be the spectrum of $B$; then, we can identify $B$ with $C(Y)$, and, using Dini's theorem on $Y$, we get that $\alpha(f_p)a$ is norm converging to $a$, and, therefore, that $a$ belongs to $\alpha(C_x(X))A_n(W)$, which finishes the proof of (iii).
\newline
Let $a\in\gN_T$, $b\in M$, analytical with respect to $\Phi$; then, by (\cite{E6}, 2.2.2), $ab$ belongs to $\gN_T$, and $\Lambda_T(ab)=J_\Phi \sigma_{-i/2}^\Phi(b^*)J_\Phi\Lambda_T(a)$, from which we get :
\[T(b^*a^*ab)=\Lambda_T(a)^*J_\Phi\sigma_{-i/2}^\Phi(b^*)^*\sigma_{-i/2}^\Phi(b^*)J_\Phi\Lambda_T(a)\]
or :
\[T((\sigma_{-i/2}^\Phi(b^*)^*a^*a\sigma_{-i/2}^\Phi(b^*))=\Lambda_T(a)^*J_\Phi b^*bJ_\Phi\Lambda_T(a)\]
let us take now a family $b_k\in A_n(W)\cap\gM_T$ increasing to $1$ (which exists, thank to \ref{Phi}(i) and \ref{thgamma}(i)). Then, the elements :
\[c_k=\sqrt{1/\pi}\int_{-\infty}^{+\infty}e^{-t^2}\sigma_t^\Phi(b_k)dt\]
are in $A_n(W)\cap\gM_T$, analytical with respect to $\Phi$, and such that $\sigma_z^\Phi(c_k)$ belongs to $A_n(W)$ for any $z\in{\bf C}$, and the family $c_k$ is increasing to $1$. So, we get that the sequence :
\[T((\sigma_{-i/2}^\Phi(c_k^*)^*a^*a\sigma_{-i/2}^\Phi(c_k^*)^*)\]
 is increasing to $T(a^*a)$. By Dini's theorem (in $C_0(X)$), we get that it is norm converging, and therefore, we have, for all $x\in X$ :
 \[lim_k\varphi^x((\sigma_{-i/2}^\Phi(c_k^*)^*a^*a\sigma_{-i/2}^\Phi(c_k^*)^*)=\varphi^x(a^*a)\]
and, therefore, for any $x\in X$ :
\[lim_k\|\pi_{\varphi^x}(a)\Lambda_{\varphi^x}((\sigma_{-i/2}^\Phi(c_k^*))\|^2=\varphi^x(a^*a)\]
So, if $a\in\gN_T$ is in the kernel of $\pi_{\varphi^x}$, we get that $\varphi^x(a^*a)=0$, which, by (iii), implies that $a^*a$ belongs to $\alpha(C_x(X))A$. Let now $e_\lambda$ be an approximate unit of $\gM_T\cap A_n(W)$; we get that $ae_\lambda$ is in $\gM_T\cap A_n(W)$, and in the kernel of $\pi_{\varphi^x}$, and, therefore, by polarisation, belongs to $\alpha(C_x(X))A$. As $ae_\lambda$ is norm converging to $a$, we get that $a\in \alpha(C_x(X))A$, which gives (iv). Then, the first part of (v) is given by (\cite{Bl2}, 3.3), and the proof for $\beta$ is made the same way. The proof of (vi) is then standard. Let's apply (vi) to the opposite measured quantum groupoid $\gG^o=(N^o, M, \beta, \alpha, \varsigma_{N}\Gamma, RTR, T, \nu^o)$ and we get (vii). \end{proof}

\subsection{Definitions}
\label{defmcf}
The left $A_n(W)$-module $A_n(W)\cap\gN_T$ is, using $\alpha_{|C_b(X)}$ (\ref{propmanag}(ii)), a right $C_b(X)$-module, and, equipped with the inner product $(a,b)\mapsto T(b^*a)$, is a inner-product $C_b(X)$-module in the sense of \cite{La}. (We write inner products left linear). 
\newline
We can see its completion $\mathcal E_\Phi$ as the norm closure of the set $\{\Lambda_T(a), a\in \gN_T\cap A_n(W)\}$; then the left-$A_n(W)$-module structure of $\mathcal E_\Phi$ gives that the restriction of $\pi_\Phi$ to $A_n(W)$ can be considered as a $C_b(X)$-linear morphism from $A_n(W)$ into $\mathcal L(\mathcal E_\Phi)$. Taking the specialization at the point $x\in X$, we obtain an Hilbert space $(\mathcal E_\Phi)_x$, which is the completion of the inner product in $A_n(W)\cap\gN_T$ given by $(a, b)\mapsto \varphi^x(b^*a)$; from which we get that $(\mathcal E_\Phi)_x=H_{\varphi^x}$, and that the representation $\pi_x$ obtained by the specialization of $\pi_{\Phi|A_n(W)}$ is equal to $\pi_{\varphi^x}$. We have obtained in \ref{continuousfield}(iii) that $\varphi^x$ is faithful on $A^x$, and in \ref{continuousfield}(iv) that $\pi_{\varphi^x}$ is a continuous field of faithful representations of $A$ (\cite{Bl2}, 2.11). 
\newline
Moreover, we had got in (\cite{E2}, 10.1) that $H_\Phi$ can be written as $\int_X^\oplus H_xd\nu(x)$, where the Hilbert spaces $H_x$ are defined, by separation and completion, from the sesquilinear positive form defined on $D(_\alpha H_\Phi, \nu)$ by $(\xi, \eta)\mapsto \langle\xi, \eta\rangle_{\alpha, \nu}(x)$. It is then starightforward to get that $H_x=H_{\varphi^x}$, and that $\pi_\Phi=\int_X^\oplus \pi_{\varphi^x}d\nu(x)$. Then, it is clear that, if $\xi$ belongs to $D(_\alpha H_\Phi, \nu)$, and $\|\xi\|=1$, that the application $a\mapsto \langle a\xi, \xi\rangle_{\alpha, \nu}$ is a continuous field of states on $A_n(W)$. 
\newline
Using $\beta$, we get another $C_0(X)$-Hilbert module $\mathcal E_{\Phi\circ R}$, and that $\pi_{\Phi\circ R|A_n(W)}$ is a $C_b(X)$-linear morphism from $A_n(W)$ into $\mathcal L(\mathcal E_{\Phi\circ R})$.

\subsection{Corollary}
\label{cortensor}
{\it Let $\gG=(N, M, \alpha, \beta, \Gamma, T, T', \nu)$ be a measured quantum groupoid, such that the von Neuman algebra $\alpha(N)$ is included into the center $Z(M)$. Let $X$ be the spectrum of ${\bf C}^*(\nu)$; using \ref{continuousfield}(i), let us denote $A_n(W)\underset{C_0(X)}{_\beta\otimes^m_\alpha}A_n(W)$ the minimal tensor product of the $C_0(X)-{\bf C}^*$-algebras $A_n(W)$ (via $\beta$) and $A_n(W)$ (via $\alpha$), which is then isomorphic to $A_n(W)\underset{C_0(X)}{_\beta\otimes_\alpha}A_n(W)$ and associative (\ref{defCXC*}). Then :
\newline
(i) the ${\bf C}^*$-algebra $A_n(W)\underset{C_0(X)}{_\beta\otimes^m_\alpha}A_n(W)$ has a faithful representation $\varpi$ on the Hilbert space $H_{\Phi}\underset{\nu}{_\beta\otimes_\alpha}H_{\Phi}$ such that, for all $a_1$ and $a_2$ in $A_n(W)$, we have :
\[\varpi(a_1\otimes a_2)=a_1\underset{N}{_\beta\otimes_\alpha}a_2\]
(iii) for any finite sum with $a_i$ and $b_i$  in $A_n(W)$, we have :
\[\|\sum_{i=1}^n a_i\otimes b_i\|_m=\|\sum_{i=1}^na_i\underset{N}{_\beta\otimes_\alpha}b_i\|\]
(iv) If $x$, $y$ are in $A_n(W)$, the application $x\otimes y\mapsto \|x\underset{N}{_\beta\otimes_\alpha}y\|$ extends to a ${\bf C}^*$-semi-norm on the algebraic tensor product $A_n(W)\odot A_n(W)$ and to a ${\bf C}^*$-norm on the quotient of this algebraic tensor product by the ideal generated by the operators of the form \[\{x\beta(f)\otimes y-x\otimes\alpha(f)y, x,y\in A_n(W), f\in N\}\]
Therefore, the ${\bf C}^*$-algebra $A_n(W)\underset{N}{_\beta\otimes_\alpha}A_n(W)$ can be considered as the min tensor product of the $C_0(X)-{\bf C}^*$-algebra $A_n(W)$ (via $\beta$) with the $C_0(X)-{\bf C}^*$-algebra $A_n(W)$ (via $\alpha$).
\newline
(v) If $\xi\in D(_\alpha H_\Phi, \nu)$ let us denote $\omega_\xi$ the continuous field of states on $A_n(W)$ introduced in \ref{defCXC*} (which is the restriction of the spatial state $\omega_\xi$ on $A_n(W)$); then it is possible to define a positive linear bounded application $id\underset{C_0(X)}{_\beta\otimes^m_\alpha}\omega_\xi$ from $A_n(W)\underset{C_0(X)}{_\beta\otimes^m_\alpha}A_n(W)$ to $A_n(W)$, (and from $M(A_n(W)\underset{C_0(X)}{_\beta\otimes^m_\alpha}A_n(W))$ to $M(A_n(W))$), which is the restriction of the conditional expectation $id\underset{N}{_\beta\otimes_\alpha}\omega_\xi$ to $A_n(W)\underset{C_0(X)}{_\beta\otimes^m_\alpha}A_n(W)$.
\newline
(vi) the application from $C_0(X)$ into $M(A_n(W)\underset{C_0(X)}{_\beta\otimes^m_\alpha}A_n(W))$ defined by :
\[f\mapsto \beta(f)\underset{C_0(X)}{_\beta\otimes^m_\alpha}1=1\underset{C_0(X)}{_\beta\otimes^m_\alpha}\alpha(f)\]
 gives to  $A_n(W)\underset{C_0(X)}{_\beta\otimes^m_\alpha}A_n(W)$ a structure of a continuous field of ${\bf C}^*$-algebras. }
\begin{proof}
Using (\cite{Bl1}, 4.1), we get that there exists a faithful $C_0(X)$-linear representation of $A_n(W)\underset{C_0(X)}{_\beta\otimes^m_\alpha}A_n(W)$ on the Hilbert $A_n(W)$-module $\mathcal E_{\Phi\circ R}\otimes_{C_0(X)}A_n(W)$, which sends the finite sum $\sum_{i=1}^n a_i\otimes b_i$ on the operator $\sum_{i=1}^n\pi_{\Phi\circ R}(a_i)\otimes_{C_0(X)}b_i$ on $\mathcal E_{\Phi\circ R}\otimes_{C_0(X)}A_n(W)$. Let's have a closer look at this last operator, and let's take finite families $x_j\in \gN_{\Phi\circ R}$, $c_j\in \gN_{\Phi}$ ($j=1, ..., m$). With a repeated use of Cauchy-Schwartz inequality, and with the same arguments as in (\cite{L1}, 1.2), one gets that the weight $\Phi$ applied to :
\[\langle\sum_{i=1, j=1}^{i=n, j=m}\pi_{\Phi\circ R}(a_i)\Lambda_{RTR}(x_j)\otimes_{C_0(X)}b_ic_j,  
\sum_{i=1, j=1}^{i=n, j=m}\pi_{\Phi\circ R}(a_i)\Lambda_{RTR}(x_j)\otimes_{C_0(X)}b_ic_j\rangle\]
is less than :
\[\|\sum_{i=1}^n\pi_{\Phi\circ R}(a_i)\otimes_{C_0(X)}b_i\|\Phi(\langle\sum_{j=1}^m\Lambda_{RTR}(x_j)\otimes_{C_0(X)}c_j, \sum_{j=1}^m\Lambda_{RTR}(x_j)\otimes_{C_0(X)}c_j\rangle)\]
But, we easily get that :
\begin{align*}
\Phi(\langle\sum_{j=1}^m\Lambda_{T}(x_j)\otimes_{C_0(X)}c_j, \sum_{j=1}^m\Lambda_{T}(x_j)\otimes_{C_0(X)}c_j\rangle)
&=
\Phi(\sum_j\alpha\circ\beta^{-1}RTR(x_j^*x_j)c_j^*c_j)\\
&=
\|\sum_{j=1}^m \Lambda_{\Phi\circ R}(x_j)\underset{\nu}{_\beta\otimes_\alpha}\Lambda_{\Phi}(c_j)\|^2
\end{align*}
and, therefore, we get that :
\[\|\sum_{i=1, j=1}^{i=n, j=m}\pi_{\Phi\circ R}(a_i)\Lambda_{\Phi\circ R}(x_j)
\underset{\nu}{_\beta\otimes_\alpha}
\pi_{\Phi}(b_i)\Lambda_{\Phi}(c_j)\|^2\]
is less than :
\[\|\sum_{i=1}^n\pi_{\Phi\circ R}(a_i)\otimes_{C_0(X)}b_i\|^2
\|\sum_{j=1}^m \Lambda_{\Phi\circ R}(x_j)\underset{\nu}{_\beta\otimes_\alpha}\Lambda_{\Phi}(c_j)\|^2\]
which, using (\cite{Bl1}, 4.1), is less than :
\[\|\sum_{i=1}^n a_i\otimes b_i\|_m^2
\|\Lambda_{\Phi\circ R}(x_j)\underset{\nu}{_\beta\otimes_\alpha}\Lambda_{\Phi}(c_j)\|^2\]
From which we induce that :
\[\|\sum_{i=1}^n\pi_{\Phi\circ R}(a_i)\underset{N}{_\beta\otimes_\alpha}\pi_{\Phi}(b_i)\|\leq\|\sum_{i=1}^n a_i\otimes b_i\|_m\]
which gives (i). 
\newline
Let us suppose now that $\sum_{i=1}^n\pi_{\Phi\circ R}(a_i)\underset{N}{_\beta\otimes_\alpha}\pi_{\Phi}(b_i)=0$; with the same calculation as above, using the faithfulness of $\Phi$, we get that, for any finite families $(x_j)_{j=1,..m}$ and $(c_j)_{j=1, ..m}$, we have :
\[\sum_{i=1, j=1}^{i=n, j=m}\pi_{\Phi\circ R}(a_i)\Lambda_{RTR}(x_j)\otimes_{C_0(X)}b_ic_j=0\]
which gives that the operator $\sum_{i=1}^n\pi_{\Phi\circ R}(a_i)\otimes_{C_0(X)}b_i$ on $\mathcal E_{\Phi\circ R}\otimes_{C_0(X)}A_n(W)$ is equal to $0$. By the faithfulness of the representation constructed in (\cite{Bl1} 4.1), we get that 
\[\|\sum_{i=1}^n a_i\otimes b_i\|_m=0\]
and, therefore, that $\sum_{i=1}^n a_i\otimes b_i$ belongs to the ideal $ J(A_n(W), A_n(W))$ introduced in (\cite{Bl1}, 2.1). But, as the semi-norm $\sum_{i=1}^n a_i\otimes b_i\mapsto \|\sum_{i=1}^n a_i\otimes b_i\|_m$ is the minimal semi-norm on $(A_n(W)\odot A_n(W))/J(A_n(W), A_n(W))$ (\cite{Bl1}, 2.9), we get (iii). Now (iv) is given by (\cite{Bl1}, 3.1) and \ref{continuousfield}(v), and (v) is trivial.  
\newline
Let's use (iv) and consider $A_n(W)\underset{C_0(X)}{_\beta\otimes^m_\alpha}A_n(W)$ as a ${\bf C}^*$-algebra on $H_\Phi\underset{\nu}{_\beta\otimes_\alpha}H_\Phi$, which is a sub-${\bf C}^*$-algebra of $M\underset{N}{_\beta*_\alpha}M$; on this von Neumann algebra, the slice map $(id\underset{N}{_\beta\otimes_\alpha}T)$ defines a normal faithful operator-valued weight from $M\underset{N}{_\beta*_\alpha}M$ onto $M\underset{N}{_\beta*_\alpha}\alpha(N)=M\underset{N}{_\beta\otimes_\alpha}1$; then, composing with $RTR$, we get a normal faithful operator-valued weight $RTR\underset{N}{_\beta\otimes_\alpha}T$ from $M\underset{N}{_\beta*_\alpha}M$ on $\beta(N)\underset{N}{_\beta\otimes_\alpha}1=1\underset{N}{_\beta\otimes_\alpha}\alpha(N)$.
\newline
Let $A$ be in $\gM_T^+$, $A=\int_X^{\oplus}a^xd\nu(x)$, and $B$ be in$\gM_{RTR}^+$,  $B=\int_X^{\oplus}b^xd\nu(x)$; using \ref{continuousfield}(vi) and (vii), we get that $(RTR\underset{N}{_\beta\otimes_\alpha}T)(B\underset{N}{_\alpha\otimes_\beta}A)$ is equal to the function $x\mapsto \overline{\varphi^x\circ R^x}(b^x)\overline{\varphi^x}(a^x)$, from which we get that this operator-valued weight $RTR\underset{N}{_\beta\otimes_\alpha}T$ is semi-finite. 
\newline
Let now $C\in A_n(W)\underset{C_0(X)}{_\beta\otimes^m_\alpha}A_n(W)\cap\gM_{(RTR\underset{N}{_\beta\otimes_\alpha}T)}^+$; $(RTR\underset{N}{_\beta\otimes_\alpha}T)(C)$ is an element of $N^+$, and therefore a positive bounded continuous function $f$ on $X$; let us suppose that $f(x)=0$; let us write, as in the proof of \ref{continuousfield}(iv), $f_p=[inf(1,f)]^{1/p}$; then $f_p\in N=C_b(X)$, and, as in 
\ref{continuousfield}, we shall obtain that $C$ belongs to the ideal in $A_n(W)\underset{C_0(X)}{_\beta\otimes^m_\alpha}A_n(W)$ generated by $1\underset{C_0(X)}{_\beta\otimes^m_\alpha}\alpha(C_x(X))=\beta(C_x(X))\underset{C_0(X)}{_\beta\otimes^m_\alpha}1$. Taking now $C\in A_n(W)\underset{C_0(X)}{_\beta\otimes^m_\alpha}A_n(W)\cap\gN_{(RTR\underset{N}{_\beta\otimes_\alpha}T)}$, in the kernel of $\pi_{\varphi^x\circ R^x}\underset{C_0(X)}{_\beta\otimes^m_\alpha}\pi_{\varphi^x}$, we shall obtain, using similar arguments, that $C$ belongs also to that ideal; finally, using again an approximate unit, we shall obtain the same result for $C\in A_n(W)\underset{C_0(X)}{_\beta\otimes^m_\alpha}A_n(W)$ in the kernel $\pi_{\varphi^x\circ R^x}\underset{C_0(X)}{_\beta\otimes^m_\alpha}\pi_{\varphi^x}$; therefore, we get that $\pi_{\varphi^x\circ R^x}\underset{C_0(X)}{_\beta\otimes^m_\alpha}\pi_{\varphi^x}$ is a continuous field of faithful representations of $A_n(W)\underset{C_0(X)}{_\beta\otimes^m_\alpha}A_n(W)$, which proves (vi). 
 \end{proof}

\subsection{Theorem}
\label{thcentral2}
{\it Let $(N, M, \alpha, \beta, \Gamma, T, T', \nu)$ be a measured quantum groupoid, and let us suppose that the von Neuman algebra $\alpha(N)$ is included into the center $Z(M)$;
then, for all $x$ in the ${\bf C}^*$-algebra $A_n(W)$, $\Gamma (x)$ belongs to the multipliers of the  ${\bf C}^*$-algebra $A_n(W)\underset{N}{_\beta\otimes_\alpha}A_n(W)$; using \ref{cortensor}, we get that the restriction of $\Gamma$ to $A_n(W)$ sends $A_n(W)$ into  $M(A_n(W)\underset{{\bf C}^*(\nu)}{_\beta\otimes_\alpha^m}A_n(W))$. }

\begin{proof}
Let $\xi$, $\eta$ in $D(_\alpha H, \nu)$; using \ref{thcentral1}, the operator $\Gamma ((id*\omega_{\xi, \eta})(W))$ is a strict limit of elements in $A_n(W)\underset{N}{_\beta\otimes_\alpha}A_n(W)$, and therefore belongs to $M(A_n(W)\underset{N}{_\beta\otimes_\alpha}A_n(W))$, from which we get the result, by definition of $A_n(W)$. \end{proof}

\subsection{Theorem}
\label{thresume}
{\it Let $\gG=(N, M, \alpha, \beta, \Gamma, T, T', \nu)$ be a measured quantum groupoid, and let us suppose that the von Neuman algebra $\alpha(N)$ is included into the center $Z(M)$; let $X$ be the spectrum of ${\bf C}^*(\nu)$, and, for $x\in X$, let $C_x(X)$ be the subalgebra of $C_0(X)$ made of functions which vanish at $x$; let $R$ be the co-inverse of $\gG$; then :
\newline
(i) thanks to the $*$-homomorphism $\alpha_{|C_0(X)}$ (resp. $\beta_{|C_0(X)}$), the ${\bf C}^*$-algebra $A_n(W)$ is a continuous field over $X$ of ${\bf C}^*$-algebras; therefore, (\cite{Bl1}, 4.1), Blanchard's minimal tensor product $A_n(W)\underset{C_0(X)}{_\beta\otimes_\alpha^m}A_n(W)$ is associative; 
\newline
(ii) the restriction of the coproduct to $A_n(W)$ sends $A_n(W)$ into  $M(A_n(W)\underset{C_0(X)}{_\beta\otimes_\alpha^m}A_n(W))$; 
\newline
(iii) for any  $a\in A_n(W)^+$, we can identify $T(a)$ with the (image by $\alpha$ of the) function $x\mapsto \varphi^x(a)$, which is lower semi-continuous (and bounded continuous is $a\in\gM_T^+$).
\newline
(iv) for any $f\in C_b(X)^+$ and $a\in A_n(W)^+$, we have $\varphi^x(\alpha (f)a)=f(x)\varphi^x(a)$, and $\varphi^x(a)=0$ if and only if $a\in\alpha(C_x(X))A_n(W)$.
\newline
(v) the representations $\pi_{\varphi^x}$ form a continuous field of faithful representations of $A_n(W)$, when considered, thanks to $\alpha$, as a continuous field over $X$ of ${\bf C}^*$-algebras;  
\newline
(vi) there exists an linear anti-$*$-isomorphism $R^x$ from $A_n(W)/\alpha(C_x(X))A_n(W)$ onto $A_n(W)/\beta(C_x(X))A_n(W)$, and, considering, thanks to $\beta$, $A_n(W)$ as a continuous field over $X$ of ${\bf C}^*$-algebras,  $\varphi^x\circ R^x$ is then a field of lower continuous faithful weights, such that $\pi_{\varphi^x\circ R^x}$ form a continuous field of faithful representations of $A_n(W)$.
\newline
(vii) the representations $\pi_{\varphi^x\circ R^x}\underset{C_0(X)}{_\beta\otimes^m_\alpha}\pi_{\varphi^x}$ form a continuous field of faithful representations of $A_n(W)\underset{C_0(X)}{_\beta\otimes^m_\alpha}A_n(W)$, which gives to that ${\bf C}^*$-algebra a structure of continuous field over $X$ of ${\bf C}^*$-algebra, thanks to the application which sends $f\in C_0(X)$ on $1\underset{C_0(X)}{_\beta\otimes^m_\alpha}\alpha(f)=\beta(f)\underset{C_0(X)}{_\beta\otimes^m_\alpha}1$. 
\newline
(viii) for any $a\in A_n(W)^+\cap\gM_T$, and $\eta\in D((H_\Phi)_\beta, \nu)$, such that $\|\eta\|=1$, we have, for all $x\in X$ :}
\[\varphi^x[(\omega_\eta\underset{C_0(X)}{_\beta\otimes^m_\alpha}id)\Gamma(a)]=\varphi^x (a)\]
\begin{proof}
(i), (iii), (iv), (v) are taken from \ref{continuousfield}(v), (ii), (iii) and (iv); (ii) is \ref{thcentral2}; (vi) is an easy corollary, (vii) was obtained in \ref{cortensor}(vi), and (vii) is just given by restriction of the formula on $\gM_T^+$. \end{proof}

\section{Abelian measured quantum groupoids}
\label{abelian}
We consider now the case of an "abelian" measured quantum groupoid (i.e. a measured quantum groupoid $\gG=(N, M, \alpha, \beta, \Gamma, T, T', \nu)$ where the underlying von Neuman algebra itself is abelian); then we prove that it is possible to put on the spectrum of the ${\bf C}^*$-algebra $A_n(W)$ a structure of a locally compact groupoid, whose basis is the spectrum of ${\bf C}^*(\nu)$ (\ref{thgroupoid}). Starting from a measured groupoid equipped with a left-invariant Haar system, we recover Ramsay's theorem which says that this groupoid is measure-equivalent to a locally compact one (\ref{ramsay}). 
\subsection{Proposition}
\label{thgroupoid}
{\it Let $\gG=(N, M, \alpha, \beta, \Gamma, T, T', \nu)$ be a measured quantum groupoid, and let us suppose that the von Neuman algebra $M$ is abelian; let us write $\mathcal G$ for the spectrum of the ${\bf C}^*$-algebra $A_n(W)$, and $\mathcal G^{(0)}$ for the spectrum of the ${\bf C}^*$-algebra ${\bf C}^*(\nu)$. Then :
\newline
(i) there exists a continuous open application $r$ from $\mathcal G$ onto $\mathcal G^{(0)}$, such that, for all $f\in C_0(\mathcal G^{(0)})$, we have $\alpha (f)=f\circ r$; there exists a continous open application $s$ from $\mathcal G$ onto $\mathcal G^{(0)}$, such that, for all $f\in C_0(\mathcal G^{(0)})$, we have $\beta (f)=f\circ s$. 
\newline
(ii) there exists a partially defined multiplication on $\mathcal G$, which gives to $\mathcal G$ a structure of locally compact groupoid, with $\mathcal G ^{(0)}$ as set of units. 
\newline
(iii) the application defined for all $F$ continuous, positive, with compact support in $\mathcal G$, by $F\mapsto \alpha^{-1}(T(F))(u)$, defines a positive Radon measure $\lambda^u$ on $\mathcal G$, whose support is $\mathcal G^u$. The measures $(\lambda^u)_{u\in \mathcal G ^{(0)}}$ are a Haar system on $\mathcal G$.
\newline
(iv) the trace $\nu$ on ${\bf C}^*(\nu)$ leads to a quasi-invariant measure (denoted again by $\nu$) on $\mathcal G ^{(0)}$. Let $\mu=\int_{\mathcal G^{(0)}}\lambda^u d\nu(u)$; then : 
\[(N, M, \alpha, \beta, \Gamma)=(L^{\infty}(\mathcal G ^{(0)}, \nu), L^{\infty}(\mathcal G , \mu), r_{\mathcal G}, s_{\mathcal G}, \Gamma_{\mathcal G})\]
where $r_{\mathcal G}$, $s_{\mathcal G}$, $\Gamma_{\mathcal G}$ have been defined in \ref{Hbimod}. Moreover, then, the operator-valued weights $T$ and $RTR$ are given, for any positive $F$ in $L^{\infty}(\mathcal G , \nu)$ by :
\[T(F)(u)=\int_{\mathcal G}Fd\lambda^u\]
\[RTR(F)(u)=\int_{\mathcal G}Fd\lambda_u\]
where $\lambda_u$ is the image of $\lambda^u$ under the application $(x\mapsto x^{-1})$. Therefore, with the notations of \ref{gd2}, we have : $\gG=\gG(\mathcal G)$. }

\begin{proof}
As $\alpha(N)\subset M(A_n(W))$, we can construct by restriction a continuous application $r$ from $\mathcal G$ into $\mathcal G^{(0)}$ such that, for all $f\in C_0(\mathcal G^{(0)})={\bf C}^*(\nu)$, we have $\alpha (f)=f\circ r$; we can construct the same way a continuous application $s$ from $\mathcal G$ into $\mathcal G^{(0)}$ such that, for all $f\in C_0(\mathcal G^{(0)})={\bf C}^*(\nu)$, we have $\beta (f)=f\circ r$. The applications $r$ and $s$ are open by (\cite{Bl2}, 3.14), which gives (i).
\newline
The application $R$ from $A_n(W)$ into itself leads to an involutive application in $\mathcal G$, we shall write $x\mapsto x^{-1}$, and, using that $R\circ\alpha=\beta$, we get that $r(x^{-1})=s(x)$ and $s(x^{-1})=r(x)$.
\newline
Thanks to \ref{continuousfield}, we may apply (\cite{Bl1}, 3.1) to $A_n(W)$, which we identify to $C_0(\mathcal G)$, and we obtain that the commutative ${\bf C}^*$-algebra $A_n(W)\underset{{\bf C}^*(\nu)}{_\beta\otimes^m_\alpha}A_n(W)$ is the quotient of the ${\bf C}^*$-algebra $A_n(W)\otimes A_n(W)$ (identified with $C_0(\mathcal G^2)$) by the ideal generated by all the functions $(x_1, x_2)\mapsto f(s(x_1))g(x_1, x_2)-f(r(x_2))g(x_1, x_2)$, where $x_1$, $x_2$ are in $\mathcal G$, $f$ in ${\bf C}^*(\nu)$ (identified with $C_0(\mathcal G^{(0)})$), and $g$ in $C_0(\mathcal G^2)$. So, a non zero character on $A_n(W)\underset{{\bf C}^*(\nu)}{_\beta\otimes^m_\alpha}A_n(W)$ is a couple $(x_1, x_2)$ in $\mathcal G^2$ such that $s(x_1)=r(x_2)$; let us write $\mathcal G^{(2)}$ for the subset of such elements of $\mathcal G^2$. 
\newline
Therefore, we can identify $A_n(W)\underset{{\bf C}^*(\nu)}{_\beta\otimes^m_\alpha}A_n(W)$ to $C_0(\mathcal G^{(2)})$. Therefore, we see that the restriction of $\Gamma$ to $A_n(W)$ leads to a continuous application from $\mathcal G^{(2)}$ into $\mathcal G$, which gives to $\mathcal G$ a structure of locally compact groupoid, which is (ii). 
\newline
As $A_n(W)\cap\gM_T$ is a dense ideal in $A_n(W)$, it contains the ideal $\mathcal K(\mathcal G)$ of continuous functions on $\mathcal G$, with compact support; for all $F$ in $\mathcal K(\mathcal G)$, $\alpha^{-1}(T(F))$ belongs to $C_b(\mathcal G^{(0)})$, and, for all $u\in\mathcal G^{(0)}$, $F\mapsto\alpha^{-1}(T(F))(u)$ defines a non zero positive Radon measure $\lambda^u$ on $\mathcal G$; it is now straightforward to get, from the left invariance of $T$, that $(\lambda^u)_{u\in \mathcal G^{(0)}}$ is a Haar system on the groupoid. Starting from $R\circ T\circ R$, we obtain measures $\lambda_u$, which are the images of $\lambda^u$ by the inverse. 
\newline
The modulus $\delta$ of the measured quantum groupoid gives that the trace $\nu$ on ${\bf C}^*(\nu)$ is a quasi-invariant measure on $\mathcal G^{(0)}$. 
\newline
Now, by density reasons, we shall identify $N$ with $L^{\infty}(\mathcal G^{(0)}, \mu)$, $M$ with $L^{\infty}(\mathcal G, \mu)$, where $\mu$ is the measure on $\mathcal G$ constructed from $\mu$ and the Haar system, $\alpha$ with $r_{\mathcal G}$, $\beta$ with $s_{\mathcal G}$, $\Gamma$ with $\Gamma_{\mathcal G}$, and we obtain the required formulae for the left and right Haar systems.
\end{proof}

\subsection{Ramsay's theorem [Ra]}
\label{ramsay}
{\it Let $\mathcal G$ be a measured groupoid, with $\mathcal G^{(0)}$ as space
of units, and $r$ and $s$ the range and source functions from $\mathcal G$ to $\mathcal G^{(0)}$, with a Haar system $(\lambda^u)_{u\in \mathcal G^{(0)}}$ and a quasi-invariant measure $\nu$ on $\mathcal G^{(0)}$. Let us write $\mu=\int_{\mathcal G^{(0)}}\lambda^ud\nu$. Let $\Gamma_{\mathcal G}$, $r_{\mathcal G}$, $s_{\mathcal G}$ be the morphisms associated in \ref{gd}. Then, there exists a locally compact groupoid $\tilde{\mathcal G}$, with set of units $\tilde{\mathcal G}^{(0)}$, with a Haar system $(\tilde{\lambda}^u)_{u\in \tilde{\mathcal G}^{(0)}}$, and a quasi-invariant measure $\tilde{\nu}$ on $\tilde{\mathcal G}^{(0)}$, such that, if $\tilde{\mu}=\int_{\tilde{\mathcal G}^{(0)}}\tilde{\lambda}^ud\tilde{\nu}$, we get that the abelian measured quantum groupoids $\gG(\mathcal G)$ and $\gG(\tilde{\mathcal G})$ are isomorphic. }
\begin{proof}
Let us apply \ref{thgroupoid} to the commutative measured quantum groupoid $\gG(\mathcal G)$ constructed from the measured groupoid $\mathcal G$. Then, we get the result. \end{proof}

\section{Measured fields of locally compact quantum groups}
\label{fieldqg}
In this chapter, we define a notion of measured field of locally compact quantum groups (\ref{deffieldqg}), which was underlying in \cite{Bl2}. We construct then from such a field a measured quantum groupoid (\ref{mqgfield}), and we show that the measured quantum groupoids obtained this way are exactly the measured quantum groupoids with a central basis, studied in chapters \ref{central} and \ref{mcfield}, such that the dual object is of the same kind (\ref{bicentral}). We finish by recalling concrete examples (\ref{SU}, \ref{ax+b}, \ref{Emu}) given by Blanchard, which give examples of measured quantum groupoids.

\subsection{Definition (\cite{L2}, 17.3)}
\label{deffieldqg}
Let $(X, \nu)$ be a $\sigma$-finite standard measure space; let us take $\{M^x, x\in X\}$ a measurable field of von Neumann algebras over $(X, \nu)$ and $\{\varphi^x, x\in X\}$ (resp. $\{\psi^x\}$) a measurable field of normal semi-finite faithful weights on $\{M^x\}$ (\cite{T}, 4.4). Moreover, let us suppose that :
\newline
(i) there exits a measurable field of injective $*$-homomorphisms $\Gamma^x$ from $M^x$ into $M^x\otimes M^x$ (which is also a measurable field of von Neumann algebras, on the measurable field of Hilbert spaces $H_{\varphi^x}\otimes H_{\varphi^x}$). 
\newline
(ii) for almost all $x\in X$, $\textbf{G}^x=(M^x, \Gamma^x, \varphi^x, \psi^x)$ is a locally compact quantum group (in the von Neumann sense (\cite{KV2}). 
\newline
In that situation, we shall say that $(M^x, \Gamma^x, \varphi^x, \psi^x, x\in X)$ is a measurable field of locally compact quantum groups over $(X, \nu)$. 

\subsection{Theorem (\cite{L2}, 17.3)}
\label{mqgfield}
{\it Let $\textbf{G}^x=(M^x, \Gamma^x, \varphi^x, \psi^x, x\in X)$ be a measurable field of locally compact quantum groups over $(X, \nu)$. Let us define :
\newline
(i) $M$ as the von Neumann algebra made of decomposable operators $\int_X^{\oplus}M^xd\nu(x)$, and $\alpha$ the $*$-isomorphism which sends $L^\infty(X, \nu)$ into the algebra of diagonalizable operators, which is included in $Z(M)$.
\newline
(ii) $\Phi$ (resp. $\Psi$) as the direct integral $\int_X^{\oplus}\varphi^x d\nu(x)$ (resp. $\int_X^{\oplus}\psi^xd\nu(x)$). Then, the Hilbert space $H_\Phi$ is equal to the direct integral $\int_X^{\oplus}H_{\varphi^x} d\nu(x)$, the relative tensor product $H_\Phi\underset{\nu}{_\alpha\otimes_\alpha}H_\Phi$ is equal to the direct integral $\int_X^{\oplus}(H_{\varphi^x}\otimes H_{\varphi^x})d\nu(x)$, and the product $M\underset{N}{_\alpha*_\alpha}M$ is equal to the direct integral $\int_X^{\oplus}(M^x\otimes M^x)d\nu(x)$.
\newline
(iii) $\Gamma$ as the decomposable $*$-homomorphism $\int_X^{\oplus}\Gamma^xd\nu(x)$, which sends $M$ into $M\underset{N}{_\alpha*_\alpha}M$.
\newline
(iv) $T$ (resp. $T'$) as an operator-valued weight from $M$ into $\alpha (L^\infty(X, \nu))$ defined the following way : $a\in M^+$ represented by the field $\{a^x\}$ belongs to $\gM_T^+$ if, for almost all $x\in X$, $a^x$ belongs to $\gM_{\varphi^x}$ (resp. $\gM_{\psi^x}$), and the function $x\mapsto \varphi^x(a^x)$ (resp. $x\mapsto \psi^x(a^x)$) is essentiallly bounded; then $T(a)$ (resp. $T'(a)$) is defined as the image under $\alpha$ of this function. 
\newline
Then, $(L^\infty(X, \nu), M, \alpha, \alpha, \Gamma, T, T', \nu)$ is a measured quantum groupoid, we shall denote by $\int_X^{\oplus}\textbf{G}^xd\nu(x)$.}

\begin{proof}
The fact that $H_\Phi=\int_X^{\oplus}H_{\varphi^x}d\nu(x)$ is given by (\cite{T}, 6.3.11). Then, we can identify, for an element $a\in\gN_\Phi$ represented by the field $\{a^x\}$, $\Lambda_\Phi (a)$ with $\int_X^{\oplus}\Lambda_{\varphi^x}(a^x)d\nu(x)$. If $\xi\in H_\Phi$, $\xi$ can be represented by a square integrable field of vectors $\{\xi^x\}$; moreover, if $\xi\in D(_\alpha H_\Phi, \nu)$, we get that there exists $C>0$ such that, for all $f\in L^\infty(X,\nu)\cap L^2(X,\nu)$:
\[\int_X\|f(x)\xi^x\|^2d\nu(x)\leq C\int_X |f(x)|^2d\nu(x)\]
which gives that the function $x\mapsto \|\xi^x\|^2$ is essentially bounded. It is then straightforward to get that this function is equal to the element $\langle\xi, \xi\rangle_{\alpha, \nu}$ of $L^\infty(X, \nu)$. 
\newline
So, the relative tensor product $H_\Phi\underset{\nu}{_\alpha\otimes_\alpha}H_\Phi$ is the completion of the algebraic tensor product $D(_\alpha H_\Phi, \nu)\odot \Lambda_\Phi(\gN_\Phi)$ by the scalar product defined by the formula (where $b\in\gN_\Phi$ is represented by the field $\{b^x\}$ and $\eta\in D(_\alpha H_\Phi, \nu)$ is represented by the vector field $\{\eta^x\}$):
\begin{align*}
(\xi\odot\Lambda_\Phi(a)|\eta\odot\Lambda_\Phi(b))
&=
(\alpha(\langle\xi, \eta\rangle_{\alpha, \nu})\Lambda_\Phi(a)|\Lambda_\Phi(b))\\
&=\int_X^{\oplus}(\xi^x\otimes\Lambda_{\varphi^x}(a^x)|\eta^x\otimes\Lambda_{\varphi^x}(b^x))d\nu(x)
\end{align*}
from which we get that $H_\Phi\underset{\nu}{_\alpha\otimes_\alpha}H_\Phi=\int_X^{\oplus}(H_{\varphi^x}\otimes H_{\varphi^x})d\nu(x)$; it is now straightforward. to obtain that $M\underset{N}{_\alpha*_\alpha}M=\int_X^{\oplus}(M^x\otimes M^x)d\nu(x)$.  We then get that the $*$-homomorphism $\Gamma$ defined in (iii) is a coassociative coproduct which makes $(L^\infty(X, \nu), M, \alpha, \alpha, \Gamma)$ a Hopf-bimodule. 
\newline
Then, $T$ as defined in (iv) is an operator-valued weight from $M$ to $\alpha(L^\infty(X, \nu))$, which verify, by definition $\nu\circ\alpha^{-1}\circ T=\Phi$; therefore, $T$ is normal, faithful, semi-finite. If $a\in\gM_T^+$ and is represented by the field $\{a^x\}$, for almost all $x\in X$, $a^x$ belongs to $\gM_{\varphi^x}$, therefore, $\Gamma^x(a^x)$ belongs to $\gM_{id\otimes\varphi^x}$, and $(id\otimes\varphi^x)\Gamma^x(a^x)=\varphi^x(a^x)1_{M^x}$. Let now $\xi\in D(_\alpha H_\Phi, \nu)$, represented by the vector field $\{\xi^x\}$. We have :
\begin{align*}
\Phi[(\omega_\xi\underset{N}{_\alpha\otimes_\alpha}id)\Gamma (a)]
&=
\int_X^{\oplus}(\omega_{\xi^x}\otimes\varphi^x)\Gamma^x(a^x)d\nu(x)\\
&=\int_X^{\oplus}\varphi^x(a^x)\omega_{\xi^x}(1_{M^x})d\nu(x)\\
&=\Phi(A)\omega_\xi (1)
\end{align*}
from which we get that $(id\underset{N}{_\alpha\otimes_\alpha}\Phi)\Gamma (a)=\Phi(a)1$, and, therefore, that :
\[(id\underset{N}{_\alpha\otimes_\alpha}T)\Gamma (a)=T(a)\underset{N}{_\alpha\otimes_\alpha}1\]
the right-invariance for $T'$ is proved the same way. 
\newline
Finally, thanks to (\cite{T}, 4.8), we have, for any $a\in M$, represented by the field $\{a^x\}$ and $t\in\mathbb{R}$, $\sigma_t^\Phi(a)=\int_X^{\oplus}\sigma_t^{\varphi^x}(a^x)d\nu(x)$ and $\sigma_t^\Psi(a)=\int_X^{\oplus}\sigma_t^{\psi^x}(a^x)d\nu(x)$. Therefore, the commutation, for almost all $x\in X$, of $\sigma_t^{\varphi^x}$ and $\sigma_t^{\psi^x}$ (\cite{KV1} 6.8) gives that $\nu$ is relatively invariant. 
\end{proof}

\subsection{Proposition}
\label{propmqgfield}
{\it Let $(X,\nu)$ be a $\sigma$-finite standard measure space, and $\{\textbf{G}^x, x\in X\}$ a measurable field of locally compact quantum groups, as defined in \ref{deffieldqg}; let $\int_X^{\oplus}\textbf{G}^xd\nu(x)$ be the measured quantum groupoid constructed in \ref{mqgfield}; then :
\newline
(i) we have $\alpha=\beta=\hat{\beta}$;
\newline
(ii) the pseudo-multiplicative unitary of the the measured quantum groupoid is a unitary on $H_\Phi\underset{\nu}{_\alpha\otimes_\alpha}H_\Phi$, which is equal to the decomposable operator $\int_X^{\oplus}W^xd\nu(x)$, where $W^x$ is the multiplicative unitary associated to the locally compact quantum group $\textbf{G}^x$. 
\newline
(iii) we have :}
\[\widehat{\int_X^{\oplus}\textbf{G}^xd\nu(x)}=\int_X^{\oplus}\widehat{\textbf{G}^x}d\nu(x)\]

\begin{proof}
The fact that $\beta=\alpha$ is given in the definition of $\int_X^{\oplus}\textbf{G}^xd\nu(x)$; moreover, as $\alpha (L^\infty(X, \nu))\subset Z(M)$, we have $\hat{\beta}=\alpha$, which is (i). Therefore, the pseudo-multiplicative unitary $W$ is a unitary on $H_\Phi\underset{\nu}{_\alpha\otimes_\alpha}H_\Phi=\int_X^{\oplus}(H_{\varphi^x}\otimes H_{\varphi^x})d\nu(x)$. Moreover, using \ref{thL1}(i), we get, for all $v\in D(_\alpha H_\Phi, \nu)$, represented by the vector field $\{v^x\}$, and $a\in\gN_\Phi$, represented by the field $\{a^x\}$ :
\[W^*(v\underset{\nu}{_\alpha\otimes_\alpha}\Lambda_\Phi(a))=\sum_i\xi_i\underset{\nu}{_\alpha\otimes_\alpha}\Lambda_\Phi[(\omega_{v, \xi_i}\underset{N}{_\alpha\otimes_\alpha}id)\Gamma(a)]\]
where $(\xi_i)_{i\in I}$ is an orthogonal $(\alpha, \nu)$-basis of $H_\Phi$; each $\xi_i$ is in $D(_\alpha H_\Phi, \nu)$, and is represented by a field $\{\xi_i^x\}$; for almost all $x\in X$, the vectors $(\xi_i^x)_{i\in I}$ (more precisely, those which are not equal to $0$) are an orthogonal basis of $H_{\varphi^x}$ (\cite{E2},10.1); therefore, we get :
\begin{align*}
W^*(v\underset{\nu}{_\alpha\otimes_\alpha}\Lambda_\Phi(a))
&=
\int_X^{\oplus}\sum_i(\xi_i^x\otimes\Lambda_{\varphi^x}[(\omega_{v^x,\xi_i^x}\otimes id)\Gamma^x(a^x)]d\nu(x)\\
&=\int_X^{\oplus}(W^x)^*(v^x\otimes\Lambda_{\varphi^x}(a^x))d\nu(x)
\end{align*}
which gives (ii). Then, we get that the Hopf-bimodules underlying to $\widehat{\int_X^{\oplus}\textbf{G}^xd\nu(x)}$ and $\int_X^{\oplus}\widehat{\textbf{G}^x}d\nu(x)$ are the same; the only result to prove is the equality of the dual operator-valued weights, which is left to the reader. \end{proof}

\subsection{Proposition}
\label{propbicentral}
{\it Let $\gG=(N, M, \alpha, \beta, \Gamma, T, T', \nu)$ be a measured quantum groupoid and $\widehat{\gG}=(N, \widehat{M}, \alpha, \hat{\beta}, \widehat{\Gamma}, \hat{T}, \widehat{R}\hat{T}\widehat{R}, \nu)$ its dual measured quantum groupoid. Then, are equivalent:
\newline
(i) $\alpha (N)\subset Z(M)\cap Z(\widehat{M})$.
\newline
(ii) $\alpha=\beta=\hat{\beta}$. }
\begin{proof} 
This is clear by using \ref{propcentral} twice (for $\gG$ and $\widehat{\gG}$). \end{proof}

\subsection{Theorem}
\label{bicentral}
{\it Let $\gG=(N, M, \alpha, \beta, \Gamma, T, T', \nu)$ be a measured quantum groupoid and $\widehat{\gG}=(N, \hat{M}, \alpha, \hat{\beta}, \widehat{\Gamma}, \hat{T}, \widehat{R}\hat{T}\widehat{R}, \nu)$ its dual measured quantum groupoid; let $W$ and $\widehat{W}$ be the pseudo-multiplicative unitaries associated, and $\Phi=\nu\circ\alpha^{-1}\circ T$ (resp. $\widehat{\Phi}=\nu\circ\alpha^{-1}\circ\hat{T}$); let us suppose that $\alpha (N)$ is central in both $M$ and $\widehat{M}$; let $X$ be the spectrum of ${\bf C}^*$algebra ${\bf C}^*(\nu)$, we shall therefore identify with $C_0(X)$; for any $x\in X$, let $C_x(X)$ be the subalgebra of $C_0(X)$ made of functions which vanish at $x$; let $A_n(W)$ be the sub-${\bf C}^*$-algebra of $M$ introduced in \ref{AW} and \ref{propmanag}(ii), which is, thanks to $\alpha_{|C_0(X)}$, a continuous field over $X$ of ${\bf C}^*$-algebras (\ref{thresume}); let $\varphi^x$ be the desintegration of $\Phi_{|A_n(W)}$ over $X$ given in \ref{continuousfield}(ii); then, by \ref{continuousfield}(iii), $\varphi^x$ is a lower semi-continuous weight  on $A_n(W)$, faithful when considered on $A_n(W)/\alpha(C_x(X))A_n(W)$, and the representation $\pi_{\varphi^x}$ form a continuous field of faithful representation of $A_n(W)$. Then :
\newline
(i) the Hilbert space $H_\Phi\underset{\nu}{_\alpha\otimes_\alpha}H_\Phi$ is equal to $\int_X^{\oplus}H_{\varphi^x}\otimes H_{\varphi^x}d\nu(x)$. 
\newline
(ii) the von Neumann algebra $M\underset{N}{_\alpha*_\alpha}M$ is equal to :
\[\int_X^{\oplus}\pi_{\varphi^x}(A_n(W)/\alpha(C_x(X))A_n(W))"\otimes\pi_{\varphi^x}(A_n(W)/\alpha(C_x(X))A_n(W))"d\nu(x)\]
\newline
(ii) the coproduct $\Gamma_{|A_n(W)}$ can be desintegrated in $\Gamma_{|A_n(W)}=\int_X^{\oplus}\Gamma^xd\nu(x)$, where $\Gamma^x$ is a continuous field of coassociative coproducts on $A_n(W)/\alpha(C_x(X))A_n(W)$. 
\newline
(iii) $R^x$ is a anti-$*$-automorphism of $A_n(W)/\alpha(C_x(X))A_n(W)$, and, for all $x\in X$, $(A_n(W)/\alpha(C_x(X))A_n(W), \Gamma^x, \varphi^x, \varphi^x\circ R^x)$ is a locally compact quantum group (in the ${\bf C}^*$-sense), we shall denote ${\bf G}^x$. We shall denote also ${\bf G}^x$ its von Neumann version. 
\newline
(iv) we have, with the notations of \ref{mqgfield}, $\gG=\int_X^{\oplus}{\bf G}^x d\nu(x)$. }

\begin{proof}
Let $a\in A_n(W)\cap\gN_\Phi$ and $b\in A_n(W)\cap\gN_\Phi$; using \ref{defmcf}, we can write $\Lambda_\Phi(a)=\int_X^{\oplus}\Lambda_{\varphi^x}(a^x)d\nu(x)$ and $\Lambda_\Phi(b)=\int_X^{\oplus}\Lambda_{\varphi^x}(b^x)d\nu(x)$, where $a^x$ (resp. $b^x$) is the image of $a$ (resp. $b$) in $A_n(W)/\alpha(C_x(X))A_n(W)$; as the linear set made of vectors of the form $\Lambda_\Phi(a)$, for all  is dense in $H_\Phi$, the relative tensor product $H_\Phi\underset{\nu}{_\alpha\otimes_\alpha}H_\Phi$ is the completion of the algebraic tensor product $D(_\alpha H_\Phi, \nu)\odot\Lambda_\Phi(\gN_\Phi\cap A_n(W))$ by the scalar product defined, if $\xi=\int_X^{\oplus}\xi^x d\nu(x)$, $\eta=\int_X^{\oplus}\eta^x d\nu(x)$ are in $D(_\alpha H_\Phi, \nu)$ by the formula :
\begin{align*}
(\xi\odot\Lambda_\Phi(a)|\eta\odot\Lambda_\Phi(b))
&=
(\alpha(\langle\xi, \eta\rangle_{\alpha, \nu})\Lambda_\Phi(a)|\Lambda_\Phi(b))\\
&=\int_X^{\oplus}(\xi^x\Lambda_{\varphi^x}(a^x)|\eta^x\Lambda_{\varphi^x}(b^x))d\nu(x)\\
&=\int_X^{\oplus}(\xi^x\otimes\Lambda_{\varphi^x}(a^x)|\eta^x\otimes\Lambda_{\varphi^x}(b^x))d\nu(x)
\end{align*}
from which we get (i). Then (ii) is a direct corollary.
\newline
As $\Gamma(\alpha(f))=\alpha(f)\underset{N}{_\alpha\otimes_\alpha}1=1\underset{N}{_\alpha\otimes_\alpha}\alpha(f)$, we get, using \ref{thresume}, that $\Gamma_{|A_n(W)}$ can be desintegrated into a continuous field of $*$-homomorphisms $\Gamma^x$ from $A_n(W)/\alpha(C_x(X))A_n(W)$ into $[A_n(W)/\alpha(C_x(X))A_n(W)]\otimes^m[A_n(W)/\alpha(C_x(X))A_n(W)]$. Moreover, if $a\in A_n(W)\cap\gM_T^+$, we have, for all $x\in X$, $(id\otimes\varphi^x)\Gamma^x(a)=\varphi^x(a)1$. 
\newline
So, if $a$ in $A_n(W)\cap\gM_T^+$ verify $\Gamma^x(a)=0$, it implies that $\varphi^x(a)=0$ and, therefore, by \ref{continuousfield}(ii), that $a\in \alpha(C_x(X))A_n(W)$. Let now $e_\lambda$ an approximate unit in $A_n(W)\cap\gM_T$, and let $a\in A_n(W)$ such that $\Gamma^x(a)=0$; as we have $\Gamma^x (ae_\lambda)=0$, we get that $ae_\lambda$ belongs to $\alpha(C_x(X))A_n(W)$, for almost all $x$, and we get the same result for $a$. Therefore, we get that $\Gamma^x$ is injective on $A_n(W)/\alpha(C_x(X))A_n(W)$. The coassociativity of $\Gamma^x$ is just a corollary of the coassociativity of $\Gamma$. 
\newline
It is clear that $R^x$ is a $*$-anti-automorphism of $A_n(W)/\alpha(C_x(X))A_n(W)$, which will be a co-inverse for $\Gamma^x$; we shall therefore get that $\varphi^x\circ R^x$ is right-invariant with respect to $\Gamma$; in order to prove that $(A_n(W)/\alpha(C_x(X))A_n(W), \Gamma^x, \varphi^x, \varphi^x\circ R^x)$ is a (${\bf C}^*$-version of a) locally compact quantum group, we shall extends all these data to $\pi_{\varphi^x}(A_n(W)\alpha(C_x(X))A_n(W))"$, and prove that the objects obtained are a locally compact quantum group in the von Neumann sense. 
\newline
In fact, from (i) and (ii), we get that $\Gamma$ can be desintegrated in $\Gamma=\int_X^{\oplus}\tilde{\Gamma}^xd\nu(x)$, where $\tilde{\Gamma}^x$ is a $*$-homomorphism from $\pi_{\varphi^x}(A_n(W)\alpha(C_x(X))A_n(W))"$ into its von Neumann tensor product by itself; moreover, by unicity of the desintegration procedure, we get that, for almost all $x\in X$, $\Gamma^x$ is equal to the restriction of 
$\tilde{\Gamma}^x$ to $A_n(W)/\alpha(C_x(X))A_n(W)$, which proves that $\Gamma^x$ extends at the von Neumann level. We shall therefore write $\Gamma^x$ instead of $\tilde{\Gamma}^x$. We obtain that, for almost all $x$, $\Gamma^x(1)=1$, and $(\Gamma^x\otimes id)\Gamma^x=(id\otimes\Gamma^x)\Gamma^x$ from the properties of $\Gamma$. Moreover, we had got that the restriction of $\Gamma^x$ to $A_n(W)/\alpha(C_x(X))A_n(W)$ is injective; so, if $a\in \pi_{\varphi^x}(A_n(W)/\alpha(C_x(X))A_n(W))"^+$ verify $\Gamma^x(a)=0$, we get that $a$ is an increasing limit of elements $a_n$ in $A_n(W)/\alpha(C_x(X))A_n(W)$ such that $\Gamma^x(a_n)=0$; so, we get that $a_n=0$, and $a=0$, which finishes the proof of (ii).  Then, (iii) is given by \ref{thresume}(vi) and similar calculations, and (iv) is straightforward. 
\end{proof}

\subsection{Theorem}
\label{cfield}
{\it Let $(X,\nu)$ be a $\sigma$-finite standard measure space, $\textbf{G}^x$ be a measurable field of locally compact quantum groups over $(X, \nu)$, ad defined in \ref{deffieldqg}, and $\int_X^{\oplus}\textbf{G}^xd\nu(x)$ be the measured quantum groupoid constructed in \ref{mqgfield}. Then :
\newline
(i) there exists a locally compact set $\tilde{X}$, and a positive Radon measure $\tilde{\nu}$ on $\tilde{X}$, such that $L^\infty(X,\nu)$ and $L^\infty(\tilde{X}, \tilde{\nu})$ are isomorphic, and such that this isomorphism sends $\nu$ on $\tilde{\nu}$. 
\newline
(ii) there exists a continuous field $(A^x)_{x\in\tilde{X}}$ of ${\bf C}^*$-algebras, and a continuous field of coassociative coproducts $\tilde{\Gamma}^x : A^x\rightarrow A^x\otimes^mA^x$;
\newline
(iii) there exists left-invariant (resp. right-invariant) weights $\tilde{\varphi}^x$ (resp. $\tilde{\psi}^x$), such that $(A^x, \tilde{\Gamma}^x, \tilde{\varphi}^x, \tilde{\psi}^x)$ is a locally compact quantum group $\tilde{\textbf{G}}^x$ (in the ${\bf C}^*$ sense). 
\newline
(iv) we have : $\int_X^{\oplus}\textbf{G}^xd\nu(x)=\int_{\tilde{X}}^{\oplus}\tilde{\textbf{G}}^xd\tilde{\nu}(x)$.}

\begin{proof} Let's apply \ref{bicentral} to \ref{mqgfield}. \end{proof}
\subsection{Example}
\label{SU}
As in (\cite{Bl1}, 7.1), let us consider the ${\bf C}^*$-algebra $A$ whose generators $\alpha$, $\gamma$ and $f$ verify :
\newline
(i) $f$ commutes with $\alpha$ and $\gamma$;
\newline
(ii) the spectrum of $f$ is $[0,1]$;
\newline
(iii) the matrix 
$\left(\begin{array}{cc}
\alpha&-f\gamma\\
\gamma & \alpha^*
\end{array}
\right)$ is unitary in $M_2(A)$. 
Then, using the sub ${\bf C}^*$-algebra generated by $f$, $A$ is a $C([0,1])$-algebra; let us consider now $A$ as a $C_0(]0,1])$-algebra.  Then, Blanchard had proved (\cite{Bl2} 7.1) that $A$ is a continuous field over $]0,1]$ of ${\bf C}^*$-algebras, and that, for all $q\in]0,1]$, we have $A^q=SU_q(2)$, where the $SU_q(2)$ are the compact quantum groups constructed by Woronowicz and $A^1=C(SU(2))$. 
\newline
Moreover, using the coproducts $\Gamma^q$ defined by Woronowicz as 
\[\Gamma^q(\alpha)=\alpha\otimes\alpha-q\gamma^*\otimes\gamma\]
\[\Gamma^q(\gamma)=\gamma\otimes\alpha+\alpha^*\otimes\gamma\]
and the (left and right-invariant) Haar state $\varphi^q$, which verifies : 
\newline
$\varphi^q(\alpha^k\gamma^{*m}\gamma^n)=0$, for all $k\geq 0$, and $m\not=n$, 
\newline
$\varphi^q(\alpha^{*|k|}\gamma^{*m}\gamma^n)=0$, for all $k<0$, and $m\not=n$, 
\newline
and $\varphi^q((\gamma^*\gamma)^m)=\frac{1-q^2}{1-q^{2m+2}}$, 
\newline
we obtain this way a continuous field of compact quantum groups (see \cite{Bl2}, 6.6 for a definition); this leads to put on $A$ a structure of ${\bf C}^*$ quantum groupoid (of compact type, in the sense of \cite{E2}, because $1\in A$). 
\newline
This structure is given by a coproduct $\Gamma$ which is $C_0(]0,1])$-linear from $A$ to $A\underset{C_0(]0,1])}{\otimes^m}A$, and given by :
\[\Gamma(\alpha)=\alpha\underset{C_0(]0,1])}{\otimes^m}\alpha-f\gamma^*\underset{C_0(]0,1])}{\otimes^m}\gamma\]
\[\Gamma(\gamma)=\gamma\underset{C_0(]0,1])}{\otimes^m}\alpha+\alpha^*\underset{C_0(]0,1])}{\otimes^m}\gamma\]
and by a conditional expectation $E$ from $A$ on $M(C_0(]0,1]))$ given by :
\newline
$E(\alpha^k\gamma^{*m}\gamma^n)=0$, for all $k\geq 0$, and $m\not=n$
\newline
$E(\alpha^{*|k|}\gamma^{*m}\gamma^n)=0$, for all $k<0$, and $m\not=n$
\newline
$E((\gamma^*\gamma)^m)$ is the bounded function $x\mapsto\frac{1-q^2}{1-q^{2m+2}}$. 
Then $E$ is both left and right-invariant with respect to $\Gamma$. This example give results at the level of ${\bf C}^*$-algebras, which are more precise than theorem \ref{mqgfield}. 

\subsection{Example}
\label{ax+b}
One can find in \cite{Bl2} another example of a continuous field of locally compact quantum group. Namely, in (\cite{Bl2}, 7.2), Blanchard constructs a ${\bf C}^*$-algebra $A$ which is a continuous field of ${\bf C}^*$-algebras over $X$, where $X$ is a compact included in $]0, 1]$, with $1\in X$. For any $q\in X$, $q\not=1$, we have $A^q=SU_q(2)$, and $A^1={\bf C}^*_r(G)$, where $G$ is the "$ax+b$" group. (\cite{Bl2}, 7.6). 
\newline
Moreover, he constructs a coproduct (denoted $\delta$) (\cite{Bl2} 7.7(c)), and "the system of Haar weights" $\Phi$ ([B2] 7.2.3), which bear left-invariant-like properties (end of remark after \cite{Bl2} 7.2.3). 
\newline
Finally, he constructs a unitary $U$ in $\mathcal L(\mathcal E_\Phi)$ (\cite{Bl2} 7.10), with which it is possible to construct a co-inverse $R$ of $(A, \delta)$, which leads to the fact that $\Phi\circ R$ is right-invariant. 
\newline
Clearly, the fact that we are here dealing with non-compact locally compact quantum groups made the results more problematic at the level of ${\bf C}^*$-algebra; at the level of von Neumann algebra, \ref{mqgfield} allow us to construct an example of measured quantum groupoid from these data. 

\subsection{Example}
\label{Emu}
Let us finish by quoting a last example given by Blanchard in (\cite{Bl2}, 7.4): for $X$ compact in $[1, \infty[$, with $1\in X$, he constructs a ${\bf C}^*$-algebra which is a continuous field over $X$ of ${\bf C}^*$ algebras, whose fibers, for $\mu\in X$, are $A^\mu=E_\mu(2)$, with a coproduct $\delta$ and a continuous field of weights $\Phi$, which is left-invariant. As in \ref{ax+b}, he then constructs a unitary $U$ on $\mathcal L(\mathcal E_\Phi)$, which will lead to a co-inverse, and, therefore, to a right-invariant ${\bf C}^*$-weight. 

\subsection{Example (\cite{L2}, 17.1)}
\label{sum}
Let us return to \ref{deffieldqg}; let $I$ be a (discrete) set, and, for all $i$ in $I$, let $\textbf{G}_i=(M_i, \Gamma_i, \varphi_i, \psi_i)$ be a locally compact quantum group; then the product $\Pi_i \textbf{G}_i$ is a measured field of locally compact quantum groups, and can be given a natural structure of measured quantum groupoid, described in (\cite{L2}, 17.1).

\section{Measured Quantum Groupoid with central basis ${\bf C}^2$.}
\label{Example}
We finish by studying the structure of measured quantum groupoids with central basis ${\bf C}^2$. This example appears in \cite{DC} as a Galois object linking two locally compact quantum groups. 

\subsection{Lemma}
\label{lemC2}
{\it Let $\alpha$ be a representation of ${\bf C}^2$ on a Hilbert space $H$; let $(e_1, e_2)$ be the canonical basis of ${\bf C}^2$, $\nu$ the faithful normal state on ${\bf C}^2$ defined by $\nu(e_1)=\nu(e_2)=1/2$. Then :
\newline
(i) all vectors in $H$ are bounded with respect to $(\alpha, \nu)$. For any $\xi$, $\eta$ in $H$, we have :
\[\langle\xi, \eta\rangle_{\alpha, \nu}=(\alpha(e_1)\xi|\eta)e_1+(\alpha(e_2)\xi|\eta)e_2\]
\newline
(ii) for any representation $\beta$ of ${\bf C}^2$ on $H$, the application which sends $\xi\underset{\nu}{_\beta\otimes_\alpha}\eta$ on the vector $[\beta(e_1)\otimes\alpha(e_1)+\beta(e_2)\otimes\alpha(e_2)](\xi\otimes\eta)$ extends to an isomorphism of the relative tensor product $H\underset{\nu}{_\beta\otimes_\alpha}H$ with the subspace of the Hilbert tensor product $H\otimes H$ which is the image of the projection $\beta(e_1)\otimes\alpha(e_1)+\beta(e_2)\otimes\alpha(e_2)$. 
\newline
(iii) let $M$ a von Neumann algebra on $H$, such that $\alpha({\bf C}^2)\subset M$ and $\beta({\bf C}^2)\subset M$; then, the isomorphism given in (i) sends $M'\underset{{\bf C}^2}{_\beta\otimes_\alpha}M'$ on the induced von Neumann algebra $(M'\otimes M')_{\beta(e_1)\otimes\alpha(e_1)+\beta(e_2)\otimes\alpha(e_2)}$ and $M\underset{{\bf C}^2}{_\beta*_\alpha}M$ on its commutant, which is the reduced von Neumann algebra $(M\otimes M)_{\beta(e_1)\otimes\alpha(e_1)+\beta(e_2)\otimes\alpha(e_2)}$.
}
\begin{proof}
For all $(\lambda, \mu)\in{\bf C}^2$, $\xi\in H$, we have :
\begin{align*}
\|\alpha(\lambda e_1+\mu e_2)\xi\|^2
&=|\lambda|^2\|\alpha(e_1)\xi\|^2+|\mu|^2\|\alpha(e_2)\xi\|^2\\
&\leq \|\xi\|^2(|\lambda|^2+|\mu|^2)\\
&=\|\xi\|^2\nu(|\lambda|^2e_1+|\mu|^2e_2)
\end{align*}
which proves that $\xi\in D(_\alpha H, \nu)$; it is straightforward then to finish the proof of (i). Then (ii) is a direct corollary of (i), and (iii) is a direct corollary of (ii).
\end{proof}

\subsection{Remark}
\label{rem}
This kind of result can be generalized to any representation of a finite dimensional ${\bf C}^*$-algebra (\cite{DC}, 5), which generalizes the results of \cite{Val4} about relative tensor product of finite-dimensional Hilbert spaces. 

\subsection{Proposition}
\label{propC2}
{\it Let us use the notations of \ref{lemC2}; let $({\bf C}^2, M, \alpha, \beta, \Gamma, T, T', \nu)$ be a measured quantum groupoid, with $\alpha({\bf C}^2)\subset Z(M)$. Let $\Phi=\nu\circ\alpha^{-1}\circ T$. Then :
\newline
(i) the fiber product $M\underset{N}{_\beta*_\alpha}M$ can be identified with the reduced von  Neumann algebra $(M\otimes M)_{\beta(e_1)\otimes\alpha(e_1)+\beta(e_2)\otimes\alpha(e_2)}$, and $\Gamma$ can be identified with an injective $*$-homomorphism from $M$ to $M\otimes M$, which satisfies :
\[\Gamma(1)=\beta(e_1)\otimes\alpha(e_1)+\beta(e_2)\otimes\alpha(e_2)\]
\[(\Gamma\otimes id)\Gamma=(id\otimes\Gamma)\Gamma\]
\[\Gamma(\alpha(e_i)\beta(e_j))=\alpha(e_i)\beta(e_1)\otimes\alpha(e_1)\beta(e_j)+ \alpha(e_i)\beta(e_2)\otimes\alpha(e_2)\beta(e_j)\]
If we write $M_{i,j}=M\alpha(e_i)\beta(e_j)$, we have :
\[\Gamma(M_{i,j})\subset (M_{i,1}\otimes M_{1,j})\oplus (M_{i,2}\otimes M_{2,j})\]
and $M_{1,1}\neq\{0\}$, $M_{2,2}\neq\{0\}$. 
\newline
(ii) the pseudo-multiplicative unitary $W$ can be identified with a partial isometry on the Hilbert tensor product $H_\Phi\otimes H_\Phi$ with initial support $[\alpha(e_1)\otimes\alpha(e_1)+\alpha(e_2)\otimes\alpha(e_2)]$, and final support $[\beta(e_1)\otimes\alpha(e_1)+\beta(e_2)\otimes\alpha(e_2)]$, satisfying the pentagonal equation, and the following intertwining relations, for all $n\in{\bf C}^2$ :
\[W(\alpha(n)\otimes 1)=(1\otimes\alpha(n))W=(\alpha(n)\otimes 1)W\]
\[W(1\otimes\alpha(n))=W(\beta(n)\otimes1)=(\beta(n)\otimes 1)W\]
\[W(1\otimes\beta(n))=(1\otimes\beta(n))W\]
(iii) there exists normal semi-finite faithful weights $\varphi_{i,j}$ on $M_{i, j}$, such that, for any $X$ in $\gM_T^+$, $X=x_{1,1}\oplus x_{1,2}\oplus x_{2,1}\oplus x_{2,2}$, with $x_{i,j}\in M_{i,j}^+$, $T(X)$ is the image under $\alpha$ of $(\varphi_{1,1}(x_{1,1})+\varphi_{1,2}(x_{1,2}))e_1+(\varphi_{2,1}(x_{2,1})+\varphi_{2,2}(x_{2,2}))e_2$.
\newline
(iv) there exists normal semi-finite faithful weights $\psi_{i,j}$ on $M_{i,j}$, such that, for any $Y$ in $\gM_{T'}^+$, $Y=y_{1,1}\oplus y_{1,2}\oplus y_{2,1}\oplus y_{2,2}$, with $y_{i,j}\in M_{i,j}^+$, $T'(Y)$ is the image inder $\beta$ of $(\psi_{1,1}(y_{1,1})+\psi_{2,1}(y_{2,1}))e_1+(\psi_{1,2}(y_{1,2})+\psi_{2,2}(y_{2,2}))e_2$. 
\newline
 (v) for any $x_{i,j}$ in $M_{i,j}$ and $k=1,2$, let us define 
 \[\Gamma_{i,j}^k(x_{i,j})=\Gamma(x_{i,j})[\alpha(e_i)\beta(e_k)\otimes\alpha(e_k)\beta(e_j)]\]
 which implies that $\Gamma_{i,j}=\sum_k\Gamma_{i,j}^k$. Then, we have, where the $\delta_{i,j}$ are the usual Kronecker symbols, for any $x_{i,j}\in\gM_{\varphi_{i,j}}$, and $y_{i,j}\in\gM_{\psi_{i,j}}$ :
 \begin{multline*}
 \delta_{i,1}\varphi_{i,j}(x_{i,j})(\alpha(e_1)\otimes 1)+\delta_{i,2}\varphi_{i,j}(x_{i,j})(\alpha(e_2)\otimes 1)\\
 =(id\otimes\varphi_{1,j})\Gamma_{i,j}^1(x_{i,j})+(id\otimes\varphi_{2,j})\Gamma_{i,j}^2(x_{i,j})
 \end{multline*}
 \begin{multline*}
 \delta_{1,j}\psi_{i,j}(y_{i,j})(1\otimes\beta(e_1))+\delta_{2,j}\psi_{i,j}(y_{i,j})(1\otimes\beta(e_2))\\
 =(\psi_{i,1}\otimes id)\Gamma_{i,j}^1(y_{i,j})+(\psi_{i,2}\otimes id)\Gamma_{i,j}^2(y_{i,j})
 \end{multline*}
 }
 \begin{proof}
 The beginning of (i) is just a corollary of \ref{lemC2}, using the fact that $\alpha({\bf C}^2)$ (and $\beta({\bf C}^2)$ by \ref{propcentral}) is included in $Z(M)$; the end of (i) is given also by this fact, using also the fact that the formulae obtained for $\Gamma(\alpha(e_1)\beta(e_1))$ and $\Gamma(\alpha(e_2)\beta(e_2))$ prove that $\alpha(e_1)\beta(e_1)\neq 0$ and $\alpha(e_2)\beta(e_2)\neq 0$. 
 \newline
 As $\alpha$ is central, we have $\hat{\beta}=\alpha$, then, the identification of $W$ with a partial isometry comes from the identification of the relative tensor Hilbert spaves made in \ref{lemC2}(ii); this identification gives as well that this partial isometry satisfies the pentagonal equation, the intertwining properties; finally, the fact that $\alpha$ and $\beta$ are central finish the proof of (ii). 
 \newline
 Results (iii) and (iv) are given by \ref{continuousfield}(vi); then result (v) is given by the left-invariance of $T$ (resp. the right-invariance of $T'$). \end{proof}
 
 \subsection{Theorem (\cite{DC}, 3.17)}
 \label{thC2}
 {\it With the notations of \ref{lemC2}, let $\gG=({\bf C}^2, M, \alpha, \beta, \Gamma, T, T', \nu)$ be a measured quantum groupoid, $R$ its co-inverse, with $\alpha({\bf C}^2)\subset Z(M)$; let us write $\widehat{\gG}=(N, \alpha, \alpha, \widehat{\Gamma}, \hat{T}, \widehat{R}\hat{T}\widehat{R}, \nu)$ its dual. Let's use the notations of \ref{propC2}; then :
\newline
(i) $\textbf{G}^1=(M_{1,1}, \Gamma_{1,1}^1, \varphi_{1,1}, \psi_{1,1})$ and $\textbf{G}^2=(M_{2,2}, \Gamma_{2,2}^2, \varphi_{2,2}, \psi_{2,2})$ are two locally compact quantum groups. The multiplicative unitary $W^1$ (resp. $W^2$) of $\textbf{G}^1$ (resp. $\textbf{G}^2$) is equal to the restriction of $W$ to $\alpha(e_1)\beta(e_1)\otimes\alpha(e_1)\beta(e_1)$ (resp. $\alpha(e_2)\beta(e_2)\otimes\alpha(e_2)\beta(e_2)$); the coinverse $R^1$ (resp. $R^2$)  of $\textbf{G}^1$ (resp. $\textbf{G}^2$) is equal to the restriction of $R$ to $M_{1,1}$ (resp. $M_{2,2}$). 
\newline
(ii) if $\alpha({\bf C}^2)\subset Z(\widehat{M})$, then $\beta=\alpha$ and $\gG=\textbf{G}^1\oplus\textbf{G}^2$ (i.e. $M_{1,2}=M_{2,1}=\{0\}$).
\newline
(iii) if $\alpha(e_1)\notin Z(\widehat{M})$, then $M_{1,2}\neq \{0\}$, $M_{2,1}=R(M_{1,2})\neq \{0\}$, $\Gamma_{1,2}^2 : M_{1,2}\rightarrow M_{1,2}\otimes M_{2,2}$ is a right action of $\textbf{G}^2$ on $M_{1,2}$, $\Gamma_{1,2}^1: M_{1,2}\rightarrow M_{1,1}\otimes M_{1,2}$ is a left action of $\textbf{G}^1$ on $M_{1,2}$, $\Gamma_{2,1}^1: M_{2,1}\rightarrow M_{2,1}\otimes M_{1,1}$ verify $\Gamma_{2,1}^1=\varsigma(R^1\otimes R)\Gamma_{1,2}^1R$, and $\Gamma_{2,1}^2: M_{2,1}\rightarrow M_{2,2}\otimes M_{2,1}$ verify $\Gamma_{2,1}^2=\varsigma(R\otimes R^2)\Gamma_{1,2}^2R$. Moreover, these actions are ergodic and integrable. }

\begin{proof}
Using \ref{propC2}, we get that $M_{1,1}\neq 0$, and that :
\[\Gamma(\alpha(e_1)\beta(e_1))[\alpha(e_1)\beta(e_1)\otimes \alpha(e_1)\beta(e_1)]=\alpha(e_1)\beta(e_1)\otimes \alpha(e_1)\beta(e_1)\]
from which we get that $\Gamma_{1,1}^1(1)=1$, when considered from $M_{1,1}$ into $M_{1,1}\otimes M_{1,1}$; by restriction, the coproduct property is straightforward from \ref{propC2}(i); using \ref{propC2}(iv), we get :
\[\alpha(e_1)(id\otimes\varphi_{1,1})\Gamma_{1,1}^1(x_{1,1})=\alpha(e_1)\varphi_{1,1}(x_{1,1})\]
and :
\[\beta(e_1)(\psi_{1,1}\otimes id)\Gamma_{1,1}^1(y_{1,1})=\beta(e_1)\psi_{1,1}(y_{1,1})\]
which proves that $\varphi_{1,1}$ is a left-invariant weight, and $\psi_{1,1}$ is a right-invariant weight, and proves that $\textbf{G}^1$ is a locally compact quantum group, in the sense of \cite{KV2}. Then, the result about the multiplicative unitary of $\textbf{G}^1$ is a straightforward calculation, and, then, by polar decomposition of the antipode, one gets the result about the coinverse of $\textbf{G}^1$. The proof for $\textbf{G}^2$ is identical, which finishes the proof of (i). 
\newline
Result (ii) is given by \ref{bicentral}; conversely, if $M_{1,2}=\{0\}$, as $M_{2,1}=R(M_{1,2})$, we have also $M_{2,1}=\{0\}$, and $M=M_{1,1}\oplus M_{2,2}$, and, using \ref{mqgfield}(iv), we get that $\alpha=\beta$, and, therefore, that $\alpha({\bf C}^2)\subset Z(\widehat{M})$. So, we get that, if $\alpha(e_1)\notin Z(\widehat{M})$, we have $M_{1,2}\neq\{0\}$, and $M_{2,1}\neq\{0\}$. Then, by restriction of the coproduct property of $\Gamma$, we obtain that $\Gamma_{1,2}^2$ is a right-action of $\textbf{G}^2$ on $M_{1,2}$, and that $\Gamma_{2,1}^1$ is a left-action of $\textbf{G}^1$ on $M_{2,1}$ (in the sense of (\cite{V},1.1); the properties of $\Gamma_{2,1}^1$ and $\Gamma_{2,1}^2$ come from the formula linking $\Gamma$ and $R$ (and the fact that $R_{|M_{1,1}}=R^1$ and $R_{|M_{2,2}}=R^2$ obtained in (i)). Moreover, using \ref{propC2}(v), one gets that, for any $x_{1,2}\in\gM_{\varphi_{1,2}}$, we have :
\[\varphi_{1,2}(x_{1,2})=(id\otimes\varphi_{2,2})\Gamma_{1,2}^2(x_{1,2})\]
But the right-hand formula is the canonical operator-valued weight $T_{\Gamma_{1,2}^2}$ (\cite{V}, 1.3) from $M_{1,2}$ on the invariants $M_{1,2}^{\Gamma_{1,2}^2}$; so we get both that this algebra $M_{1,2}^{\Gamma_{1,2}^2}$is equal to ${\bf C}$ (which means that $\Gamma_{1,2}^2$ is ergodic), and that this operator-valued weight is semi-finite (which means that $\Gamma_{1,2}^2$ is integrable). The proof for $\Gamma_{2,1}^1$ is identical. 
\end{proof}

\subsection{Remark}
\label{rem2}
In \cite{DC} is given a very interesting interpretation of these actions, and of the link between $\textbf{G}^1$ and $\textbf{G}^2$ which occur in that situation, in term of Morita-Rieffel equivalence. Let's have a look at what happens when $\gG$ is abelian (resp. symmetric) : 
\newline
If $\gG$ is abelian, by \ref{ramsay}, we have $\gG=\gG(\mathcal G)$, where $\mathcal G$ is a locally compact groupoid, with a two-points basis. Then, if $X=\{x\in\mathcal G, s(g)=1, r(g)=2\}$ is not empty, it is clear that, in the construction given in \ref{thC2}, we obtain two locally compact groups which are isomorphic, and acts left and right on $X$; if $X$ is empty, we obtain that $\mathcal G$ is the disjoint union of the two locally compact groups $G^1$ and $G^2$ (\ref{sum})).
\newline
If $\gG$ is symmetric, then $\alpha({\bf C}^2)\subset Z(\widehat{M})$, and $\gG=\widehat{G^1}\oplus\widehat{G^2}$, where $G^i$ are locally compact groups, and $\widehat{G^i}$ their duals as symmetric locally compact quantum groups. 
\newline
So, we see that this construction, which is completely trivial in the case of groupoids, give very rich information in the case of quantum groupoids. 


\end{document}